\documentclass[11pt,leqno,twoside]{amsart}

\usepackage{amsmath, amsfonts, amssymb, amsthm, bm}
\usepackage{iftex}
\ifPDFTeX
  \usepackage[utf8]{inputenc} 
  \usepackage[T1]{fontenc} 
  \usepackage{lmodern} 
\else
  \usepackage{fontspec} 
\fi
\usepackage{accents}

\newcommand{\thickbar}[1]{\accentset{\rule{0.45em}{0.8pt}}{#1}}

\usepackage{booktabs}
\newtheorem{thm}{Theorem}[section] 
\newtheorem{rem}[thm]{Remark} 
\usepackage[scale=0.8]{geometry}

\usepackage{xcolor}
\usepackage{morefloats}
\usepackage{mathtools} 
\usepackage{enumitem}
\usepackage[english]{babel}

\usepackage{caption}
\usepackage{subcaption}
\usepackage{svg}
\usepackage{mathrsfs}

\numberwithin{equation}{section}

\usepackage{tikz}
\usepackage{multirow}

\usepackage{algorithmic}
\usepackage{algorithm}

\usepackage{todonotes}
\usepackage{amsmath, amsthm}
\usepackage{amsfonts, amssymb}
\usepackage{bm}
\usepackage[colorlinks=true, allcolors=blue]{hyperref}

	\newtheorem{theorem}{Theorem}[section]
	\newtheorem*{theorem*}{Theorem}
	\newtheorem{lemma}{Lemma}[section]
	
	\newtheorem*{proposition*}{Proposition}

\theoremstyle{remark}

\theoremstyle{definition}

\newcommand{\RR}{\mathbb{R}}

\newcommand{\Q}{\mathcal{Q}}

\title[Kinetic-Fluid Domain Decomposition for Multitemperature Gas Mixtures]{Realizability-Based Kinetic-Fluid Domain Decomposition for Multitemperature Gas Mixtures on Unstructured Meshes}

\author[D. Caparello]{Domenico Caparello}
\address{Domenico Caparello \\
Université Côte d’Azur, CNRS, LJAD, Parc Valrose, F-06108 Nice, France\\
Department of Mathematics and Computer Science, University of Ferrara \\
}
\email{Domenico.CAPARELLO@univ-cotedazur.fr}
\email{domenico.caparello@unife.it}

\author[L. Pareschi]{Lorenzo Pareschi}
\address{Lorenzo Pareschi \\
Maxwell Institute for Mathematical Sciences and Department of Mathematics,
Heriot-Watt University \\
Department of Mathematics and Computer Science, University of Ferrara \\
}
\email{L.Pareschi@hw.ac.uk}

\author[T. Rey]{Thomas Rey}
\address{Thomas Rey \\
Université Côte d’Azur, CNRS, LJAD, Parc Valrose, F-06108 Nice, France
}
\email{Thomas.REY@univ-cotedazur.fr}

\keywords{Gas mixtures, kinetic--fluid coupling, moment realizability, dynamic domain decomposition, asymptotic-preserving schemes, IMEX Runge--Kutta methods, unstructured meshes, finite volume methods}
	
\subjclass[2020]{Primary: 76P05,
  82C40;
  Secondary: 65M08,
  65M55 
}

\begin{document}
\newcommand{\figrow}[3]{%
\begin{subfigure}{0.32\textwidth}\includegraphics[width=\linewidth]{#1}\end{subfigure}\hfill
\begin{subfigure}{0.32\textwidth}\includegraphics[width=\linewidth]{#2}\end{subfigure}\hfill
\begin{subfigure}{0.32\textwidth}\includegraphics[width=\linewidth]{#3}\end{subfigure}\par\smallskip}

\begin{abstract}
We present a dynamic domain decomposition method for gas mixtures described by a multi-species BGK model in the regime where collisions between particles of the same species are dominant. In this regime the hydrodynamic limit is a system of Euler equations in which every species keeps its own velocity and temperature, the species being coupled through momentum and energy exchange terms. To decide locally whether this multi-velocity, multi-temperature Euler description is accurate, we extend the moment realizability criterion of Levermore, Morokoff and Nadiga \cite{lev-mor-nad-1998}: for each species, a realizability matrix is derived from the Chapman--Enskog expansion of the kinetic model, and its zeroth- and first-order approximations are compared. The indicator involves only macroscopic quantities and is evaluated species by species, so that in the same cell one species may be described by the Euler equations and another one by the kinetic equation. The kinetic solver is therefore used only where, and for the species for which, the indicator detects non-equilibrium effects, such as across shocks or near obstacles in rarefied flows.
Both solvers are discretized on unstructured triangular meshes with second-order accuracy in space and time. The kinetic equation is integrated by an asymptotic-preserving IMEX Runge--Kutta scheme, whose limit for vanishing Knudsen number is a consistent discretization of the multi-velocity, multi-temperature Euler system. In the numerical tests the hybrid solution stays close to the kinetic one, while the computational time is reduced by a factor between about three and almost five.
\end{abstract}

\maketitle

\section{Introduction}\label{sec1}
\setcounter{equation}{0}
Rarefied gas mixtures appear in many applications, from high-altitude atmospheric flows and nuclear engineering to evaporation--condensation problems \cite{jin-qin-2013} and plasma physics, where electrons and ions are modelled as different interacting species \cite{bot-cow-sch-2024}. In several of these problems the flow is in a transitional regime, as during the re-entry of space vehicles \cite{dill-cor-2008} or in micro-electro-mechanical systems (MEMS). At high altitude non-equilibrium effects become important and macroscopic models lose accuracy: the Navier--Stokes equations, for instance, do not give correct shock profiles in the hypersonic regime when rarefaction effects are significant, whereas a kinetic description remains applicable. Near strong shocks or in kinetic boundary layers fluid models are therefore not reliable, and for mixtures a multi-species kinetic description is needed.

The Boltzmann equation for gas mixtures has been studied in several works, see e.g. \cite{mar-2016, bri-dau-2016}. In \cite{rey-ten-2024} the analysis was extended to granular mixtures, proving well-posedness in the spatially homogeneous case and studying the long-time behaviour. All these models share the same structure: a linear transport operator and one or more stiff, nonlinear collision operators, acting on distribution functions defined in a high-dimensional phase space \cite{dim-par-2014}.

The main difficulty in the numerical solution of kinetic equations is the dimension of the phase space, up to seven including time, which makes kinetic solvers much more expensive than fluid solvers. Moreover, kinetic models contain small scales, such as the mean free path. Asymptotic expansions in these small parameters lead to hydrodynamic equations, such as the Euler equations, which describe the average behaviour of the gas. Solving the kinetic equation everywhere in a near-continuum regime is very expensive, and wherever the hydrodynamic model is valid it is much cheaper to solve it instead.

Hybrid methods, or dynamic domain decomposition methods \cite{dim-par-2014, bur-boy-2009, fil-rey-2014, tiw-1998, cap-par-rey-2025, cap-ten-2025, cap-par-rey-2026}, are based on this observation: at every time step the spatial domain is divided into regions where different models are solved, according to the local degree of non-equilibrium. Compared with a purely hydrodynamic computation the solution is more accurate, and compared with a purely kinetic one it is cheaper.
Other kinetic--fluid approaches for gas mixtures include the particle-based micro--macro method of \cite{cre-kli-pir-2020}. In \cite{dat-shu-2023} a multi-species kinetic model is coupled with a multi-fluid model, and in \cite{chu-sag-mur-haa-2024} a multi-species Vlasov--BGK model is coupled with its single-velocity, single-temperature hydrodynamic limit through buffer regions; in both works the kinetic region is fixed in time. Here, instead, the partition is dynamic: it is updated at every time step with a realizability-based indicator, so that the kinetic and fluid regions of each species evolve in time with the flow.

In collision-dominated regions the collisions drive the gas quickly towards local equilibrium, and the collision operator is stiff. Explicit schemes then require time steps of the order of the relaxation time. Asymptotic-preserving (AP) schemes \cite{cor-per-1991, jin-1995, caf-shi-rus-1995, jin-par-tos-1998, jin-1999, jin-par-2001, klar-1998, gos-tos-2002, diaz-chen-2015, fil-jin-2011, fil-jin-2010, jin-qin-2013, hu-wan-2015, bos-cho-rus-yun-2021, bos-par-rus-2024} remove this restriction: they are stable for time steps independent of the Knudsen number and reduce, in the limit, to a consistent discretization of the hydrodynamic equations \cite{shi-2010}. We use this property for the kinetic solver.

In this work we consider mixtures in which collisions between particles of the same species are much more frequent than collisions between different species. This happens, for instance, in mixtures of light and heavy gases, where the energy exchange between species is slow \cite{gal-mak-1994, bis-gro-mar2021}. Each species then relaxes to its own Maxwellian, and the hydrodynamic limit is a system of Euler equations with one velocity and one temperature per species \cite{bis-gro-mar2021}. We propose a hybrid method that couples this multi-velocity, multi-temperature Euler system with the multi-species BGK model and switches between the two automatically. To our knowledge, this is the first realizability-based dynamic spatial domain decomposition method for a kinetic mixture model whose hydrodynamic limit has a velocity and a temperature for each species. Since the fluid model keeps these quantities separate, the choice between fluid and kinetic description can be made species by species.

The switching criterion is based on a moment realizability matrix that we derive for this regime. The computations are carried out for the BGK model of \cite{bob-bis-gro-spi-pot-2018, bis-gro-mar2021}. The collision model enters the indicator only through the Chapman--Enskog correction terms, so the same construction can be used with the full Boltzmann operator for mixtures once these corrections are known. Since the indicator is based on the Chapman--Enskog expansion, it can also be extended to more than two levels, as done for a single gas in \cite{cap-par-rey-2025, cap-par-rey-2026}. In \cite{cap-par-rey-2025} the hierarchy consists of the Euler equations, the ellipsoidal statistical BGK (ES-BGK) model \cite{and-let-per-per-2000} and the full Boltzmann equation on Cartesian meshes. In \cite{cap-par-rey-2026}, on unstructured meshes, the intermediate level is given by the Navier--Stokes equations, which correspond to the ES-BGK model at first order in the Chapman--Enskog expansion: on unstructured meshes the cost is dominated by transport, and the ES-BGK model does not give a clear computational advantage.

The coupling is non-intrusive: the two solvers only exchange moments and Maxwellians, and they only need to use the same Runge--Kutta stages. The kinetic solver is used only in the cells, and for the species, selected by the indicator. We implement the method on unstructured triangular meshes with a second-order finite volume reconstruction, so that problems with non-trivial geometries can be treated.

The opposite situation, in which all collisions occur on the same time scale and the mixture relaxes to a single velocity and a single temperature, leads to a different hydrodynamic limit and requires a different indicator. It is considered in the companion paper \cite{PaperB}.

The paper is organized as follows. Section~\ref{secModel} presents the kinetic and macroscopic models. In Section~\ref{secIndicators} we derive the moment realizability matrices and the regime indicators. Section~\ref{secNumSch} describes the numerical schemes, the AP property and the coupling between kinetic and fluid solvers. Numerical tests are presented in Section~\ref{secSimulations}, and the computational cost is discussed in Section~\ref{secNumEff}. Section~\ref{secCon} contains some conclusions.

\section{Multi-temperature, multi-velocity BGK model for gas mixtures}\label{secModel}
We recall the kinetic model considered in this work, following \cite{bis-gro-mar2021, bob-bis-gro-spi-pot-2018}. We consider a mixture of $N$ monatomic gases whose particles undergo elastic collisions. Species $i$ is described by its distribution function $f_i(t,x,v)$, where $x\in\RR^3$ is the position, $v\in\RR^3$ the velocity and $t\ge0$ the time. The evolution is given by the Boltzmann equations
\begin{equation}\label{eq1}
\begin{cases}
\displaystyle\frac{\partial f_i}{\partial t} + v\cdot\nabla_x f_i = \mathcal{Q}_i[f],\qquad \mathcal{Q}_i[f]=\sum_{k=1}^N\mathcal{Q}_{i,k}[f_i,f_k],\\[6pt]
f_i(0,x,v)=f_{i,0}(x,v),
\end{cases}
\end{equation}
where $\mathcal{Q}_i$ is the collision operator of species $i$ and the binary operators read
\begin{equation}\label{eqQBoltzMulti}
\Q_{i,k}[f_i,f_k] = \int_{\RR^3}\int_{\mathbb S^2} g_{i,k}(|y|,\hat y\cdot\Omega)\,\big[f_i(v')f_k(w')-f_i(v)f_k(w)\big]\,d\Omega\,dw .
\end{equation}
Here $y=v-w$ is the relative velocity, $v'$ and $w'$ are the pre-collisional velocities, and the collision kernel $g_{i,k}$ depends on $|y|$ and on $\mu=\hat y\cdot\Omega$, with $\hat y=y/|y|$. The number density, mean velocity and temperature of species $i$, whose particles have mass $m_i$, are
\begin{equation}\label{eqMoments}
n_i=\int_{\RR^3}f_i\,dv,\qquad u_i=\frac{1}{n_i}\int_{\RR^3}v\,f_i\,dv,\qquad T_i=\frac{m_i}{3n_i}\int_{\RR^3}|v-u_i|^2f_i\,dv .
\end{equation}
The collision equilibria of \eqref{eq1} are local Maxwellians with a common velocity $u$ and a common temperature $T$,
\begin{equation}\label{eqFeq}
f_i^{eq}=n_i\left(\frac{m_i}{2\pi T}\right)^{3/2}\exp\left(-\frac{m_i}{2T}|v-u|^2\right),
\end{equation}
where
\begin{equation}\label{eqGlobal}
\begin{aligned}
n&=\sum_{i=1}^N n_i,\qquad \rho=\sum_{i=1}^N m_in_i,\qquad
u=\frac{1}{\rho}\sum_{i=1}^N m_in_iu_i,\\
T&=\frac{1}{3n}\sum_{i=1}^N m_i\int_{\RR^3}|v-u|^2f_i\,dv\\
 &=\frac{1}n\sum_{i=1}^N n_iT_i+\frac{1}{3n}\sum_{i=1}^N m_in_i|u_i-u|^2.
\end{aligned}
\end{equation}
In \cite{bis-gro-mar2021} the authors consider the BGK model of \cite{bob-bis-gro-spi-pot-2018}, in which the collision operators are replaced by relaxation operators,
\begin{equation}\label{eq8}
\mathcal{Q}_i[f]=\sum_{k=1}^N\nu_{ik}\,(n_iM_{ik}-f_i),
\end{equation}
where $\nu_{ik}$ are collision frequencies, independent of $v$, and $M_{ik}$ are Maxwellians with mass $m_i$, velocity $u_{ik}$ and temperature $T_{ik}$,
\begin{equation}
M_{ik}=\left(\frac{m_i}{2\pi T_{ik}}\right)^{3/2}\exp\left(-\frac{m_i}{2T_{ik}}|v-u_{ik}|^2\right).
\end{equation}
The parameters $u_{ik}$, $T_{ik}$ are fixed by requiring that the exchange rates of momentum and energy coincide with those of the Boltzmann operator for Maxwellian interactions; this guarantees the conservation of total momentum and total energy \cite{bis-gro-mar2021}. One obtains
\begin{equation}\label{eqUT}
u_{ik}=(1-a_{ik})u_i+a_{ik}u_k,\qquad T_{ik}=(1-b_{ik})T_i+b_{ik}T_k+\gamma_{ik}|u_i-u_k|^2,
\end{equation}
with
\begin{equation}\label{eq_aik_bik_gammaik}
a_{ik}=\frac{\lambda_{ik}n_km_k}{\nu_{ik}(m_i+m_k)},\qquad b_{ik}=\frac{2a_{ik}m_i}{m_i+m_k},\qquad \gamma_{ik}=\frac{m_ia_{ik}}{3}\left(\frac{2m_k}{m_i+m_k}-a_{ik}\right),
\end{equation}
and $\lambda_{ik}=2\pi\int_{-1}^1(1-\mu)\,g_{i,k}(\mu)\,d\mu$; for Maxwellian interactions $g_{i,k}$ does not depend on $|y|$ and $\lambda_{ik}$ is a constant. We assume $\lambda_{ik}=\lambda_{ki}$, which gives pairwise conservation of momentum and energy. Since $u_{ii}=u_i$ and $T_{ii}=T_i$, collisions between particles of the same species are described by the usual BGK operator for a single gas.

We are interested in inert mixtures in which collisions between particles of the same species are much more frequent than collisions between particles of different species, as in $\varepsilon$-mixtures of heavy and light gases, where energy exchange between species is slower \cite{gal-mak-1994, bis-gro-mar2021}. The rescaled model reads \cite{bis-gro-mar2021}
\begin{equation}\label{rescaledBGK}
\frac{\partial f_i}{\partial t}+v\cdot\nabla_xf_i=\frac{1}{\varepsilon}\,\nu_{ii}(n_iM_{ii}-f_i)+\sum_{\substack{k=1\\k\ne i}}^N\nu_{ik}(n_iM_{ik}-f_i),\qquad i=1,\dots,N,
\end{equation}
where $\varepsilon$ is the Knudsen number, the ratio between the mean free path and a macroscopic length. For $\varepsilon\ll1$ the gas is in the fluid regime, while for $\varepsilon\gtrsim1$ it is rarefied.

\section{Computation of the indicators}\label{secIndicators}
In this section we construct the indicator used to decide, for each species and in each cell, whether the multi-temperature, multi-velocity Euler equations give an adequate description. We follow the approach based on moment realizability, introduced in \cite{lev-mor-nad-1998} for a single gas and later used as a breakdown criterion in \cite{fil-rey-2014, cap-ten-2025, fil-xio-2018, tiw-2000, xio-qiu-2017}. Departure from equilibrium is estimated by comparing successive approximations of a moment matrix built from macroscopic quantities; loss of positive semi-definiteness is a sufficient warning of non-realizability. Evaluating this indicator does not require reconstructing the distribution function \cite{lev-mor-nad-1998, fil-rey-2014}. To our knowledge, realizability criteria have not been used before for a domain decomposition in the multi-temperature, multi-velocity limit of a mixture. The matrix is obtained from the Chapman--Enskog expansion of \eqref{rescaledBGK}, which we recall first.

\subsection{Chapman--Enskog expansion}
We follow \cite{bis-gro-mar2021}. We write $f_i=f_i^{(0)}+\varepsilon f_i^{(1)}+\mathcal O(\varepsilon^2)$ and expand in the same way
\[
M_{ik}=M_{ik}^{(0)}+\varepsilon M_{ik}^{(1)},\quad \nu_{ik}=\nu_{ik}^{(0)}+\varepsilon\nu_{ik}^{(1)},\quad u_{ik}=u_{ik}^{(0)}+\varepsilon u_{ik}^{(1)},\quad T_{ik}=T_{ik}^{(0)}+\varepsilon T_{ik}^{(1)} .
\]
In this regime the moments of every species are hydrodynamic variables and are not modified by the expansion:
\begin{equation}
n_i^{(0)}=n_i,\quad u_i^{(0)}=u_i,\quad T_i^{(0)}=T_i,\qquad n_i^{(1)}=u_i^{(1)}=T_i^{(1)}=0 .
\end{equation}
Moreover $M_{ii}=M_i$, with
\begin{equation}
M_i=\left(\frac{m_i}{2\pi T_i}\right)^{3/2}\exp\left(-\frac{m_i}{2T_i}|v-u_i|^2\right),
\end{equation}
so that $M_{ii}^{(0)}=M_i$ and $M_{ii}^{(1)}=0$. At zeroth order one finds \cite{bis-gro-mar2021}
\begin{equation}\label{eq4}
f_i^{(0)}=n_iM_i ,
\end{equation}
and the parameters of the mixed Maxwellians expand as
\begin{equation}\label{eq6}
\begin{aligned}
u_{ik}^{(0)}&=(1-a_{ik}^{(0)})u_i+a_{ik}^{(0)}u_k, & u_{ik}^{(1)}&=a_{ik}^{(1)}(u_k-u_i),\\
T_{ik}^{(0)}&=(1-b_{ik}^{(0)})T_i+b_{ik}^{(0)}T_k+\gamma_{ik}^{(0)}|u_i-u_k|^2, & T_{ik}^{(1)}&=b_{ik}^{(1)}(T_k-T_i)+\gamma_{ik}^{(1)}|u_i-u_k|^2,
\end{aligned}
\end{equation}
with
\begin{equation}
\begin{aligned}
a_{ik}^{(0)}&=\frac{\lambda_{ik}n_km_k}{\nu_{ik}^{(0)}(m_i+m_k)}, & a_{ik}^{(1)}&=-\frac{\nu_{ik}^{(1)}}{\nu_{ik}^{(0)}}\,a_{ik}^{(0)},\\
b_{ik}^{(0)}&=\frac{2a_{ik}^{(0)}m_i}{m_i+m_k}, & b_{ik}^{(1)}&=-\frac{\nu_{ik}^{(1)}}{\nu_{ik}^{(0)}}\,b_{ik}^{(0)},\\
\gamma_{ik}^{(0)}&=\frac{m_ia_{ik}^{(0)}}{3}\left(\frac{2m_k}{m_i+m_k}-a_{ik}^{(0)}\right), & \gamma_{ik}^{(1)}&=-\frac{\nu_{ik}^{(1)}}{\nu_{ik}^{(0)}}\left(\gamma_{ik}^{(0)}-\frac{1}{3}m_i\big(a_{ik}^{(0)}\big)^2\right).
\end{aligned}
\end{equation}
The hydrodynamic equations for $n_i$, $u_i$, $T_i$ follow from the weak form of \eqref{rescaledBGK} \cite{bis-gro-mar2021}. Inserting the truncated expansion and integrating against $\phi_i(v)\in\{1,\,m_iv,\,\frac{1}{2}m_i|v|^2\}$ gives, neglecting $\mathcal O(\varepsilon^2)$ terms,
\begin{equation}\label{weakForm}
\begin{aligned}
&\frac{\partial}{\partial t}\int_{\RR^3}\phi_i\,(f_i^{(0)}+\varepsilon f_i^{(1)})\,dv+\nabla_x\cdot\int_{\RR^3}\phi_i\,v\,(f_i^{(0)}+\varepsilon f_i^{(1)})\,dv
=-(\nu_{ii}^{(0)}+\varepsilon\nu_{ii}^{(1)})\int_{\RR^3}\phi_i\,f_i^{(1)}\,dv\\
&\qquad+\sum_{k\ne i}(\nu_{ik}^{(0)}+\varepsilon\nu_{ik}^{(1)})\int_{\RR^3}\phi_i\,\big(n_iM_{ik}^{(0)}+\varepsilon n_iM_{ik}^{(1)}-f_i^{(0)}-\varepsilon f_i^{(1)}\big)\,dv .
\end{aligned}
\end{equation}
Retaining only $f_i^{(0)}$ gives the multi-temperature, multi-velocity Euler equations; including $f_i^{(1)}$ gives the corresponding Navier--Stokes equations, with viscosity and heat conduction.

\subsection{Euler equations}
At zeroth order, \eqref{weakForm} gives the Euler equations for a mixture with one velocity and one temperature per species \cite{bis-gro-mar2021}:
\begin{equation}\label{eq7}
\begin{cases}
\displaystyle\frac{\partial n_i}{\partial t}+\nabla_x\cdot(n_iu_i)=0,\\[8pt]
\displaystyle\frac{\partial}{\partial t}(\rho_iu_i)+\nabla_x\cdot(\rho_iu_i\otimes u_i)+\nabla_x(n_iT_i)=\sum_{k\ne i}R_{ik},\\[8pt]
\displaystyle\frac{\partial}{\partial t}\Big(\frac{1}{2}\rho_i|u_i|^2+\frac{3}{2}n_iT_i\Big)+\nabla_x\cdot\Big[\Big(\frac{1}{2}\rho_i|u_i|^2+\frac{5}{2}n_iT_i\Big)u_i\Big]=\sum_{k\ne i}S_{ik},
\end{cases}
\end{equation}
with $\rho_i=m_in_i$. The exchange terms
\[
R_{ik}=\nu_{ik}^{(0)}\int_{\RR^3}m_iv\,(n_iM_{ik}^{(0)}-f_i^{(0)})\,dv,\qquad S_{ik}=\nu_{ik}^{(0)}\int_{\RR^3}\frac{1}{2}m_i|v|^2(n_iM_{ik}^{(0)}-f_i^{(0)})\,dv
\]
can be computed explicitly \cite{bis-gro-mar2021}:
\begin{equation}\label{eqRS}
R_{ik}=\lambda_{ik}m_{ik}n_in_k(u_k-u_i),\qquad S_{ik}=\lambda_{ik}\frac{m_{ik}}{m_i+m_k}n_in_k\big[(m_iu_i+m_ku_k)\cdot(u_k-u_i)+3(T_k-T_i)\big],
\end{equation}
where $m_{ik}=m_im_k/(m_i+m_k)$ is the reduced mass.

\subsection{Navier--Stokes equations}
To obtain the Navier--Stokes level we need the first-order correction $f_i^{(1)}$. Following \cite{bis-gro-mar2021}, with $c_i=v-u_i$,
\begin{equation}\label{eq100}
\begin{aligned}
f_i^{(1)}=\frac{1}{\nu_{ii}^{(0)}}\Bigg\{&n_i\sum_{k\ne i}\nu_{ik}^{(0)}\big(M_{ik}^{(0)}-M_i\big)-n_iM_i\Bigg[\frac{1}{n_iT_i}\,c_i\cdot\sum_{k\ne i}R_{ik}
+m_i\Big(\frac{c_i\otimes c_i}{T_i}-\frac{1}{3}\frac{|c_i|^2}{T_i}I\Big):\nabla_xu_i\\
&+\frac{2}{3}\frac{1}{n_iT_i}\Big(\frac{m_i}{2}\frac{|c_i|^2}{T_i}-\frac{3}{2}\Big)\Big(\sum_{k\ne i}S_{ik}-u_i\cdot\sum_{k\ne i}R_{ik}\Big)
+\Big(\frac{m_i}{2}\frac{|c_i|^2}{T_i}-\frac{5}{2}\Big)\frac{c_i}{T_i^{1/2}}\cdot\frac{\nabla_xT_i}{T_i^{1/2}}\Bigg]\Bigg\},
\end{aligned}
\end{equation}
where $I$ is the identity matrix and $u_{ik}^{(0)}$, $T_{ik}^{(0)}$ are given by \eqref{eq6}. Inserting $f_i^{(0)}+\varepsilon f_i^{(1)}$ into \eqref{weakForm} and neglecting $\mathcal O(\varepsilon^2)$ terms gives \cite{bis-gro-mar2021}
\begin{equation}
\begin{cases}
\displaystyle\frac{\partial n_i}{\partial t}+\nabla_x\cdot(n_iu_i)=0,\\[8pt]
\displaystyle\frac{\partial}{\partial t}(\rho_iu_i)+\nabla_x\cdot(\rho_iu_i\otimes u_i)+\nabla_x(n_iT_i)+\varepsilon\nabla_x\cdot P_i^{(1)}=\sum_{k\ne i}R_{ik},\\[8pt]
\displaystyle\frac{\partial}{\partial t}\Big(\frac{1}{2}\rho_i|u_i|^2+\frac{3}{2}n_iT_i\Big)+\nabla_x\cdot\Big[\Big(\frac{1}{2}\rho_i|u_i|^2+\frac{5}{2}n_iT_i\Big)u_i\Big]+\varepsilon\nabla_x\cdot\big(P_i^{(1)}u_i\big)+\varepsilon\nabla_x\cdot q_i^{(1)}=\sum_{k\ne i}S_{ik},
\end{cases}
\end{equation}
with
\[
P_i^{(1)}=m_i\int_{\RR^3}c_i\otimes c_i\,f_i^{(1)}\,dv,\qquad q_i^{(1)}=\frac{1}{2}m_i\int_{\RR^3}c_i|c_i|^2f_i^{(1)}\,dv .
\]

\subsection{Regime indicators}
We now define the regime indicators, extending to \eqref{rescaledBGK} the approach of \cite{lev-mor-nad-1998}, developed further in \cite{fil-rey-2014, fil-xio-2018, tiw-2000, li-son-wan-2021, xio-qiu-2017, cap-par-rey-2025, cap-ten-2025, cap-par-rey-2026}. For each species we introduce the polynomials
\begin{equation}\label{soninePol}
A_i=\frac{m_i}{n_i}\Big[\frac{c_i\otimes c_i}{T_i}-\frac{1}{3}\frac{|c_i|^2}{T_i}I\Big],\qquad
B_i=\frac{\sqrt{m_i}}{n_i}\Big[\frac{m_i}{2}\frac{|c_i|^2}{T_i}-\frac{5}{2}\Big]\frac{c_i}{T_i^{1/2}},\qquad
C_i=\frac{1}{n_i}\Big[\frac{m_i}{2}\frac{|c_i|^2}{T_i}-\frac{5}{2}\Big]^2 .
\end{equation}
With this notation \eqref{eq100} becomes
\begin{equation}\label{eq3}
\begin{aligned}
f_i^{(1)}=\frac{1}{\nu_{ii}^{(0)}}\Bigg\{&n_i\sum_{k\ne i}\nu_{ik}^{(0)}\big(M_{ik}^{(0)}-M_i\big)-n_iM_i\Bigg[\frac{1}{n_iT_i}\,c_i\cdot\sum_{k\ne i}R_{ik}+n_iA_i:\nabla_xu_i\\
&+\frac{2}{3}\frac{1}{n_iT_i}\Big(\frac{m_i}{2}\frac{|c_i|^2}{T_i}-\frac{3}{2}\Big)\Big(\sum_{k\ne i}S_{ik}-u_i\cdot\sum_{k\ne i}R_{ik}\Big)+\frac{n_i}{\sqrt{m_i}}\,B_i\cdot\frac{\nabla_xT_i}{T_i^{1/2}}\Bigg]\Bigg\}.
\end{aligned}
\end{equation}

The idea of \cite{lev-mor-nad-1998} is the following. If $f_i\ge0$, then $\int\psi^2f_i\,dv\ge0$ for every function $\psi=\psi(v)$. Taking $\psi=a^\intercal K_i$, where $K_i(v)$ is a vector of $p$ polynomials and $a\in\RR^p$, we get
\begin{equation}
a^\intercal\Big(\int_{\RR^3}K_iK_i^\intercal f_i\,dv\Big)a=\int_{\RR^3}(a^\intercal K_i)^2f_i\,dv\ge0\qquad\forall a\in\RR^p ,
\end{equation}
so the matrix
\begin{equation}
\mathcal W_i=\frac{1}{n_i}\int_{\RR^3}K_iK_i^\intercal f_i\,dv
\end{equation}
is positive semi-definite for every non-negative $f_i$. When $f_i$ is replaced by a truncated Chapman--Enskog expansion, $\mathcal W_i$ becomes an explicit function of macroscopic quantities. If it is not positive semi-definite, the truncated expansion cannot correspond to a non-negative distribution function, and the corresponding fluid closure should not be used. In \cite{lev-mor-nad-1998} this was used to switch between the Euler and Navier--Stokes equations for a single gas; as in \cite{fil-rey-2014, cap-par-rey-2025, tiw-2000, cap-par-rey-2026}, we use it to switch between fluid and kinetic descriptions, here separately for each species.

We choose
\begin{equation}
K_i=\begin{bmatrix}1\\[4pt]\sqrt{m_i}\,\dfrac{c_i}{\sqrt{T_i}}\\[8pt]\dfrac{m_i}{2}\dfrac{|c_i|^2}{T_i}-\dfrac52\end{bmatrix}.
\end{equation}
Using \eqref{eqMoments}, in particular $\frac{1}{n_i}\int\frac{m_i|c_i|^2}{2T_i}f_i\,dv=\frac{3}{2}$, one obtains
\begin{equation}
\mathcal W_i=\begin{pmatrix}1&0^\intercal&-1\\0&I+\thickbar A_i&\thickbar B_i\\-1&\thickbar B_i^\intercal&\thickbar C_i\end{pmatrix},\qquad
\thickbar A_i=\int_{\RR^3}A_if_i\,dv,\quad \thickbar B_i=\int_{\RR^3}B_if_i\,dv,\quad \thickbar C_i=\int_{\RR^3}C_if_i\,dv .
\label{eq2}
\end{equation}
The quantities $\thickbar A_i$, $\thickbar B_i$, $\thickbar C_i$ depend on $f_i$, hence on the order at which the Chapman--Enskog expansion is truncated. We will compare the matrices obtained at two different orders: if they are close, the two descriptions are consistent; if not, the first-order corrections are important and the zeroth-order closure is not accurate.

We first derive a lower bound for $\thickbar C_i$, which is needed below. Expanding the square and using $\frac{1}{n_i}\int\frac{m_i|c_i|^2}{2T_i}f_i\,dv=\frac{3}{2}$, we get
\begin{equation}\label{eqTempC}
\thickbar C_i=\frac{1}{n_i}\int_{\RR^3}\Big[\frac{m_i}{2}\frac{|c_i|^2}{T_i}-\frac{5}{2}\Big]^2f_i\,dv
=\frac{1}{n_i}\Big[\frac{1}{4T_i^2}\int_{\RR^3}\big(m_i|c_i|^2\big)^2f_i\,dv-\frac{5}{4}\,n_i\Big].
\end{equation}
A lower bound for $\int(m_i|c_i|^2)^2f_i\,dv$ follows from H\"older's inequality
\begin{equation}\label{eq::Holder}
\int|g(v)h(v)|\,dv\le\Big(\int|g(v)|^p\,dv\Big)^{1/p}\Big(\int|h(v)|^q\,dv\Big)^{1/q},
\end{equation}
with $p=q=2$, $g=\sqrt{f_i}$ and $h=\sqrt{f_i}\,m_i|c_i|^2$.
By the equality condition in H\"older's inequality \cite{mit-1970}, equality holds only if $|g|^2$ and $|h|^2$ are proportional almost everywhere, that is, only if $m_i|c_i|^2$ is constant on the support of $f_i$. Under the non-degeneracy assumption that $f_i$ is not concentrated on a sphere $|c_i|=\mathrm{const}$, in both the continuous and the discrete setting, the inequality is therefore strict.
Since $\int g^2\,dv=n_i$ and $\int gh\,dv=\int m_i|c_i|^2f_i\,dv=3n_iT_i$, this gives
\[
9n_i^2T_i^2<n_i\int_{\RR^3}\big(m_i|c_i|^2\big)^2f_i\,dv,\qquad\text{i.e.}\qquad 9n_iT_i^2<\int_{\RR^3}\big(m_i|c_i|^2\big)^2f_i\,dv .
\]
Inserting this bound in \eqref{eqTempC} we obtain
\begin{equation}\label{eqTempC2}
\thickbar C_i>\frac{1}{n_i}\Big[\frac{9n_iT_i^2}{4T_i^2}-\frac{5}{4}\,n_i\Big]=1 .
\end{equation}
Since $\thickbar C_i>0$, we can take the Schur complement \cite{horn2012} of $\mathcal W_i$ with respect to $\thickbar C_i$: $\mathcal W_i$ is positive semi-definite if and only if
\begin{equation}
\mathcal R_i=\begin{pmatrix}1-\dfrac1{\thickbar C_i}&\dfrac{\thickbar B_i^\intercal}{\thickbar C_i}\\[8pt]\dfrac{\thickbar B_i}{\thickbar C_i}&I+\thickbar A_i-\dfrac{\thickbar B_i\thickbar B_i^\intercal}{\thickbar C_i}\end{pmatrix}
\end{equation}
is positive semi-definite. By \eqref{eqTempC2} the entry $1-1/\thickbar C_i$ is positive, and a second Schur complement gives
\begin{equation}\label{eq5}
\mathcal S_i=I+\thickbar A_i-\frac{\thickbar B_i\thickbar B_i^\intercal}{\thickbar C_i}-\frac{\thickbar B_i\thickbar B_i^\intercal}{\thickbar C_i^2\,(1-1/\thickbar C_i)}=I+\thickbar A_i-\frac{\thickbar B_i\thickbar B_i^\intercal}{\thickbar C_i-1}.
\end{equation}
Hence $\mathcal W_i$ is positive semi-definite if and only if the $3\times3$ matrix $\mathcal S_i$ is.

We evaluate \eqref{eq2} with $f_i=f_i^{(0)}$ from \eqref{eq4} and with $f_i=f_i^{(0)}+\varepsilon f_i^{(1)}$ from \eqref{eq3}, and denote by $\mathcal S_i^{(0)}$ and $\mathcal S_i^{(1)}$ the corresponding matrices \eqref{eq5}. At zeroth order
\begin{equation}
\thickbar A_i^{(0)}=0,\qquad \thickbar B_i^{(0)}=0,\qquad \thickbar C_i^{(0)}=\frac{5}{2} ,
\end{equation}
and at first order
\begin{equation}
\begin{aligned}
\thickbar A_i^{(1)}&=\varepsilon\frac{m_i}{T_i\nu_{ii}^{(0)}}\sum_{k\ne i}\nu_{ik}^{(0)}\Big[(u_{ik}^{(0)}-u_i)\otimes(u_{ik}^{(0)}-u_i)-\frac{1}{3}\sigma_{ik}I\Big]-\frac{\varepsilon}{\nu_{ii}^{(0)}}D(u_i),\\
\thickbar B_i^{(1)}&=\varepsilon\frac{\sqrt{m_i}}{T_i^{3/2}\nu_{ii}^{(0)}}\sum_{k\ne i}\nu_{ik}^{(0)}(u_{ik}^{(0)}-u_i)\Big[\frac{5}{2}(T_{ik}^{(0)}-T_i)+\frac{1}{2}m_i\sigma_{ik}\Big]-\varepsilon\frac{5}{2}\frac{1}{\sqrt{m_i}\,\nu_{ii}^{(0)}}\frac{\nabla_xT_i}{\sqrt{T_i}},\\
\thickbar C_i^{(1)}&=\thickbar C_i^{(0)}+\varepsilon\frac{m_i^2}{T_i^2}\sum_{k\ne i}\frac{\nu_{ik}^{(0)}}{\nu_{ii}^{(0)}}\Big[\frac{\sigma_{ik}^2}{4}+\frac{5}{4}\,\frac{3(T_{ik}^{(0)})^2-6T_{ik}^{(0)}T_i+5T_i^2}{m_i^2}+\frac{5(T_{ik}^{(0)}-T_i)\sigma_{ik}}{2m_i}\Big]-\varepsilon\frac{5}{2}\sum_{k\ne i}\frac{\nu_{ik}^{(0)}}{\nu_{ii}^{(0)}},
\end{aligned}
\end{equation}
where $D(u_i)=\nabla_xu_i+(\nabla_xu_i)^\intercal-\frac{2}{3}(\nabla_x\cdot u_i)I$ and $\sigma_{ik}=|u_{ik}^{(0)}-u_i|^2$. Inserting these quantities in \eqref{eq5} gives $\mathcal S_i^{(0)}$ and $\mathcal S_i^{(1)}$. In particular $\mathcal S_i^{(0)}=I$.

The Euler description of species $i$ is considered adequate as long as $\mathcal S_i^{(1)}$ stays close to $\mathcal S_i^{(0)}$, that is, it remains positive definite and its eigenvalues stay close to those of $\mathcal S_i^{(0)}$. Species $i$ is switched from the Euler equations (Layer 0) to the BGK equation (Layer 1) when
\begin{equation}\label{crit:CompEuler}
\exists\,j\in\{1,2,3\}\quad\text{such that}\quad|\theta_{i,j}^1-\theta_{i,j}^0|>\eta_i ,
\end{equation}
where $\theta_{i,j}^0$ and $\theta_{i,j}^1$ are the eigenvalues of $\mathcal S_i^{(0)}$ and $\mathcal S_i^{(1)}$ (ordered consistently), and $\eta_i>0$ is a threshold. Here $\mathcal S_i^{(0)}=I$, so $\theta_{i,j}^0=1$ for $j=1,2,3$.

Since the indicator is computed for each species, the partition of the domain is in general different for different species.
\begin{rem}
    The bound \eqref{eqTempC2} relies on the non-negativity of $f_i$. Since the truncated expansion $f_i^{(0)}+\varepsilon f_i^{(1)}$ may take negative values, in principle the condition $\thickbar C_i^{(1)}>1$ is not guaranteed. However, $\thickbar C_i^{(1)}=\frac52+\mathcal O(\varepsilon)$, so for $\varepsilon$ small, as required by the Chapman--Enskog expansion, we still have $\thickbar C_i^{(1)}>1$, and the reduction \eqref{eq5} is well defined. The condition can fail only if the first-order correction to $\thickbar C_i$ is comparable with its zeroth-order value $\frac52$, that is, in extreme cases where the truncated expansion is no longer accurate. In our computations the indicator is evaluated through \eqref{eq5} without further checks on $\thickbar C_i^{(1)}$. In all the numerical tests of Section~\ref{secSimulations} the indicator was computed in every cell and at every time step, and no additional checks were needed, which indicates that the approach is robust in the regimes considered.
\end{rem}

\begin{rem}
The same criterion is used for the reverse transition: species $i$ is moved from Layer 1 (BGK) back to Layer 0 (Euler) only when
\begin{equation}
|\theta_{i,j}^1(t,x)-\theta_{i,j}^0(t,x)|\le\eta_i\qquad\forall\,j\in\{1,2,3\}.
\end{equation}
\end{rem}

\begin{rem}\label{rem:hierarchy}
The indicator is based on the Chapman--Enskog expansion, and the realizability matrices can be computed at any truncation order. The method is therefore not limited to two levels: it can be used to build a hierarchy of models including, for instance, Navier--Stokes or Burnett-type closures, as in \cite{cap-par-rey-2025, cap-par-rey-2026} for a single gas. In such a hierarchy, the passage from level $l$ to level $l+1$ is decided as in \eqref{crit:CompEuler}, by comparing the eigenvalues of the matrices $\mathcal S_i^{(l)}$ and $\mathcal S_i^{(l+1)}$ obtained at two consecutive orders of the expansion, with a threshold for each level.
\end{rem}

The spectral indicator measures the Chapman--Enskog corrections represented in these moments; it is not a bound on the distance of an arbitrary kinetic distribution from a Maxwellian. In particular, spatially homogeneous non-equilibrium modes with equilibrium hydrodynamic moments need not be detected. The switching strategy is therefore assessed here for the initial data and flow configurations of the numerical tests. When the initial data are not at equilibrium, a simple remedy is to compute the first time steps with the kinetic model everywhere: if the distribution functions are out of equilibrium, their evolution can produce non-equilibrium effects in the macroscopic quantities used by the indicator. The indicator can then identify the non-equilibrium regions, and the kinetic model is kept only there, as done for instance in \hyperref[Test4]{Test~4}.

\section{Numerical schemes}\label{secNumSch}
In this section we describe the velocity and space discretizations, the time integration and the coupling between the solvers, and the treatment that ensures the conservation of macroscopic quantities at the kinetic level.

\subsection{Velocity discretization}
The velocity space is truncated to a bounded domain, which neglects the tails of the distribution functions. Since the exact distributions do not have compact support, conservation of the macroscopic quantities is lost, and the same happens for the collision terms. To recover exact conservation we use the method of \cite{gam-tha-2009}, described in Section~\ref{secConMas}. The bounded domain $[-v_{\max},v_{\max}]^{d_v}$ is discretized by a uniform Cartesian grid with $N_v$ points per direction and step $\Delta v=2v_{\max}/(N_v-1)$; we denote by $v_m$, $m=1,\dots,N_V$, with $N_V=N_v^{d_v}$, the points of the grid. The distribution function of species $i$ is replaced by the values $f_{i,m}(t,x)\approx f_i(t,x,v_m)$, $m=1,\dots,N_V$, which solve
\begin{equation}
\partial_tf_{i,m}+v_m\cdot\nabla_xf_{i,m}=\frac{1}{\varepsilon}\nu_{ii}(n_iM_{ii,m}-f_{i,m})+\sum_{k\ne i}\nu_{ik}(n_iM_{ik,m}-f_{i,m}),\qquad m=1,\dots,N_V,
\end{equation}
where $M_{ik,m}$ denotes $M_{ik}$ evaluated at $v_m$.

\subsection{Spatial discretization}
We use a finite volume method with a MUSCL reconstruction on triangular meshes \cite{buf-cla-2010}. The domain is covered by a Delaunay triangulation made of $N_P$ triangles $P_c$, $c=1,\dots,N_P$, with barycentres $x_c$. The cell averages
\begin{equation}
\thickbar f_{i,m,c}^n=\frac{1}{|P_c|}\int_{P_c}f_{i,m}(t^n,x)\,dx,\qquad
\thickbar U_{i,c}^n=\frac{1}{|P_c|}\int_{P_c}U_i(t^n,x)\,dx,
\end{equation}
are the unknowns of the kinetic and fluid solvers, respectively, where $U_i=(n_i,\rho_iu_i,E_i)^\intercal$, $E_i=\frac{1}{2}\rho_i|u_i|^2+\frac{3}{2}n_iT_i$, are the conservative variables of species $i$ in \eqref{eq7}. In each cell we build second-order MUSCL reconstructions (as in \cite{buf-cla-2010}) of the distribution functions $f_{i,m,c}$ and of the primitive variables $(n_i,u_i,T_i)$, from which the conservative variables are computed at the edges; the numerical fluxes are computed with the local Lax--Friedrichs flux from the reconstructed values.

\subsection{IMEX relaxation step}\label{sec:numerical_properties}
We now study the implicit relaxation step of the IMEX scheme for the $N$-species model \eqref{rescaledBGK}. We first explain how the collision parameters are treated, then derive the linear systems for the velocities and temperatures (see also \cite{bis-bos-dim-gro-mar-2022}) and prove that they are uniquely solvable with positive temperatures. The AP property is proved in Section~\ref{subsec:AP_Euler}.

\subsubsection{Treatment of the collision parameters}\label{subsec:linearization}
In general the collision frequencies $\nu_{ik}$ may depend on the macroscopic moments. A fully implicit treatment would make the moment systems nonlinear and would require an iterative solver. In what follows the collision frequencies are frozen at the explicit state. The density equations decouple and are solved first ($n_k^{n+1}=n_k^*$, see below), so that the coefficients \eqref{eq_aik_bik_gammaik} can be computed with the new densities and the frozen frequencies. The relaxation steps for velocities and temperatures then reduce to linear systems.

\subsubsection{Macroscopic moment systems}\label{subsec:solvability}
Consider the first-order IMEX update for species $i$,
\begin{equation}\label{eq:imex_update_formal}
f_i^{n+1}=f_i^*+\frac{\Delta t}{\varepsilon}\nu_{ii}\big(n_i^{n+1}M_{ii}^{n+1}-f_i^{n+1}\big)+\Delta t\sum_{k\ne i}\nu_{ik}\big(n_i^{n+1}M_{ik}^{n+1}-f_i^{n+1}\big),
\end{equation}
where $f_i^*=f_i^n-\Delta t\,v\cdot\nabla_xf_i^n$. We integrate \eqref{eq:imex_update_formal} against $\phi_i(v)=(1,\,m_iv,\,\frac{1}{2}m_i|v|^2)^\intercal$. Since $M_{ii}^{n+1}$ has the same moments as $f_i^{n+1}$, the stiff term disappears from the moment equations.

\paragraph{\bf Densities.} With $\phi_i=1$,
\begin{equation}
n_i^{n+1}=n_i^*,\qquad i=1,\dots,N,
\end{equation}
which fixes $a_{ik}^{n+1}$, $b_{ik}^{n+1}$ and $\gamma_{ik}^{n+1}$ through \eqref{eq_aik_bik_gammaik}.

\paragraph{\bf Velocities.} With $\phi_i=m_iv$,
\begin{equation}
\rho_i^{n+1}u_i^{n+1}=\rho_i^*u_i^*+\Delta t\sum_{k\ne i}\nu_{ik}\rho_i^{n+1}\big(u_{ik}^{n+1}-u_i^{n+1}\big).
\end{equation}
Using $u_{ik}^{n+1}=(1-a_{ik}^{n+1})u_i^{n+1}+a_{ik}^{n+1}u_k^{n+1}$ and $\rho_i^{n+1}=\rho_i^*$, we obtain for each spatial component the $N\times N$ linear system
\begin{equation}\label{eq:ns_velocity_system}
\Big(1+\Delta t\sum_{k\ne i}\nu_{ik}a_{ik}^{n+1}\Big)u_i^{n+1}-\Delta t\sum_{k\ne i}\nu_{ik}a_{ik}^{n+1}u_k^{n+1}=u_i^*,\qquad i=1,\dots,N.
\end{equation}

\paragraph{\bf Temperatures.} With $\phi_i=\frac{1}{2}m_i|v|^2$ we obtain the energy balance. We write $E_i=\mathcal K_i+\frac{3}{2}n_iT_i$, with $\mathcal K_i=\frac{1}{2}\rho_i|u_i|^2$, and similarly $E_{ik}=\mathcal K_{ik}+\frac{3}{2}n_iT_{ik}$, $\mathcal K_{ik}=\frac{1}{2}\rho_i|u_{ik}|^2$. The velocities $u_i^{n+1}$ are known from \eqref{eq:ns_velocity_system}. Using \eqref{eqUT} for $T_{ik}^{n+1}$ we obtain
\begin{equation}\label{eq:ns_temperature_system}
\Big(1+\Delta t\sum_{k\ne i}\nu_{ik}b_{ik}^{n+1}\Big)T_i^{n+1}-\Delta t\sum_{k\ne i}\nu_{ik}b_{ik}^{n+1}T_k^{n+1}=r_i^*,\qquad i=1,\dots,N,
\end{equation}
with
\begin{equation}\label{eq:Sstar_def}
r_i^*=\frac{2}{3n_i^{n+1}}\Big[E_i^*-\mathcal K_i^{n+1}+\Delta t\sum_{k\ne i}\nu_{ik}\Big(\mathcal K_{ik}^{n+1}-\mathcal K_i^{n+1}+\frac{3}{2}n_i^{n+1}\gamma_{ik}^{n+1}|u_i^{n+1}-u_k^{n+1}|^2\Big)\Big].
\end{equation}

\subsubsection{Well-posedness of the relaxation systems}
Both \eqref{eq:ns_velocity_system} and \eqref{eq:ns_temperature_system} have the form $\mathcal Ax=y$ with
\begin{equation}
\mathcal A_{ii}=1+\sum_{k\ne i}c_{ik},\qquad\mathcal A_{ik}=-c_{ik}\quad(k\ne i),
\end{equation}
where $c_{ik}=\Delta t\,\nu_{ik}a_{ik}^{n+1}$ or $c_{ik}=\Delta t\,\nu_{ik}b_{ik}^{n+1}$. We first collect some properties of the parameters.

\begin{lemma}[Admissibility of the parameters]\label{lem:parameters}
Let $g_{i,k}(\mu)\ge0$, so that $\lambda_{ik}\ge0$, and assume that the collision frequencies $\nu_{ik}>0$ satisfy
\begin{equation}\label{eq:nu_bound}
\nu_{ik}\ge\frac{\lambda_{ik}n_k^{n+1}}{2}
\end{equation}
(see also \cite{bis-gro-mar2021}). Then, for every $n_k^{n+1}\ge0$, the parameters \eqref{eq_aik_bik_gammaik} satisfy
\begin{equation}
a_{ik}^{n+1}\ge0,\qquad\gamma_{ik}^{n+1}\ge0,\qquad0\le b_{ik}^{n+1}\le1 .
\end{equation}
\end{lemma}

\begin{proof}
Since $g_{i,k}\ge0$ and $1-\mu\ge0$ on $[-1,1]$, $\lambda_{ik}\ge0$; with $\nu_{ik}>0$ and $n_k^{n+1}\ge0$ we get $a_{ik}^{n+1}\ge0$ and therefore $b_{ik}^{n+1}\ge0$. By \eqref{eq_aik_bik_gammaik}, $\gamma_{ik}^{n+1}\ge0$ if and only if $a_{ik}^{n+1}\le2m_k/(m_i+m_k)$, that is
\[
\frac{\lambda_{ik}n_k^{n+1}m_k}{\nu_{ik}(m_i+m_k)}\le\frac{2m_k}{m_i+m_k}\iff\nu_{ik}\ge\frac{\lambda_{ik}n_k^{n+1}}{2},
\]
which is \eqref{eq:nu_bound}. Finally, from $a_{ik}^{n+1}\le2m_k/(m_i+m_k)$,
\[
b_{ik}^{n+1}\le\frac{2m_i}{m_i+m_k}\,\frac{2m_k}{m_i+m_k}=\frac{4m_im_k}{(m_i+m_k)^2}\le1 .
\]
\end{proof}

\begin{theorem}[Unique solvability and positivity]\label{teo1}
Under the assumptions of Lemma~\ref{lem:parameters}, assume that $n_i^*>0$ and $T_i^*>0$ for all $i$. Then for every $\Delta t>0$ the systems \eqref{eq:ns_velocity_system} and \eqref{eq:ns_temperature_system} have a unique solution, and $T_i^{n+1}>0$ for $i=1,\dots,N$.
\end{theorem}

\begin{proof}
\emph{Step 1 (sign of the coefficients).} Since $n_i^{n+1}=n_i^*>0$, by Lemma~\ref{lem:parameters} $c_{ik}\ge0$ for $i\ne k$.

\emph{Step 2 (invertibility).} For every row, $|\mathcal A_{ii}|=1+\sum_{k\ne i}c_{ik}>\sum_{k\ne i}|\mathcal A_{ik}|$, so $\mathcal A$ is strictly diagonally dominant and, by the L\'evy--Desplanques theorem \cite{tyr-1997}, non-singular for every $\Delta t>0$.
Hence \eqref{eq:ns_velocity_system} has a unique solution $u^{n+1}$ for each spatial component; with $u^{n+1}$ known, the right-hand side $r_i^*$ is determined, and \eqref{eq:ns_temperature_system} has a unique solution as well.

\emph{Step 3 (lower bound for the temperatures).} Let $T_m^{n+1}=\min_iT_i^{n+1}$. Row $m$ of \eqref{eq:ns_temperature_system} gives
\[
T_m^{n+1}=r_m^*+\sum_{k\ne m}c_{mk}\big(T_k^{n+1}-T_m^{n+1}\big)\ge r_m^*,
\]
so $T_i^{n+1}\ge\min_j r_j^*$ for all $i$.

\emph{Step 4 (positivity of $r_i^*$).} Let $V_i=\rho_i^{n+1}(u_i^{n+1}-u_i^*)=\Delta t\sum_{k\ne i}\nu_{ik}\rho_i^{n+1}(u_{ik}^{n+1}-u_i^{n+1})$. Since $E_i^*=\mathcal K_i^*+\frac{3}{2}n_i^*T_i^*$,
\[
E_i^*-\mathcal K_i^{n+1}=\frac{3}{2}n_i^*T_i^*-u_i^{n+1}\cdot V_i+\frac{|V_i|^2}{2\rho_i^{n+1}} .
\]
Using $\frac{1}{2}(|x|^2-|y|^2)=y\cdot(x-y)+\frac{1}{2}|x-y|^2$,
\[
\Delta t\sum_{k\ne i}\nu_{ik}\big(\mathcal K_{ik}^{n+1}-\mathcal K_i^{n+1}\big)=u_i^{n+1}\cdot V_i+\frac{\Delta t}{2}\sum_{k\ne i}\nu_{ik}\rho_i^{n+1}|u_{ik}^{n+1}-u_i^{n+1}|^2 .
\]
Inserting both relations in \eqref{eq:Sstar_def}, the terms $u_i^{n+1}\cdot V_i$ cancel. Moreover $u_{ik}^{n+1}-u_i^{n+1}=-a_{ik}^{n+1}(u_i^{n+1}-u_k^{n+1})$, and by \eqref{eq_aik_bik_gammaik}
\[
\frac{1}{2}\rho_i(a_{ik})^2|u_i-u_k|^2+\frac{3}{2}n_i\gamma_{ik}|u_i-u_k|^2=\rho_ia_{ik}\frac{m_k}{m_i+m_k}|u_i-u_k|^2 .
\]
Therefore
\begin{equation}\label{eq:Sstar}
r_i^*=T_i^*+\frac{2}{3n_i^{n+1}}\Big[\frac{|V_i|^2}{2\rho_i^{n+1}}+\Delta t\sum_{k\ne i}\nu_{ik}\rho_i^{n+1}a_{ik}^{n+1}\frac{m_k}{m_i+m_k}|u_i^{n+1}-u_k^{n+1}|^2\Big]\ge T_i^*>0 ,
\end{equation}
since $a_{ik}^{n+1}\ge0$. Together with Step 3 this gives $T_i^{n+1}>0$.
\end{proof}

\begin{rem}[Mixed Maxwellians]
The update requires $T_{ik}^{n+1}>0$. By \eqref{eqUT}, $T_{ik}^{n+1}=(1-b_{ik}^{n+1})T_i^{n+1}+b_{ik}^{n+1}T_k^{n+1}+\gamma_{ik}^{n+1}|u_i^{n+1}-u_k^{n+1}|^2$. By Theorem~\ref{teo1} and Lemma~\ref{lem:parameters}, the first two terms are a convex combination of positive temperatures and the last one is non-negative; hence $T_{ik}^{n+1}>0$ and all the Maxwellians $M_{ik}^{n+1}$ are well defined.
\end{rem}

\begin{rem}[Other scalings]\label{rem:scalings}
The intra-species term does not enter \eqref{eq:ns_velocity_system}--\eqref{eq:ns_temperature_system}. If the collisions between different species are scaled by a factor $\beta>0$, that is, $\sum_{k\ne i}\nu_{ik}(n_iM_{ik}-f_i)$ is replaced by $\beta\sum_{k\ne i}\nu_{ik}(n_iM_{ik}-f_i)$, the two systems keep the same form with $\Delta t\,\nu_{ik}$ replaced by $\beta\,\Delta t\,\nu_{ik}$ for $k\ne i$. The parameters \eqref{eq_aik_bik_gammaik} and condition \eqref{eq:nu_bound} do not depend on $\Delta t$ or $\beta$, so Lemma~\ref{lem:parameters}, Theorem~\ref{teo1} and formula \eqref{eq:Sstar} hold without changes. In particular they hold for $\beta=1/\varepsilon$, when collisions between different species are also stiff. This case is considered in \cite{PaperB}.
\end{rem}

Once the macroscopic state is known, the Maxwellians $M_{ii}^{n+1}$, $M_{ik}^{n+1}$ are fixed and \eqref{eq:imex_update_formal} is solved pointwise in velocity:
\begin{equation}\label{eq:final_pointwise_update}
f_i^{n+1}=\frac{f_i^*+\frac{\Delta t}\varepsilon\nu_{ii}n_i^{n+1}M_{ii}^{n+1}+\Delta t\sum_{k\ne i}\nu_{ik}n_i^{n+1}M_{ik}^{n+1}}{1+\frac{\Delta t}\varepsilon\nu_{ii}+\Delta t\sum_{k\ne i}\nu_{ik}} .
\end{equation}

The positivity result concerns the local implicit relaxation step with positive input densities and temperatures. It does not by itself prove positivity of the preceding high-order transport and Runge--Kutta stages.

\subsection{Asymptotic-preserving property}\label{subsec:AP_Euler}

\begin{theorem}[Euler limit, first order]
Under the assumptions of Theorem~\ref{teo1}, for $\varepsilon\to0$ and fixed $\Delta t>0$, $f_i^{n+1}$ given by \eqref{eq:final_pointwise_update} converges to $n_i^{n+1}M_{ii}^{n+1}$, and the scheme reduces to a consistent time discretization of \eqref{eq7}, with no restriction of order $\varepsilon$ on $\Delta t$.
\end{theorem}

\begin{proof}
Multiplying numerator and denominator of \eqref{eq:final_pointwise_update} by $\varepsilon$ and letting $\varepsilon\to0$,
\[
\lim_{\varepsilon\to0}f_i^{n+1}=\lim_{\varepsilon\to0}\frac{\varepsilon f_i^*+\Delta t\,\nu_{ii}n_i^{n+1}M_{ii}^{n+1}+\varepsilon\Delta t\sum_{k\ne i}\nu_{ik}n_i^{n+1}M_{ik}^{n+1}}{\varepsilon+\Delta t\,\nu_{ii}+\varepsilon\Delta t\sum_{k\ne i}\nu_{ik}}=n_i^{n+1}M_{ii}^{n+1} .
\]
We take the moments of \eqref{eq:imex_update_formal} with respect to $\phi_i$:
\begin{equation}\label{eq:discrete_moment_system}
\int_{\RR^3}\phi_i\Big(\frac{f_i^{n+1}-f_i^n}{\Delta t}+v\cdot\nabla_xf_i^n\Big)dv=\int_{\RR^3}\phi_i\Big[\frac{\nu_{ii}}{\varepsilon}\big(n_i^{n+1}M_{ii}^{n+1}-f_i^{n+1}\big)+\sum_{k\ne i}\nu_{ik}\big(n_i^{n+1}M_{ik}^{n+1}-f_i^{n+1}\big)\Big]dv .
\end{equation}
The stiff term vanishes for every $\varepsilon>0$, since $\phi_i$ are collision invariants of the intra-species operator. When $f_i^n=n_i^nM_{ii}^n$, which holds after the first step in the limit, the pressure tensor and the heat flux reduce to their equilibrium values,
\[
\mathbb P_i=\int_{\RR^3}m_i(v-u_i)\otimes(v-u_i)\,n_iM_{ii}\,dv=n_iT_iI,\qquad q_i=\frac{1}{2}\int_{\RR^3}m_i|v-u_i|^2(v-u_i)\,n_iM_{ii}\,dv=0,
\]
and \eqref{eq:discrete_moment_system} becomes
\begin{equation}
\begin{cases}
\displaystyle\frac{n_i^{n+1}-n_i^n}{\Delta t}+\nabla_x\cdot(n_iu_i)^n=0,\\[8pt]
\displaystyle\frac{(\rho_iu_i)^{n+1}-(\rho_iu_i)^n}{\Delta t}+\nabla_x\cdot\big(\rho_iu_i\otimes u_i+n_iT_iI\big)^n=\sum_{k\ne i}\nu_{ik}\rho_i^{n+1}\big(u_{ik}^{n+1}-u_i^{n+1}\big),\\[8pt]
\displaystyle\frac{E_i^{n+1}-E_i^n}{\Delta t}+\nabla_x\cdot\big((E_i+n_iT_i)u_i\big)^n=\sum_{k\ne i}\nu_{ik}\big(E_{ik}^{n+1}-E_i^{n+1}\big).
\end{cases}
\end{equation}
Using \eqref{eqUT}--\eqref{eq_aik_bik_gammaik} one checks that $\nu_{ik}\rho_i(u_{ik}-u_i)=R_{ik}$ and $\nu_{ik}(E_{ik}-E_i)=S_{ik}$, with $R_{ik}$, $S_{ik}$ given in \eqref{eqRS}. The scheme is therefore a consistent discretization of \eqref{eq7}, with exchange terms treated implicitly.
\end{proof}

\begin{theorem}[Euler limit for ARS(2,2,2)]
Let the initial data be well prepared, $f_i^n=n_i^nM_{ii}^n$, and assume that the explicit states entering the two implicit stages satisfy the assumptions of Theorem~\ref{teo1}. For $\varepsilon\to0$ the stages $f_i^{(2)}$ and $f_i^{n+1}$ of the ARS(2,2,2) scheme \cite{bos-dim-2021, asc-ruu-spi-1997,bos-par-rus-2024} converge to the local Maxwellians $n_i^{(2)}M_{ii}^{(2)}$ and $n_i^{n+1}M_{ii}^{n+1}$, and the scheme reduces to a second-order IMEX discretization of \eqref{eq7}, with no restriction of order $\varepsilon$ on $\Delta t$ (see also \cite{bis-bos-dim-gro-mar-2022}).
\end{theorem}

\begin{proof}
Let $g_i(f)=-v\cdot\nabla_xf_i$, $\mathcal Q_{ii}(f_i)=\nu_{ii}(n_iM_{ii}-f_i)$ and $\mathcal Q_{i,\rm inter}(f)=\sum_{k\ne i}\nu_{ik}(n_iM_{ik}-f_i)$. With the tableaux of Section~\ref{sub:UnifiedIMEX}, the stages are $f_i^{(1)}=f_i^n$,
\begin{align}
f_i^{(2)}&=f_i^n+\gamma\Delta t\,g_i(f^{(1)})+\gamma\Delta t\Big[\frac{1}{\varepsilon}\mathcal Q_{ii}(f_i^{(2)})+\mathcal Q_{i,\rm inter}(f^{(2)})\Big],\label{eq:ars_stage2}\\
f_i^{n+1}&=f_i^n+\Delta t\big[\delta\,g_i(f^{(1)})+(1-\delta)\,g_i(f^{(2)})\big]\nonumber\\
&\quad+\Delta t\Big[(1-\gamma)\Big(\frac{1}{\varepsilon}\mathcal Q_{ii}(f_i^{(2)})+\mathcal Q_{i,\rm inter}(f^{(2)})\Big)+\gamma\Big(\frac{1}{\varepsilon}\mathcal Q_{ii}(f_i^{n+1})+\mathcal Q_{i,\rm inter}(f^{n+1})\Big)\Big],\label{eq:ars_stage3}
\end{align}
where we used that the scheme is globally stiffly accurate, so that $f_i^{n+1}=f_i^{(3)}$. Multiplying \eqref{eq:ars_stage2} by $\varepsilon$ and letting $\varepsilon\to0$ gives $\mathcal Q_{ii}(f_i^{(2)})=0$, that is $f_i^{(2)}\to n_i^{(2)}M_{ii}^{(2)}$; moreover \eqref{eq:ars_stage2} shows that $\mathcal Q_{ii}(f_i^{(2)})=\mathcal O(\varepsilon)$. Multiplying \eqref{eq:ars_stage3} by $\varepsilon$ we then get $\mathcal Q_{ii}(f_i^{n+1})=0$ in the limit, i.e. $f_i^{n+1}\to n_i^{n+1}M_{ii}^{n+1}$.

Projecting the stages on $\phi_i$, the terms $\mathcal Q_{ii}$ vanish for every $\varepsilon>0$, and since all stages are Maxwellians the fluxes are the Euler fluxes $F_i(U_i)=\big(n_iu_i,\ \rho_iu_i\otimes u_i+n_iT_iI,\ (E_i+n_iT_i)u_i\big)^\intercal$, with $U_i=(n_i,\rho_iu_i,E_i)^\intercal$. We obtain
\begin{align*}
&\frac{U_i^{(2)}-U_i^n}{\gamma\Delta t}+\nabla_x\cdot F_i(U_i^n)=\mathcal X_{i,\rm inter}^{(2)},\\
&\frac{U_i^{n+1}-U_i^n}{\Delta t}+\delta\,\nabla_x\cdot F_i(U_i^n)+(1-\delta)\,\nabla_x\cdot F_i(U_i^{(2)})=(1-\gamma)\,\mathcal X_{i,\rm inter}^{(2)}+\gamma\,\mathcal X_{i,\rm inter}^{n+1},
\end{align*}
with
\[
\mathcal X_{i,\rm inter}^{(s)}=\int_{\RR^3}\phi_i\,\mathcal Q_{i,\rm inter}(f^{(s)})\,dv=\Big(0,\ \sum_{k\ne i}R_{ik}^{(s)},\ \sum_{k\ne i}S_{ik}^{(s)}\Big)^\intercal .
\]
This is the ARS(2,2,2) discretization of \eqref{eq7}.
\end{proof}

\subsection{Time integration of the hybrid scheme}\label{sec:UnifiedFramework}
We now describe the time integration of the coupled problem and the treatment of the interfaces between regions where different models are used.

\subsubsection{IMEX scheme for the two models}\label{sub:UnifiedIMEX}
We consider two models: $\mathscr M_0$, the Euler equations \eqref{eq7} (Layer 0), and $\mathscr M_1$, the BGK model \eqref{rescaledBGK} (Layer 1). The construction below extends to more models. For every species $i$ and every cell $P_c$, the active model at time $t^n$ is $\mathscr M_l$, with $l=l_{i,c}^n\in\{0,1\}$; different species can use different models in the same cell. For each species, the solution in a cell is advanced with the active model only, that is, with the distribution function if the species is kinetic and with the conservative variables if it is in the Euler regime; the variables of the other model are recomputed from these when they are needed by the coupling described in Section~\ref{sub:coupling}. We write the equations of species $i$ in cell $c$ as $\partial_tU_{i,c}=F_{i,c}^{l}(U)+G_{i,c}^{l}(U)$, where $U$ collects the variables of all species (distribution functions or conservative variables) and:
\begin{itemize}
\item in $\mathscr M_1$, the transport term $-v\cdot\nabla_xf_i$ is in the explicit part $F_{i,c}^1$, and the collision operator (intra- and inter-species) is in the implicit part $G_{i,c}^1$, solved as in Section~\ref{sec:numerical_properties};
\item in $\mathscr M_0$, the fluxes are in the explicit part $F_{i,c}^0$, and the exchange terms are in the implicit part $G_{i,c}^0=\big(0,\sum_{k\ne i}R_{ik},\sum_{k\ne i}S_{ik}\big)^\intercal$.
\end{itemize}
We use a $q$-stage IMEX Runge--Kutta scheme with explicit tableau $\tilde A=(\tilde a_{sj})$, $\tilde a_{sj}=0$ for $j\ge s$, and diagonally implicit tableau $A=(a_{sj})$, $a_{sj}=0$ for $j>s$, with weights $\tilde\omega$, $\omega$ and nodes $\tilde c$, $c$:
\begin{center}
\begin{tabular}{c|c}$\tilde c$&$\tilde A$\\\hline&$\tilde\omega^\intercal$\end{tabular}\qquad\qquad
\begin{tabular}{c|c}$c$&$A$\\\hline&$\omega^\intercal$\end{tabular}
\end{center}
Writing $F_{i,c,j}=F_{i,c}^l(t^n+\tilde c_j\Delta t,U^{(j)})$ and $G_{i,c,j}=G_{i,c}^l(t^n+c_j\Delta t,U^{(j)})$, each stage $s=1,\dots,q$ is computed in two steps:
\begin{enumerate}
\item explicit step, in which all the terms are known from the previous stages:
\begin{equation}\label{eq:IMEX_system}
U_{i,c}^{(s)*}=U_{i,c}^n+\Delta t\sum_{j=1}^{s-1}\tilde a_{sj}F_{i,c,j}+\Delta t\sum_{j=1}^{s-1}a_{sj}G_{i,c,j};
\end{equation}
\item implicit step, solved as in Section~\ref{sec:numerical_properties} with $\Delta t$ replaced by $a_{ss}\Delta t$:
\begin{equation}\label{eq:ImplicitCorrector}
U_{i,c}^{(s)}=U_{i,c}^{(s)*}+\Delta t\,a_{ss}\,G_{i,c}^l(t^n+c_s\Delta t,U^{(s)}) .
\end{equation}
In Layer 0 this step only involves the exchange terms. The densities do not change and, since $\nu_{ik}\rho_i(u_{ik}-u_i)=R_{ik}$ and $\nu_{ik}(E_{ik}-E_i)=S_{ik}$, the velocities and temperatures solve \eqref{eq:ns_velocity_system}--\eqref{eq:ns_temperature_system}, with $\Delta t$ replaced by $a_{ss}\Delta t$ and the starred quantities computed from $U_{i,c}^{(s)*}$. These are the same systems as in Layer 1, so Theorem~\ref{teo1} applies.
\end{enumerate}
The solution is then
\begin{equation}
U_{i,c}^{n+1}=U_{i,c}^n+\Delta t\sum_{j=1}^q\tilde\omega_jF_{i,c,j}+\Delta t\sum_{j=1}^q\omega_jG_{i,c,j}.
\end{equation}
In Layer 1 the implicit step acts on the BGK operator. Since we use a second-order cell-centred finite volume scheme, the collision operator is evaluated at the cell average only, and the implicit step is local to each cell and can be solved explicitly as in Section~\ref{sec:numerical_properties}, without nonlinear solvers.

\begin{rem}[Higher order]\label{rem:penalization}
With a higher-order reconstruction, such as third-order WENO with several quadrature points per cell, or with the full Boltzmann operator, the implicit evaluation of the collision term at all quadrature points would be too expensive. One can then use a penalization \cite{fil-jin-2010, bos-dim-2021,bos-par-rus-2024}: a simple cell-centred BGK operator is subtracted from the high-order collision term and added back, and only this simple operator is treated implicitly. To control the explicit stiff residual, the IMEX scheme must then be globally stiffly accurate (GSA), i.e. $\omega_s=a_{q s}$ and $\tilde\omega_s=\tilde a_{q s}$. With the second-order scheme used here this is not needed. A third-order version of the hybrid method, with penalization, is developed in \cite{PaperB} for the regime in which the mixture relaxes to a single velocity and a single temperature.
\end{rem}

We use the ARS(2,2,2) scheme \cite{bos-dim-2021, asc-ruu-spi-1997,bos-par-rus-2024}, which is globally stiffly accurate:
\begin{center}
\begin{tabular}{c|ccc}
$0$ & $0$ & $0$ & $0$\\
$\gamma$ & $\gamma$ & $0$ & $0$\\
$1$ & $\delta$ & $1-\delta$ & $0$\\
\hline&$\delta$ & $1-\delta$ & $0$
\end{tabular}\qquad\qquad
\begin{tabular}{c|ccc}
$0$ & $0$ & $0$ & $0$\\
$\gamma$ & $0$ & $\gamma$ & $0$\\
$1$ & $0$ & $1-\gamma$ & $\gamma$\\
\hline& $0$ & $1-\gamma$ & $\gamma$
\end{tabular}
\end{center}
with
\begin{equation}
\gamma=1-\frac{1}{\sqrt2},\qquad\delta=1-\frac{1}{2\gamma}.
\end{equation}
This scheme was used for the Boltzmann equation on unstructured meshes in \cite{bos-dim-2021}.

\subsubsection{Interface coupling}\label{sub:coupling}
Let $\mathcal T$ be the mesh, made of the cells $P_c$, $c=1,\dots,N_P$. Since every species can change regime independently, we define for each species a regime flag
\begin{equation}
\chi_{i,c}\in\{\mathrm K,\mathrm E\},\qquad i=1,\dots,N,\ c=1,\dots,N_P,
\end{equation}
where $\mathrm K$ denotes the kinetic description (Layer 1) and $\mathrm E$ the Euler description (Layer 0). Let $\mathcal N_c$ be the stencil of cell $c$. Since different species may be in different regimes in the same cell, both the fluxes and the inter-species terms require some care.

\paragraph{\bf Kinetic cells.} Let $\chi_{i,c}=\mathrm K$.
\begin{enumerate}
\item \emph{Transport.} The flux of species $i$ requires its distribution function in the whole stencil $\mathcal N_c$. If $\chi_{i,j}=\mathrm E$ for some $j\in\mathcal N_c$, the distribution function in cell $j$ is replaced by the Maxwellian with the moments $U_{i,j}$:
\begin{equation}\label{eq:kinetic_coupling}
\widetilde f_{i,m,j}=\begin{cases}f_{i,m,j},&\chi_{i,j}=\mathrm K,\\
n_{i,j}M_{i,m}(u_{i,j}, T_{i,j}),&\chi_{i,j}=\mathrm E.\end{cases}
\end{equation}
Only quantities of species $i$ are involved.
\item \emph{Inter-species collisions.} The implicit step is local to the cell and is solved once for all species, with the systems \eqref{eq:ns_velocity_system}--\eqref{eq:ns_temperature_system}. Their right-hand sides contain the starred moments, that is, the densities, velocities and temperatures after the explicit step \eqref{eq:IMEX_system} of the current stage: for species $k$ they are computed by quadrature of $f_{k,c}^*$ if $\chi_{k,c}=\mathrm K$, and from $U_{k,c}^*$ if $\chi_{k,c}=\mathrm E$. The solution gives the new velocities and temperatures of all species in the cell. For the kinetic species they define the Maxwellians $M_{ik}$, towards which the distribution functions relax through \eqref{eq:final_pointwise_update}; for the fluid species they are directly the new macroscopic state.
\end{enumerate}

\paragraph{\bf Fluid cells.} Let $\chi_{i,c}=\mathrm E$.
\begin{enumerate}
\item \emph{Transport.} The Euler fluxes of species $i$ use only the conservative variables of species $i$ in $\mathcal N_c$. If $\chi_{i,j}=\mathrm K$, the moments are computed from the distribution function:
\begin{equation}\label{eq:macro_coupling}
\widetilde U_{i,j}=\begin{cases}U_{i,j},&\chi_{i,j}=\mathrm E,\\[4pt]\displaystyle\int_{\RR^{d_v}}\phi_i(v)\,f_{i,j}(v)\,dv,&\chi_{i,j}=\mathrm K.\end{cases}
\end{equation}
\item \emph{Exchange terms.} They are treated in the implicit step, together with all the species in the cell, as described above for kinetic cells; if $\chi_{k,c}=\mathrm K$, the moments of species $k$ are obtained by integrating $f_{k,c}^*$ in velocity.
\end{enumerate}

\paragraph{\bf Regime update and synchronization of the stages.} The coupling above is a local projection: moments are obtained from kinetic data by quadrature, and kinetic data are obtained from moments through Maxwellians. It is applied at two points of each time step.
\begin{enumerate}
\item \emph{Beginning of the time step.} The flags are updated from $\chi_{i,c}^n$ to $\chi_{i,c}^{n+1}$ with the criterion of Section~\ref{secIndicators}. Cells that change regime for species $i$ update their variables: if $\mathrm K\to\mathrm E$, $U_{i,c}^n$ is computed by quadrature of $f_{i,c}^n$; if $\mathrm E\to\mathrm K$, $f_{i,m,c}^n=n_{i,c}^nM_{i,m}(u_{i,c}^n,T_{i,c}^n)$. The first stage of ARS(2,2,2) has no implicit part and coincides with the solution at $t^n$, so its fluxes use the data at $t^n$, taken according to the previous flags. If $\chi_{i,c}^{n+1}=\mathrm K$, the transport of species $i$ in cell $c$ needs the distribution function of species $i$ in the stencil $\mathcal N_c$: for a neighbour $j$ with $\chi_{i,j}^{n}=\mathrm E$ the Maxwellian $n_{i,j}^nM_{i,m}(u_{i,j}^n,T_{i,j}^n)$, built from the moments of cell $j$, is used, and otherwise $f_{i,j}^n$. If $\chi_{i,c}^{n+1}=\mathrm E$, the Euler fluxes of species $i$ need its moments in the stencil: for a neighbour $j$ with $\chi_{i,j}^{n}=\mathrm K$ they are computed by quadrature of $f_{i,j}^n$, and otherwise $U_{i,j}^n$ is used. In particular, a neighbour that has just switched from $\mathrm K$ to $\mathrm E$ for species $i$ still provides its distribution function at $t^n$, since it was kinetic at that time.
\item \emph{After each stage.} After each implicit stage the coupling is applied again, now with the new flags $\chi_{i,c}^{n+1}$: in kinetic cells the moments are computed by quadrature of the current distribution functions, and in fluid cells the Maxwellians $n_{i,c}M_{i,m}(u_{i,c},T_{i,c})$ are built from the current moments. In particular, from the second stage on, a cell that has switched from $\mathrm K$ to $\mathrm E$ for species $i$ no longer provides $f_{i,c}^n$ to its kinetic neighbours, but the Maxwellian built from its current moments. All the other cells provide data according to their current regime. In this way, at each stage all the cells of a stencil provide data at the same stage. Exchanging data only at the end of the time step would introduce a time lag at the interfaces and reduce the accuracy there to first order \cite{cap-ten-2025}.
\end{enumerate}
The stage-wise coupling preserves the order of the IMEX scheme when the partition is fixed in time \cite{cap-ten-2025}. When the partition changes this is no longer guaranteed, but the first-order time lag at the interfaces is avoided.

\paragraph{\bf Lifting.} We lift fluid states with $n_iM_{i,m}(u_i,T_i)$. Non-equilibrium reconstructions, for example based on the Chapman--Enskog expansion, would be more accurate in principle but are not non-negative in general, and negative values of $f_{i,m}$ can make the computation unstable. The analytic Maxwellian used for lifting is non-negative. The subsequent conservative $L^2$ correction of Section~\ref{secConMas}, however, does not in general preserve pointwise non-negativity.

\subsection{Conservation of macroscopic quantities}\label{secConMas}
The discretization of the velocity space destroys the exact conservation of the macroscopic quantities \cite{mie-2000, dim-par-2014}. Taking a large velocity domain makes the truncation error small, but quadrature errors remain. Several remedies have been proposed, based on entropy minimization \cite{mie-2000} or on $L^2$ minimization \cite{gam-tha-2009, bos-cho-rus-2022}. The entropy approach gives a conservative discrete Maxwellian, while the $L^2$ approach applies to general distributions. We use the $L^2$ minimization of \cite{gam-tha-2009}, and apply it to every Maxwellian appearing in the scheme, for each species. A weighted version of this projection was recently used for gas mixtures in \cite{cap-par-rey-2025-2}.

\section{Numerical tests}\label{secSimulations}
We now present numerical tests of the hybrid scheme, including problems with non-trivial geometries. All computations were run on a laptop with an Intel Core Ultra 7 155H processor and 16 GB of RAM, without dedicated high-performance computing resources.

The solvers are second-order accurate in space and time. All tests consider a binary mixture with different molecular masses, specified in each test. Space and velocity are two-dimensional, and the velocity space is discretized with a uniform Cartesian grid. All the quantities of Sections~\ref{secModel} and~\ref{secNumSch} are computed accordingly in velocity dimension $d_v=2$; in particular $E_i=\frac{1}{2}\rho_i|u_i|^2+n_iT_i$ and the adiabatic exponent of the Euler equations is $2$. The indicators of Section~\ref{secIndicators} are computed with the three-dimensional formulas of that section, with the velocities and the gradients extended by a zero third component. These three-dimensional expressions are used as macroscopic switching indicators for the two-velocity-dimensional tests.
Unless stated otherwise, the time step is
\[
\Delta t=0.9\,\frac{\Delta x_{\min}}{v_{\max}},
\]
where $\Delta x_{\min}$ is the smallest edge length of the triangulation and the velocity grid covers $[-v_{\max},v_{\max}]^2$. The factor $0.9$ is replaced by $0.7$ in \hyperref[Test3]{Test~3}, and $\Delta t$ is varied in the convergence test (\hyperref[Test1]{Test~1}). In all tests $\nu_{ik}=\nu_{ik}^{(0)}=1$ and $\lambda_{ik}=1$.

\subsection{Test 1: order of convergence in time}\label{Test1}
We check the second-order convergence in time of the ARS(2,2,2) scheme, for both solvers and for the coupling of Section~\ref{sec:UnifiedFramework}. We consider two species with mass ratio $m_2/m_1=2$. The regime flag of every cell is drawn at random, independently for the two species, and kept fixed during the computation. As proved in \cite{cap-ten-2025}, with a frozen partition the coupling preserves the order of the IMEX scheme.

\begin{figure}[htbp]
\centering
\includegraphics[width=0.7\linewidth]{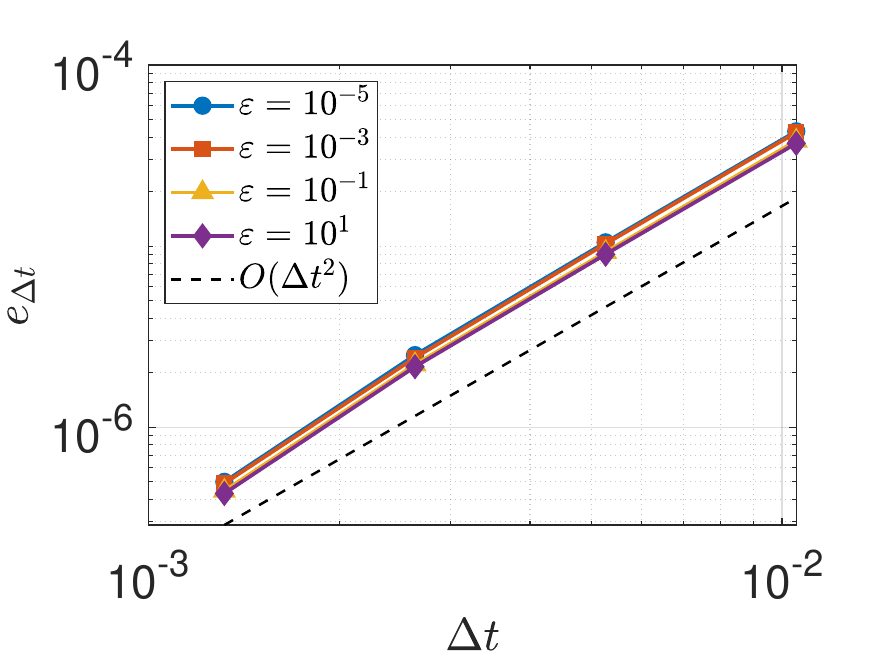}
\caption{\hyperref[Test1]{Test~1}. $L^1$ error $e_{\Delta t}$ of the total mass density as a function of $\Delta t$, for $\varepsilon=10,10^{-1},10^{-3},10^{-5}$, with a random partition fixed in time. The ARS(2,2,2) scheme is second-order accurate in all regimes.}
\label{Fig1}
\end{figure}

The domain $[0,7]\times[0,7]$ is discretized with about 2800 triangles. The macroscopic moments $W_i=(\rho_i,u_i^x,u_i^y,T_i)$ of the two species are initialized as
\begin{equation}
W_1(x,y)=W_2(x,y)=\big(1+0.1\exp(-\kappa r^2),\,0,\,0,\,1\big),
\end{equation}
with $r^2=(x-x_c)^2+(y-y_c)^2$, $x_c=y_c=7/2$ and $\kappa=3$; in kinetic cells $f_1$ and $f_2$ are the corresponding Maxwellians. The velocity domain $[-7,7]^2$ is discretized with $32$ points per direction. We take $\varepsilon\in\{10^{-5},10^{-3},10^{-1},10\}$ and, for each value, compute the $L^1$ error $e_{\Delta t}$ on the total mass density at the final time with respect to a reference solution obtained with the smallest time step. Figure~\ref{Fig1} shows second-order convergence in all regimes.

\subsection{Test 2: shock tube}\label{Test2}
This test, inspired by \cite{bis-bos-dim-gro-mar-2022}, is a Riemann problem for a binary mixture with $\varepsilon=10^{-5}$. In this collision-dominated regime the Euler equations are accurate, and the test is used to check that each solver is correct for mixtures and that the coupling does not produce spurious waves at the interfaces, and to show where the indicator activates the kinetic solver.

The domain $[0,1]\times[0,0.021]$ is discretized with 402 triangles, with zero-gradient boundary conditions. The problem is one-dimensional; the unstructured two-dimensional mesh, whose edges are not aligned with the flow direction, is used to check that the solution remains one-dimensional. The velocity domain $[-20,20]^2$ is discretized with $100$ points per direction, and the final time is $t_f=0.075$. The initial moments $W_i=(\rho_i,u_i^x,u_i^y,T_i)$ are
\begin{equation}
\begin{cases}
W_1(x<0.5,y)=(1,0,0,5),& W_2(x<0.5,y)=(1,0,0,15),\\
W_1(x\ge0.5,y)=(0.125,0,0,4),& W_2(x\ge0.5,y)=(0.125,0,0,12),
\end{cases}
\end{equation}
and in kinetic cells $f_1$, $f_2$ are the corresponding Maxwellians. The mass ratio is $m_2/m_1=10$, a value that already makes the kinetic computation demanding, and the thresholds are $\eta_1=\varepsilon$ and $\eta_2=\varepsilon/50$.

\begin{figure}
\centering
\includegraphics[width=0.96\textwidth]{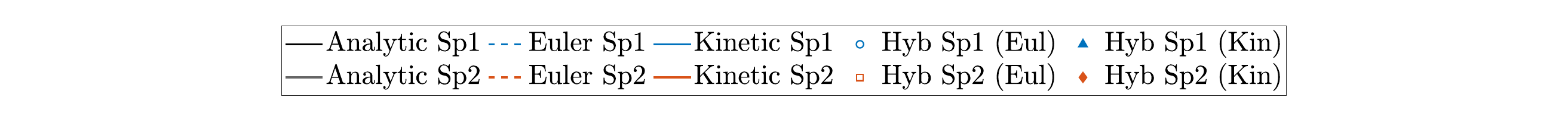}\par\smallskip
\figrow{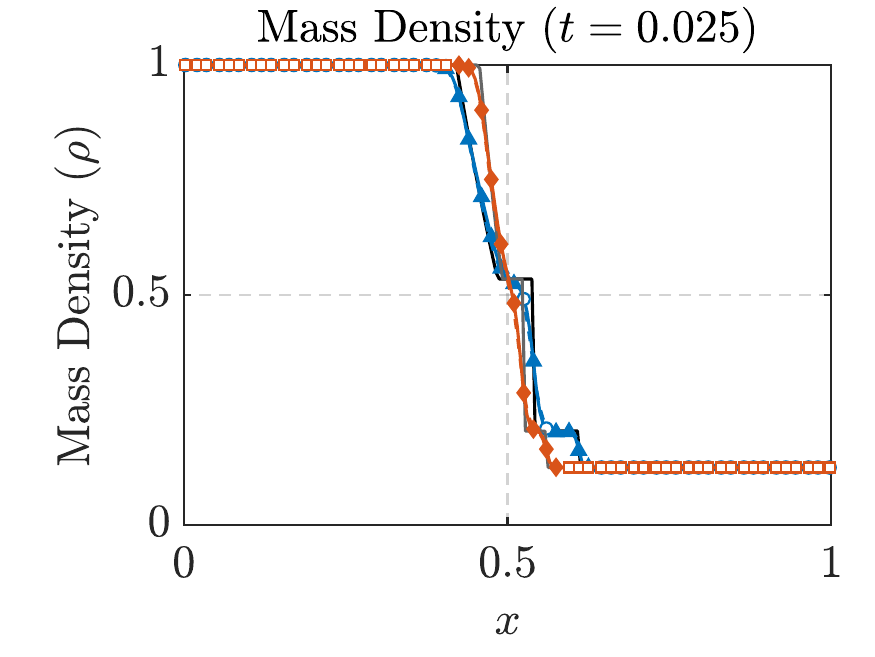}{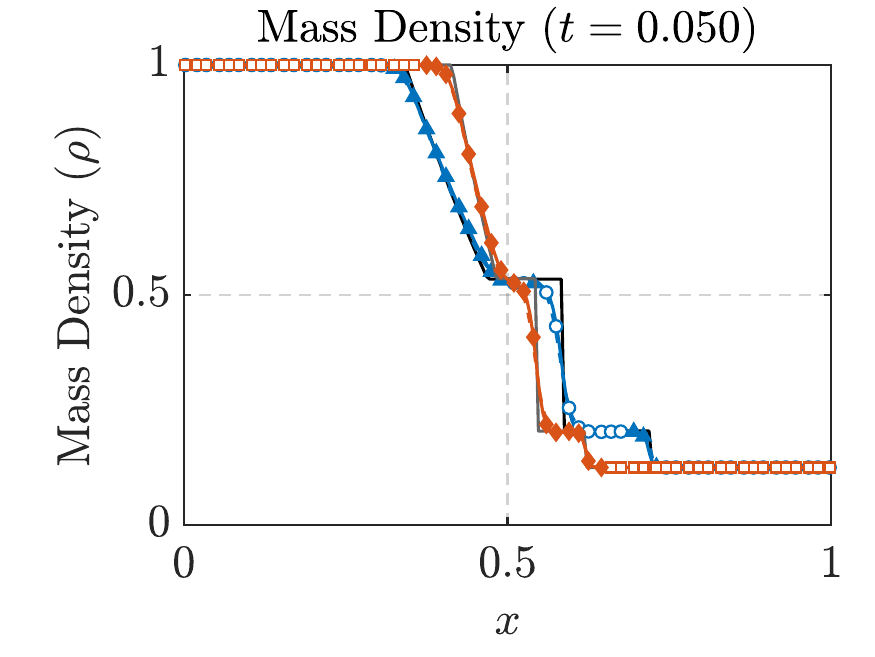}{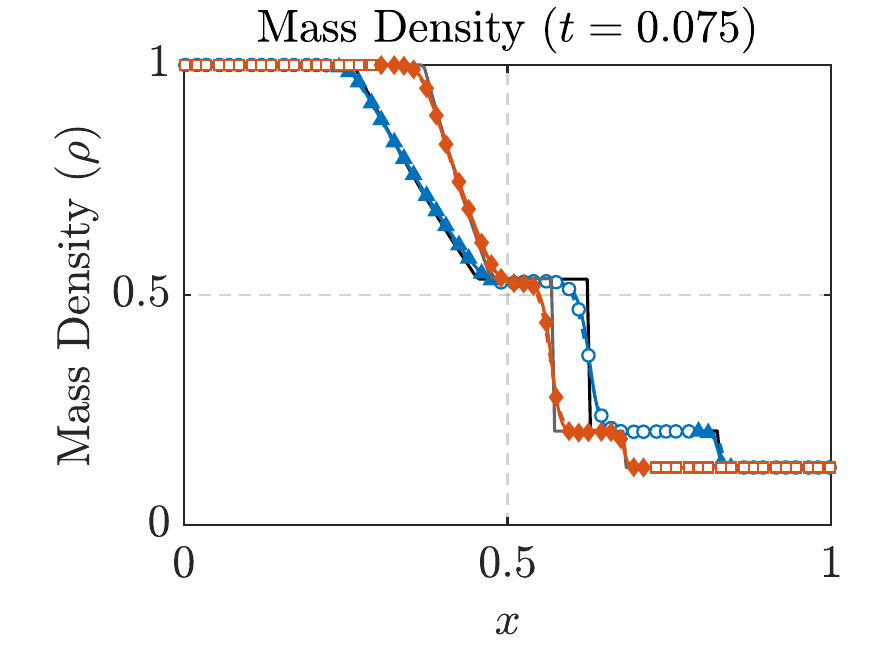}
\figrow{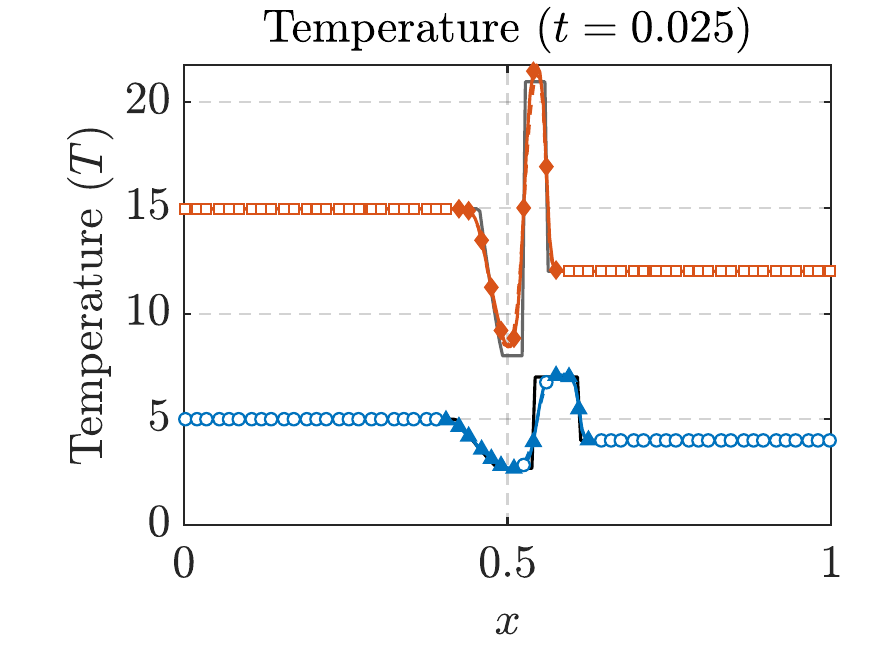}{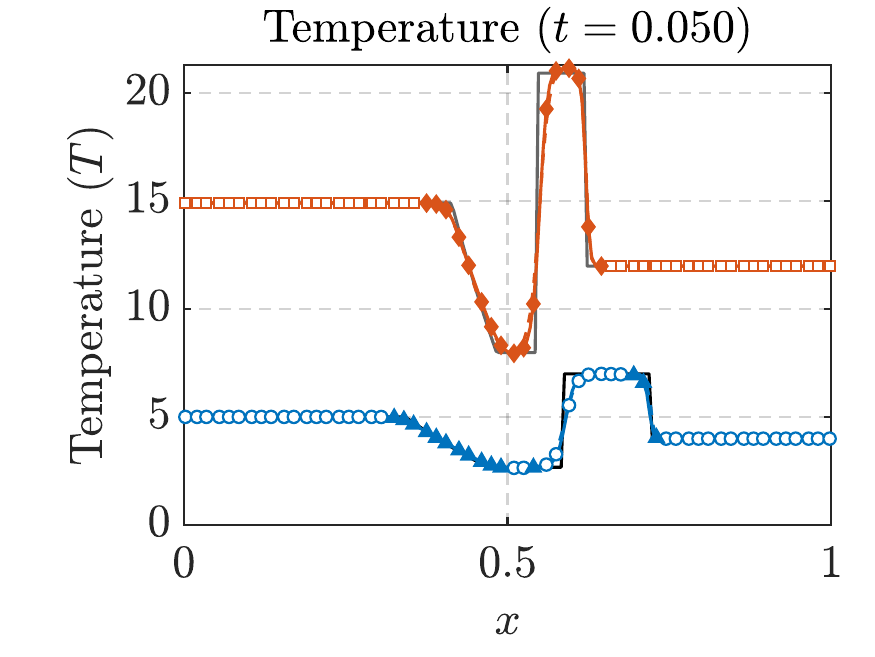}{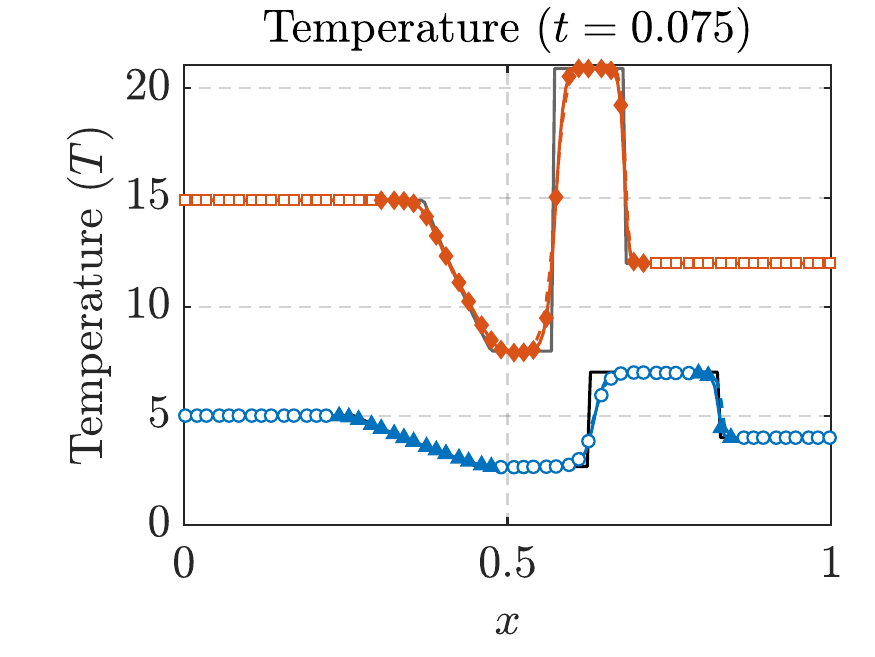}
\figrow{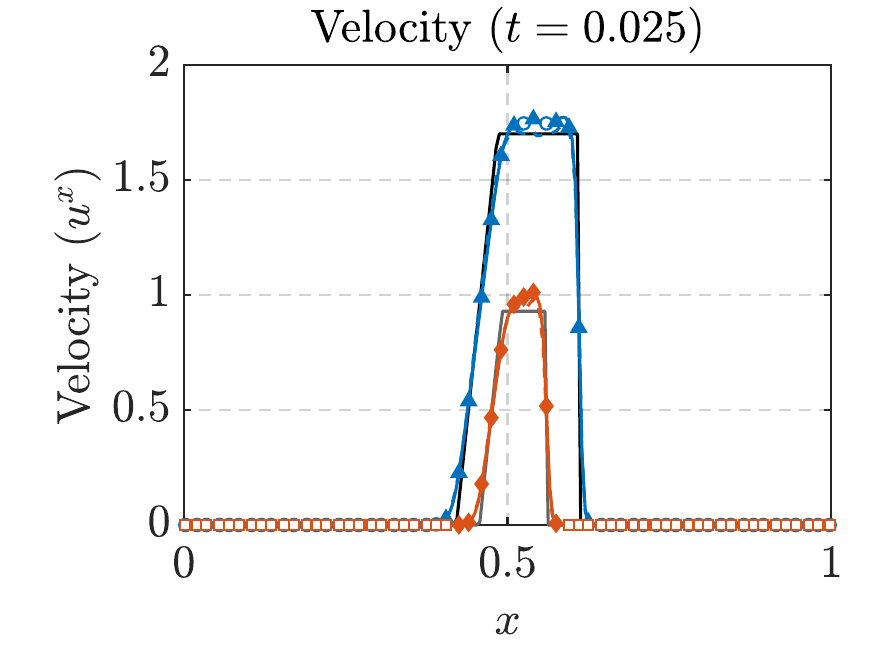}{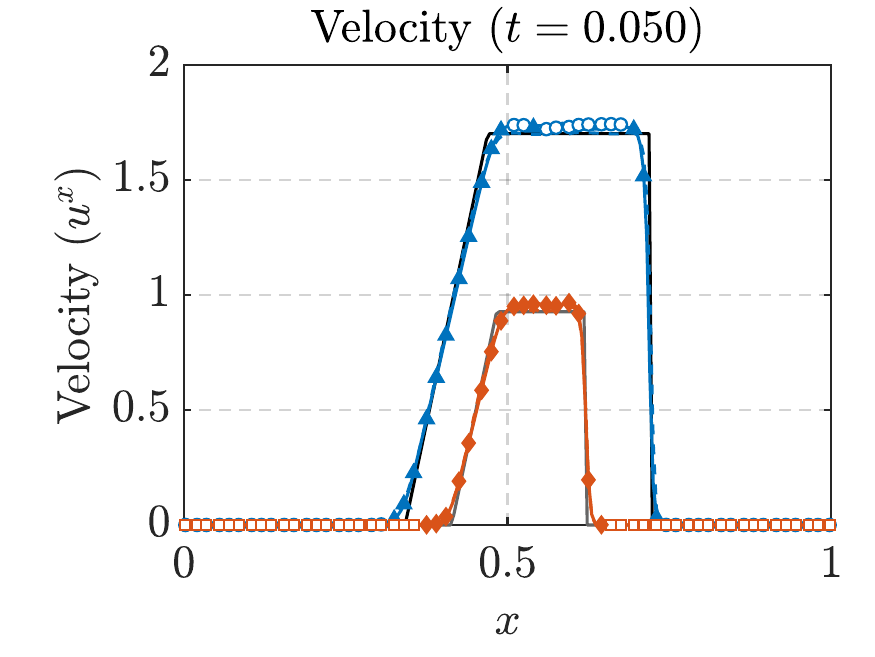}{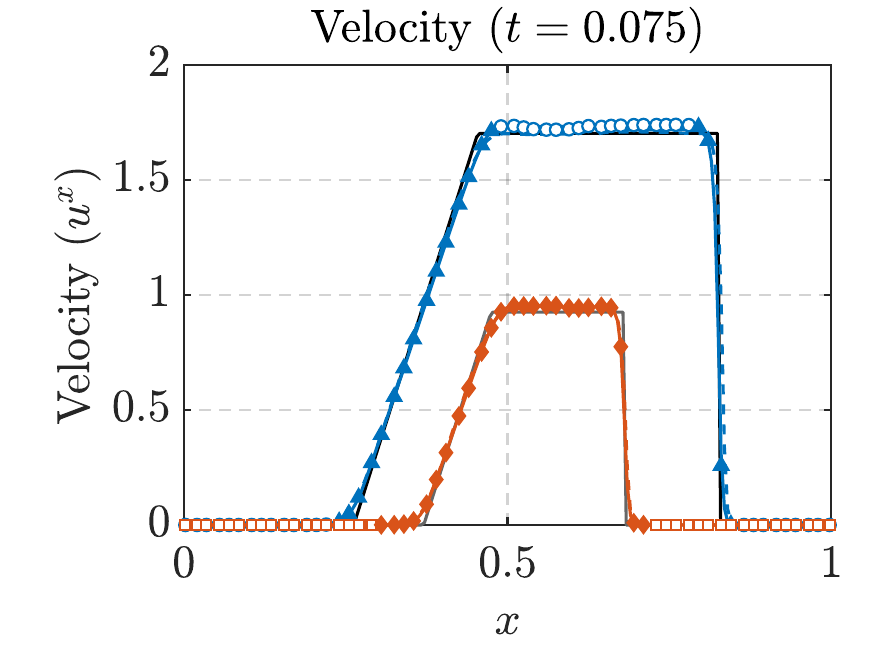}
\caption{\hyperref[Test2]{Test~2}. Mass density, temperature and velocity of both species along $x$ (projection along $y$), for the hybrid, Euler and kinetic solvers and for the analytic solutions of the uncoupled species Euler equations, used here as reference profiles. For readability, the hybrid solution is shown only in a subset of cells, roughly equally spaced in $x$; open markers denote cells in the Euler regime and filled markers cells in the kinetic regime. The kinetic solver is active in regions where non-equilibrium effects would appear.}
\label{Fig3}
\end{figure}

Figure~\ref{Fig3} shows density, temperature and velocity of both species along $x$ for the hybrid, Euler and kinetic solvers, together with the analytic solution. In the present scaling the exchange between species is not stiff and acts on time scales much longer than $t_f=0.075$, so that the analytic solution of the Euler equations for each species taken separately is a good reference, as in \cite{bis-bos-dim-gro-mar-2022}. The indicator activates the kinetic solver in regions where non-equilibrium effects would appear. When the waves move away and leave constant states behind, the cells go back to the Euler description. The three numerical solutions agree with each other and with the analytic solution, and no artefacts are visible at the interfaces. Since in this regime the three solvers give the same macroscopic profiles, we do not report errors with respect to the kinetic solution.

The computational times are 15.63~s for the Euler solver, 1403~s for the kinetic solver and 410.5~s for the hybrid one, which is therefore about 3.4 times faster than the kinetic solver.

\subsection{Test 3: flow past an inclined NACA 0012 airfoil}\label{Test3}
We consider a rarefied binary mixture flowing past a NACA 0012 airfoil with an angle of attack of $20^\circ$, a setting inspired by \cite{bos-dim-2020}. The flow around the airfoil is asymmetric and far from equilibrium, and the two species, which have different masses, behave differently.

The domain $[0,3]\times[0,3]$ is discretized with 1032 triangles, and the velocity domain $[-20,20]^2$ with $60$ points per direction. The mass ratio is $m_2/m_1=4$. The gas is initially in equilibrium and the initial moments $W_i=(\rho_i,u_i^x,u_i^y,T_i)$ are
\begin{equation}
W_1(x,y)=(m_1,\,1,\,0,\,3),\qquad W_2(x,y)=(m_2,\,1,\,0,\,3),
\end{equation}
i.e. unit number density for both species, and $f_1$, $f_2$ are the corresponding Maxwellians. Zero-gradient conditions are imposed on the external boundary and specular reflection on the airfoil. We take $\varepsilon=0.1$, $t_f=0.1$ and $\eta_1=\eta_2=\varepsilon/3$; the partition is updated at the beginning of every time step.

\begin{figure}[htbp]
\centering
\figrow{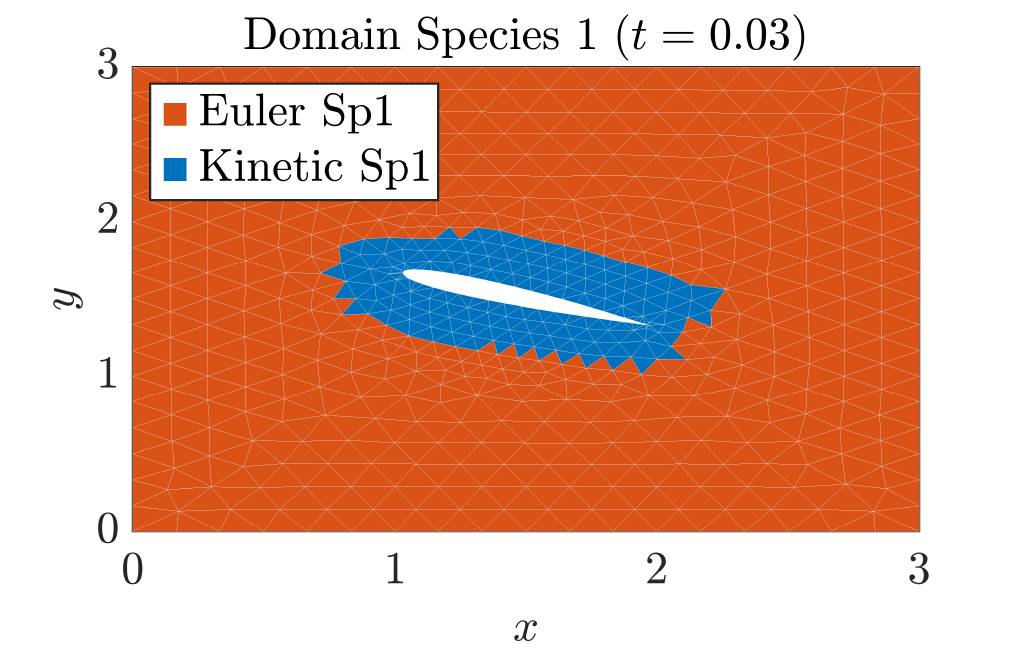}{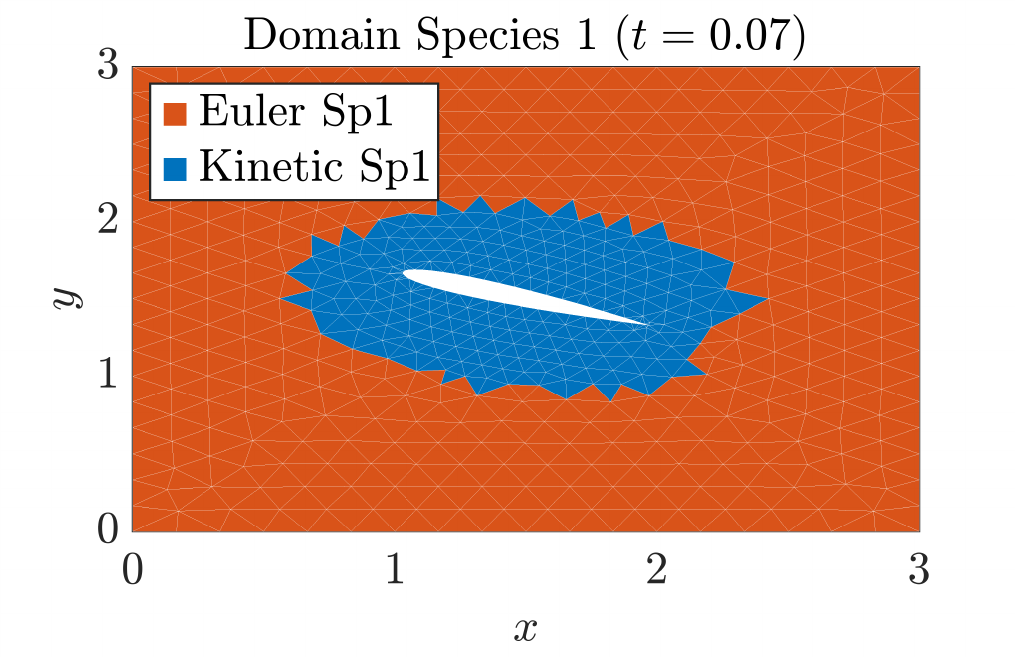}{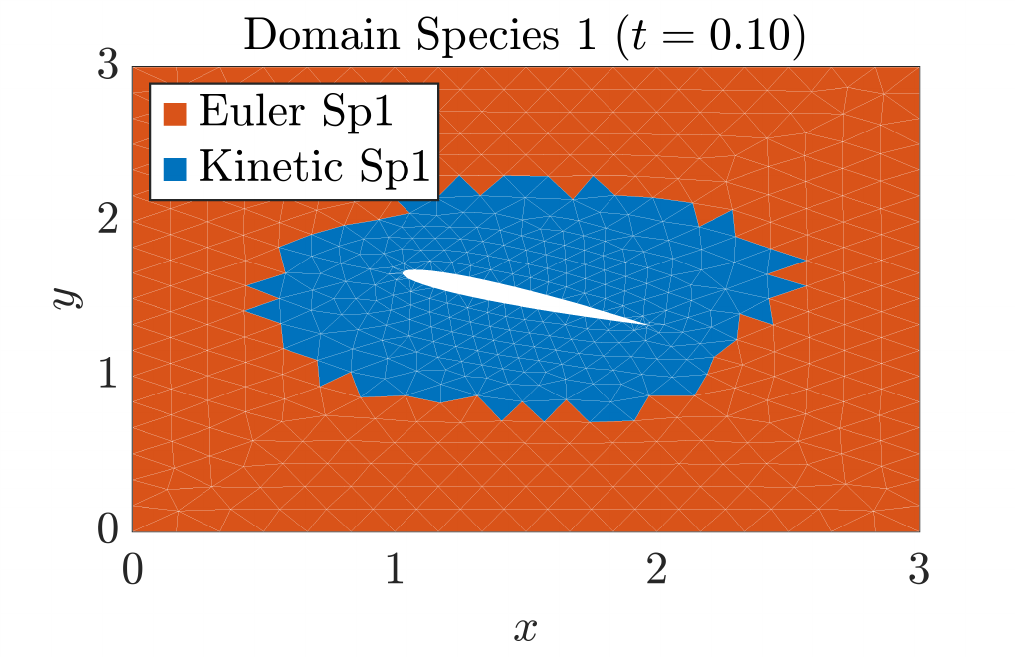}
\figrow{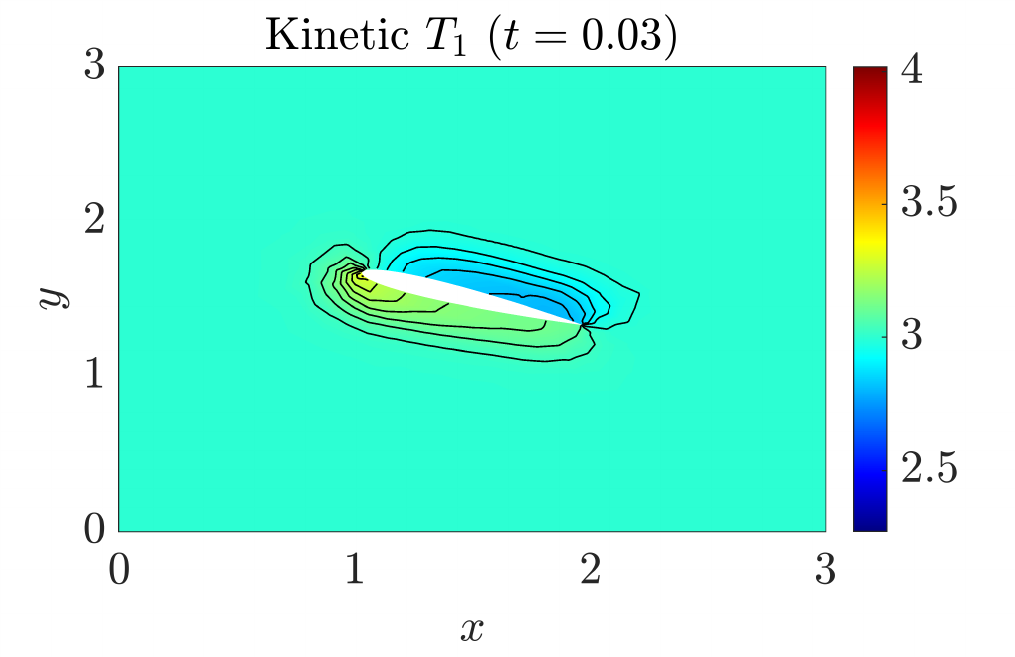}{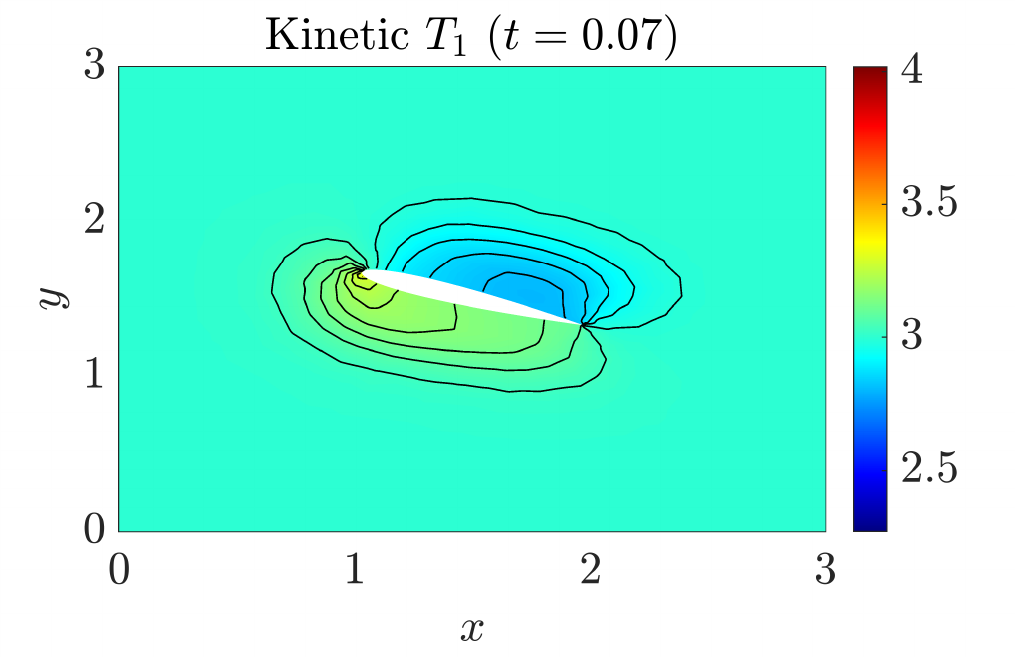}{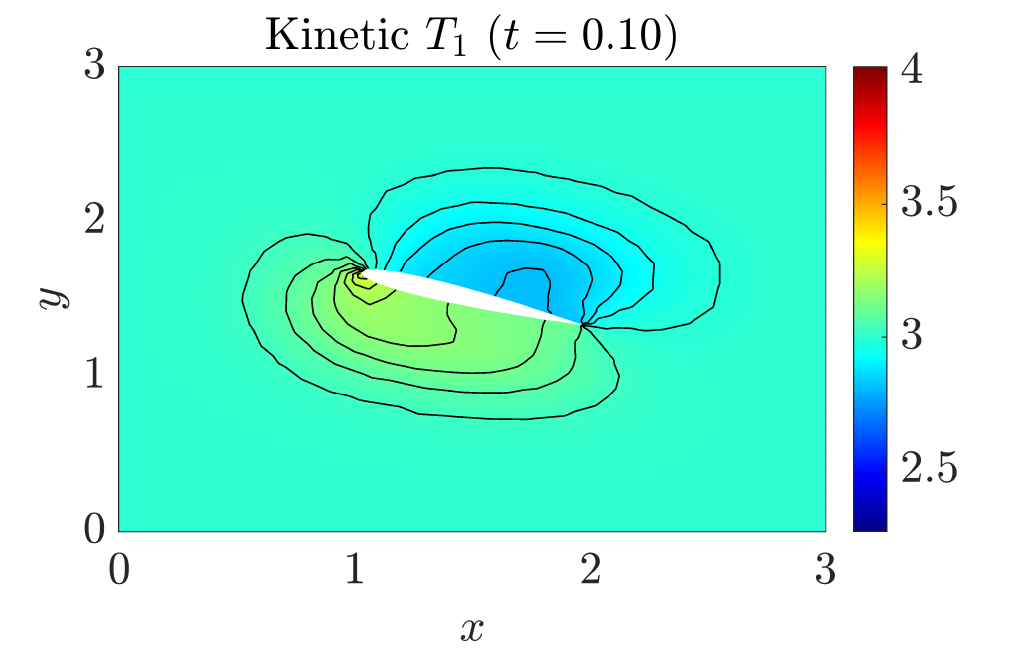}
\figrow{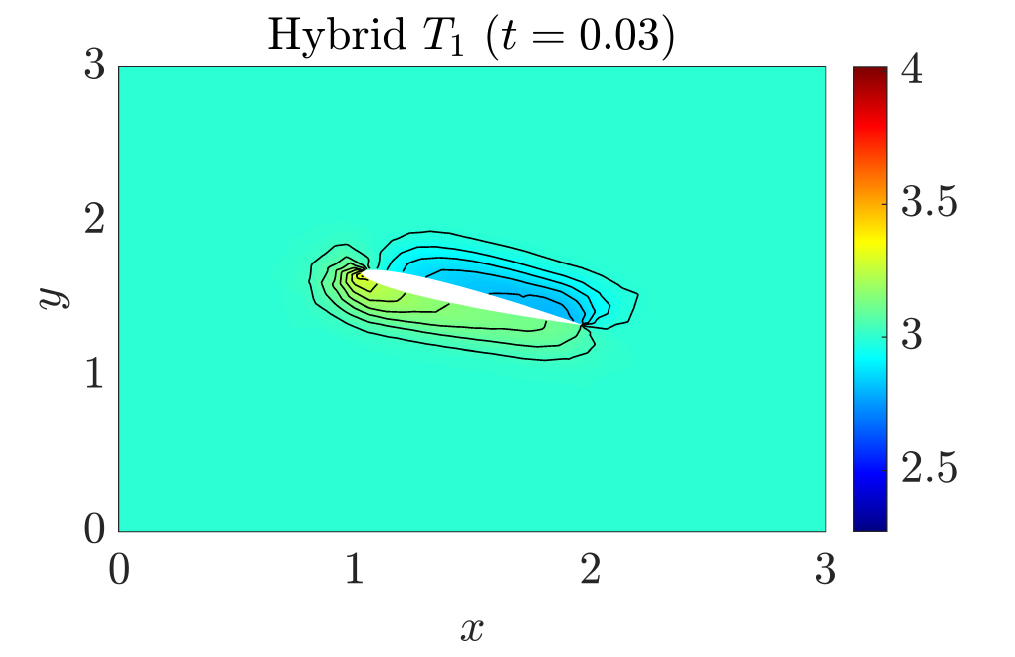}{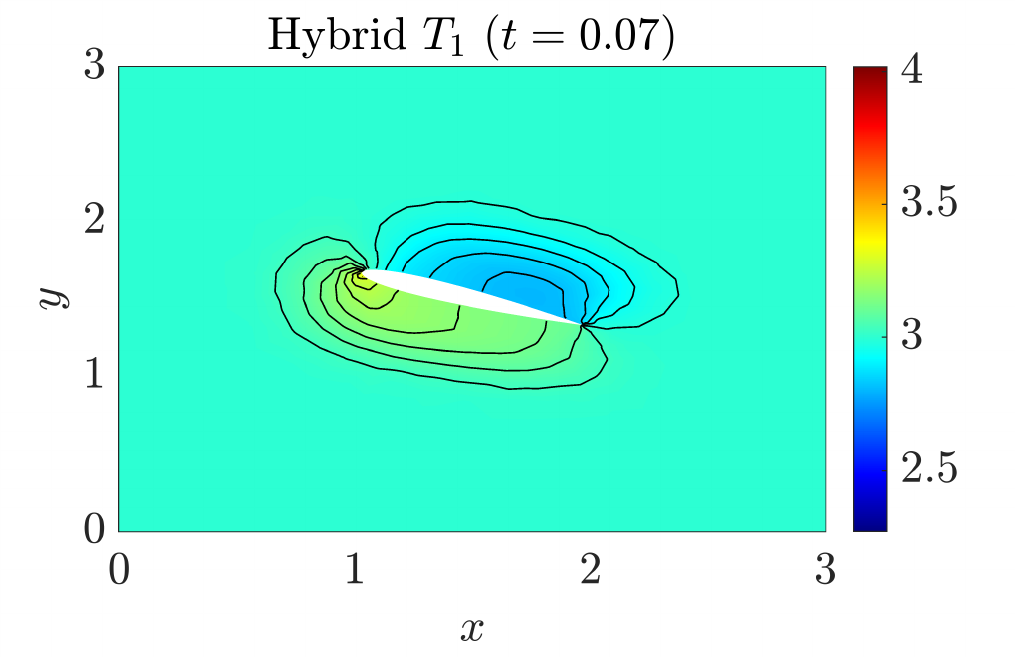}{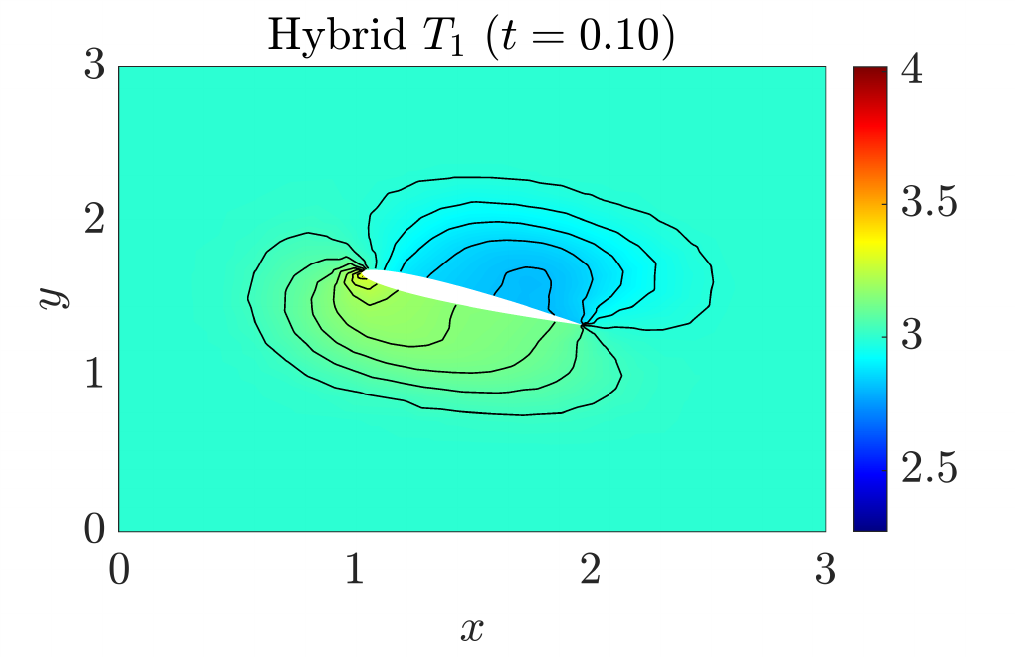}
\figrow{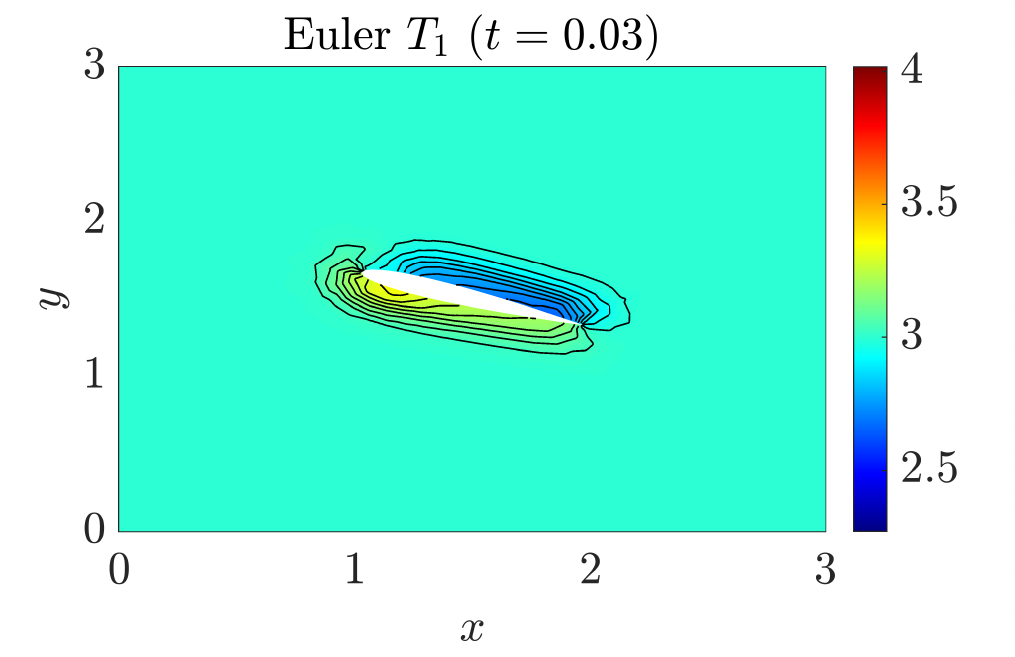}{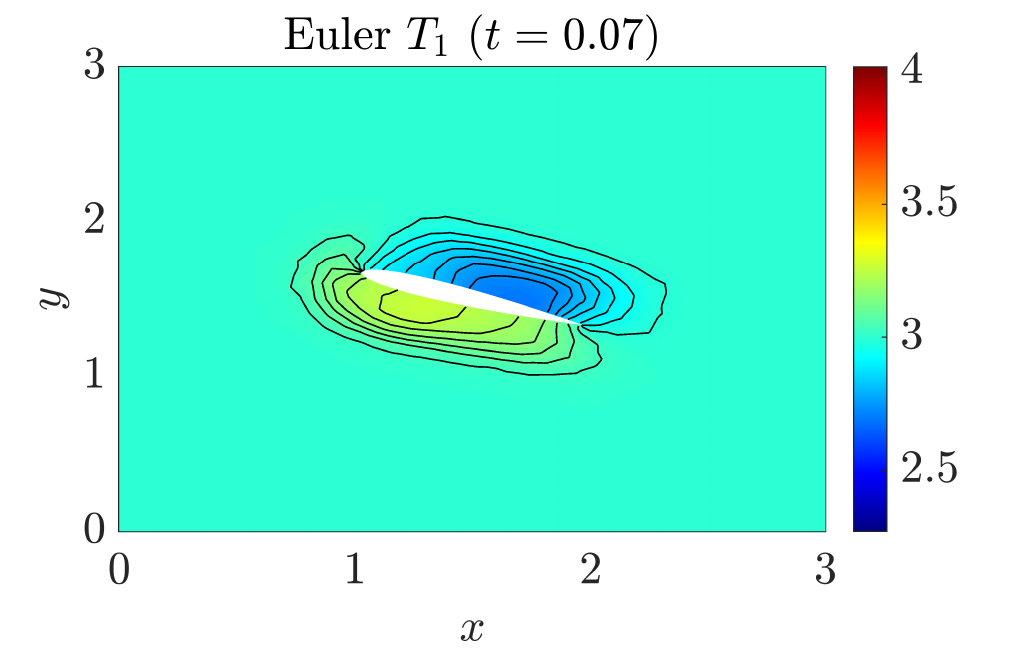}{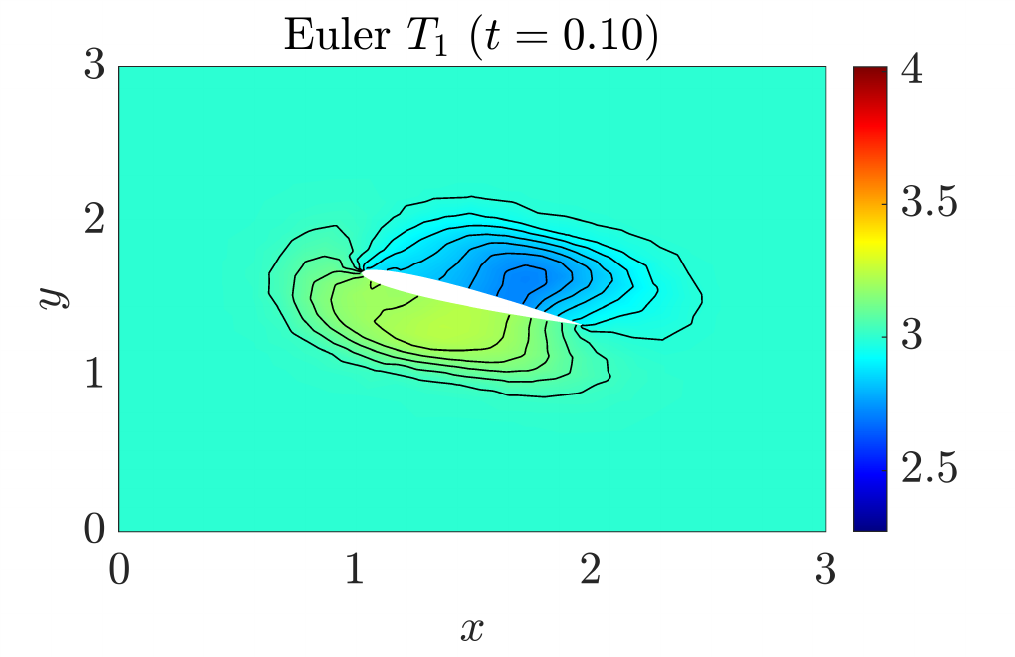}
\caption{\hyperref[Test3]{Test~3}, species 1. First row: partition of the domain. Second to fourth row: temperature computed with the kinetic, hybrid and Euler solvers.}
\label{Fig7-Temperature-Species1}
\end{figure}

\begin{figure}[htbp]
\centering
\figrow{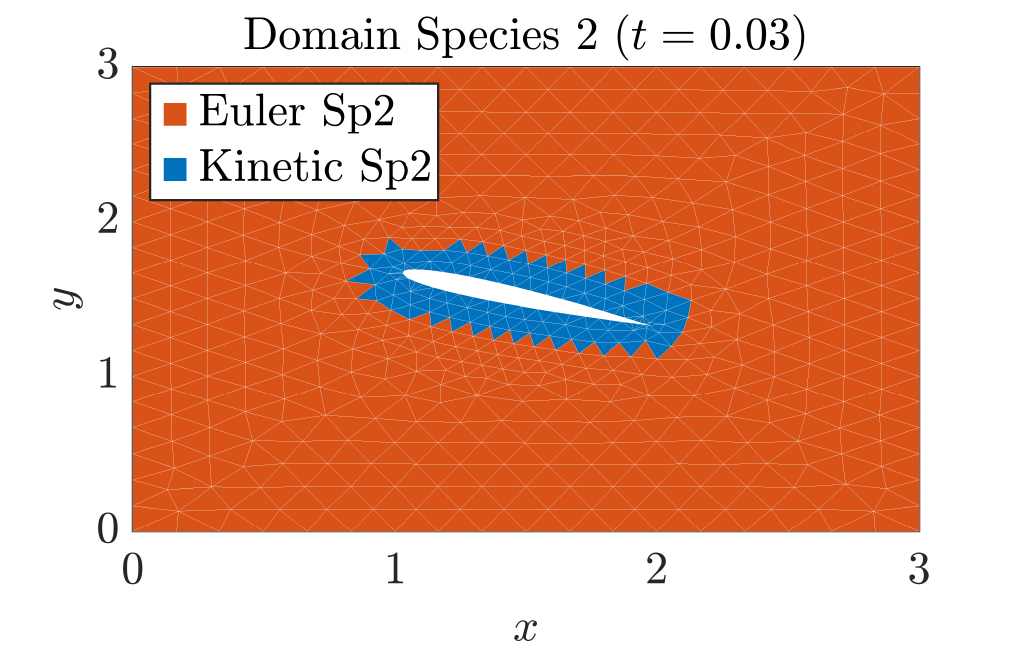}{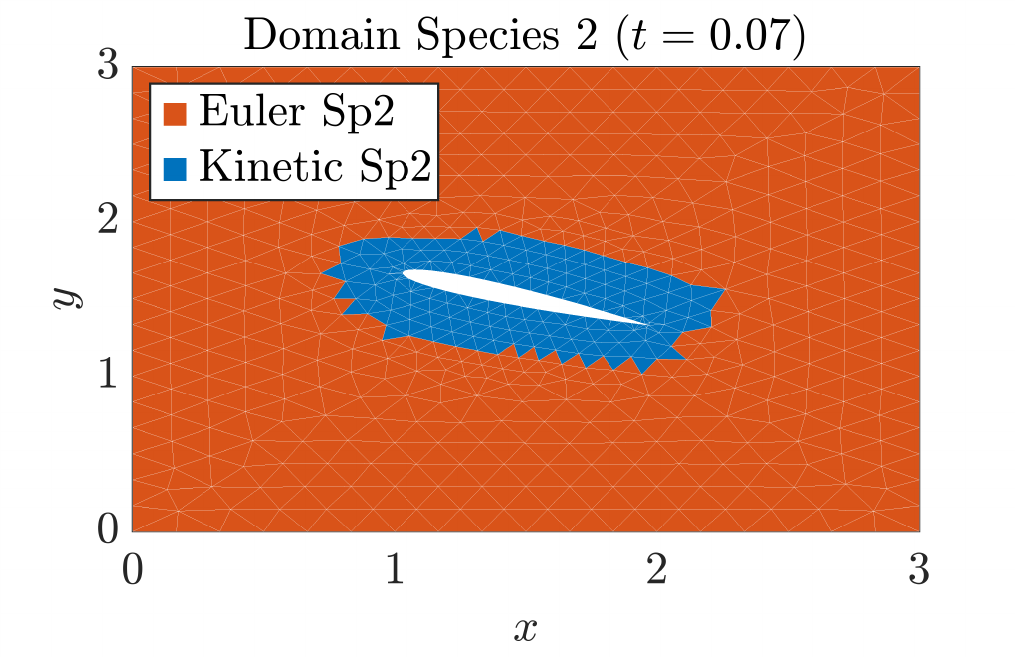}{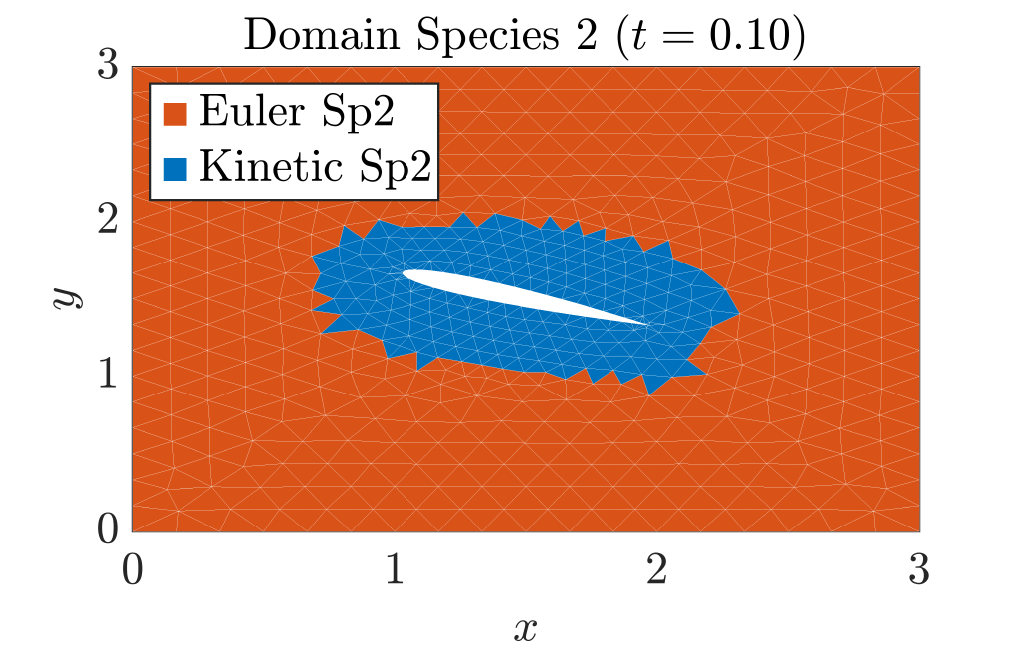}
\figrow{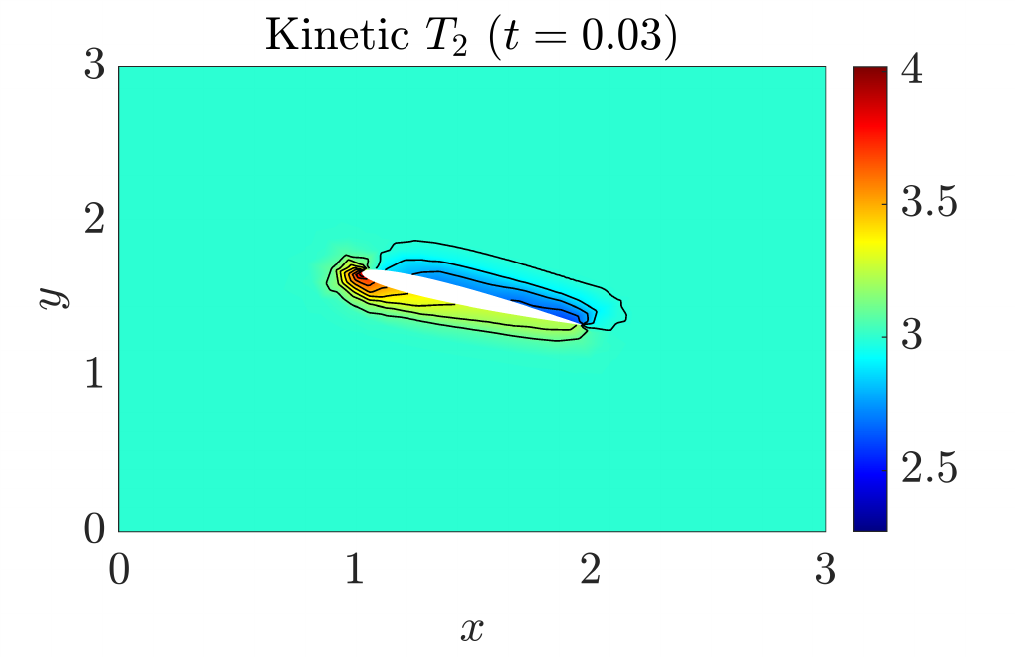}{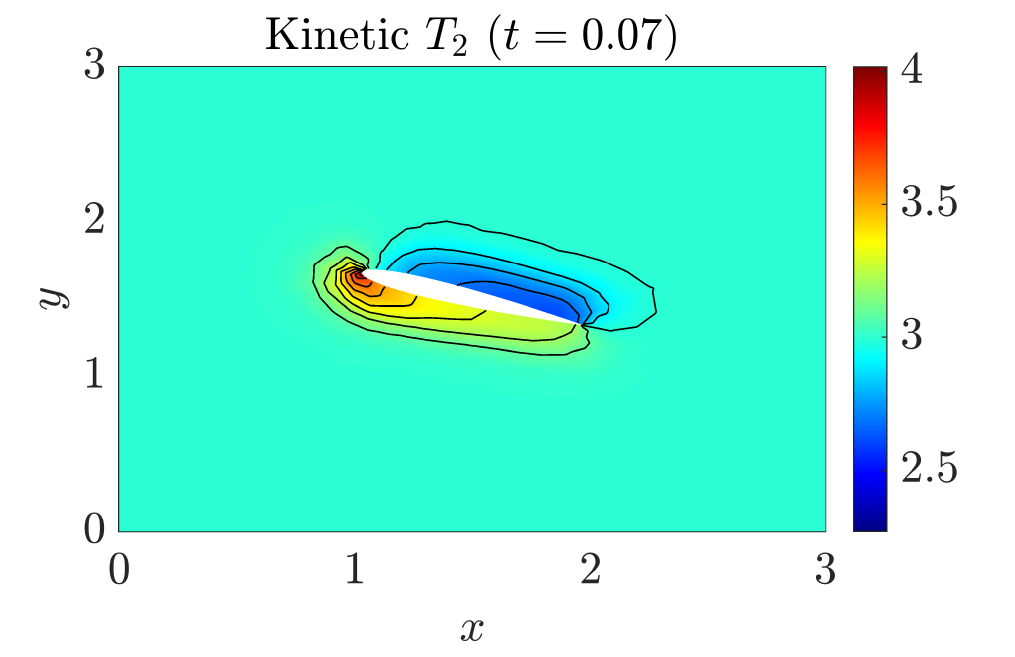}{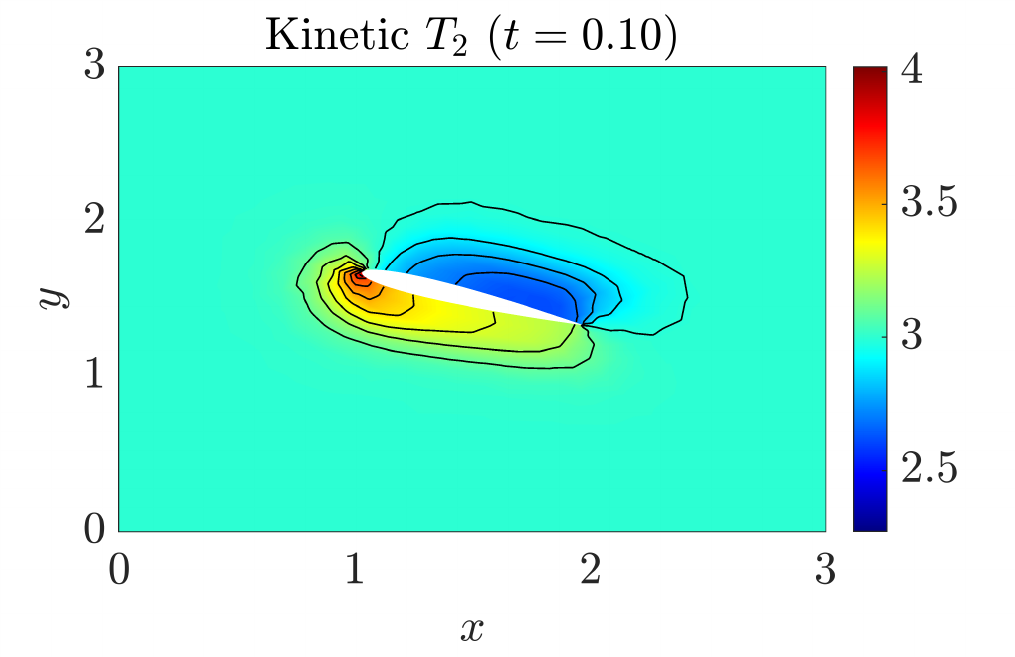}
\figrow{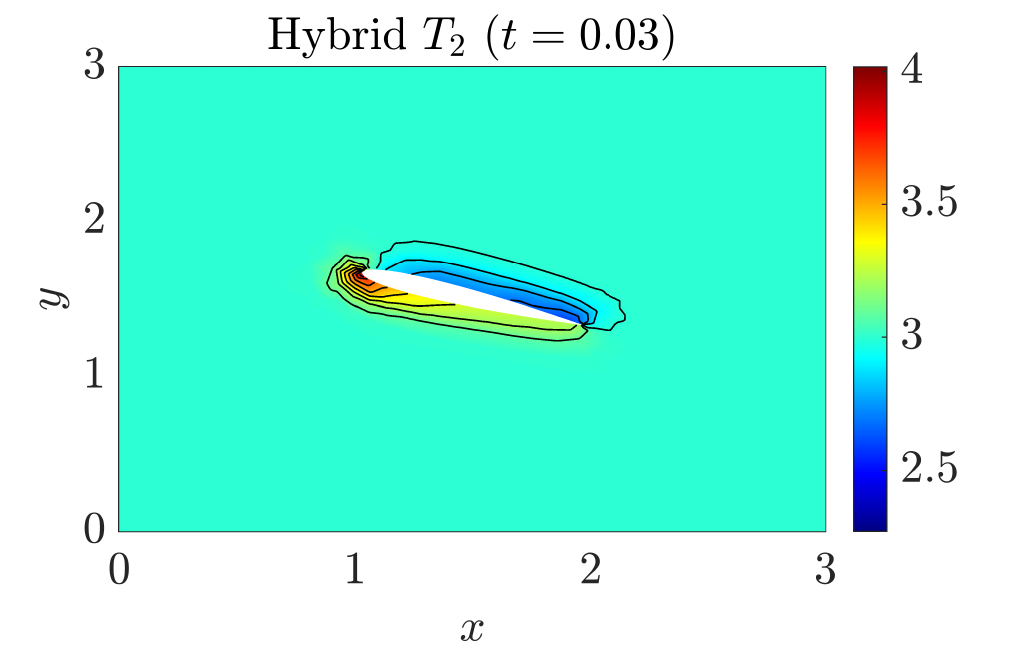}{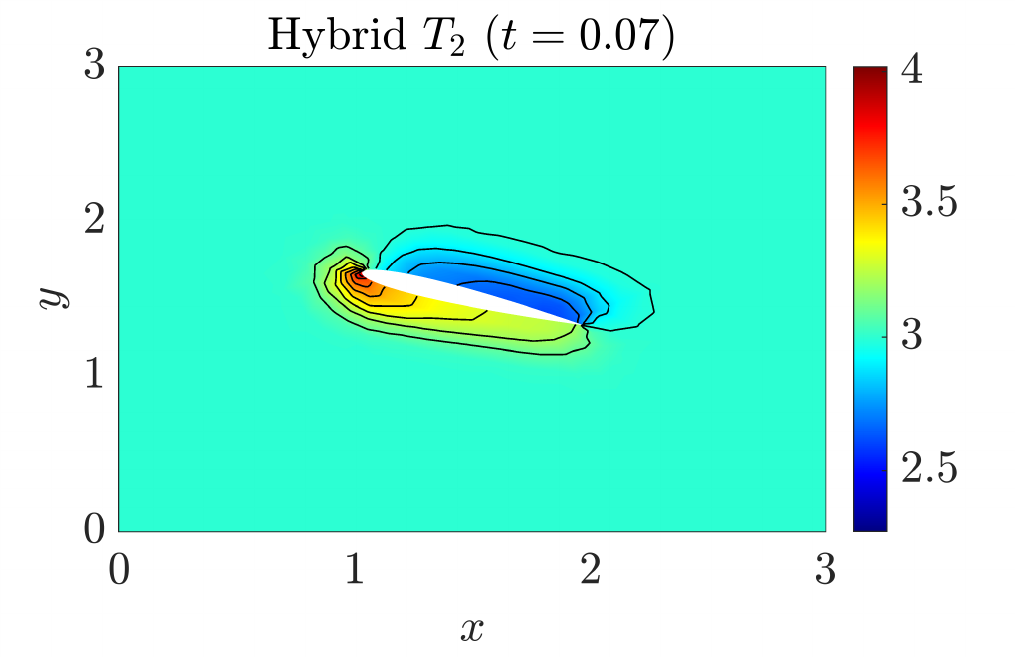}{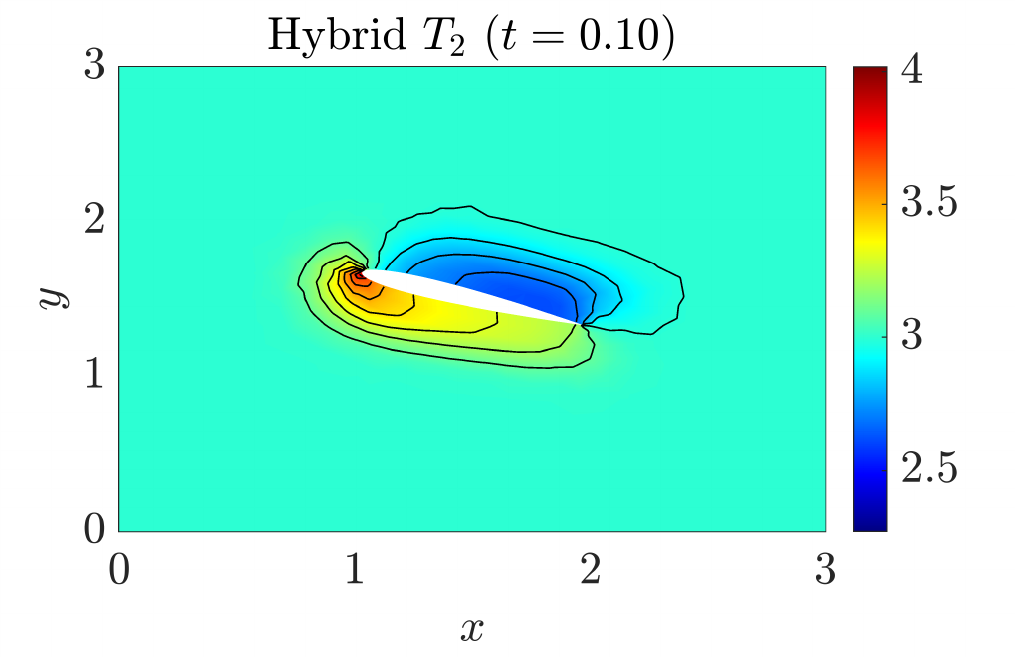}
\figrow{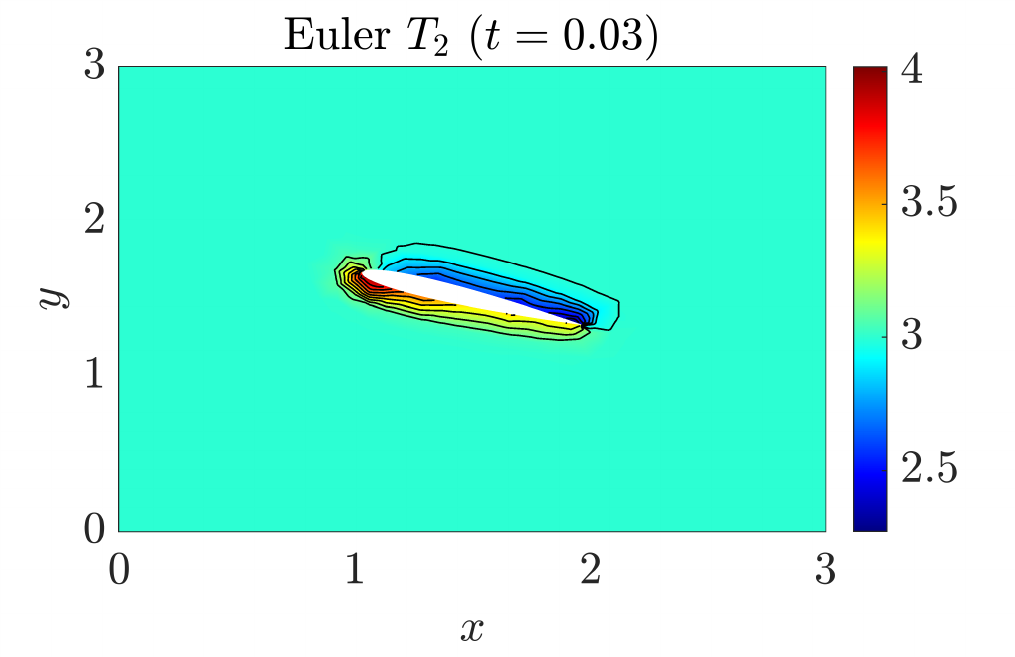}{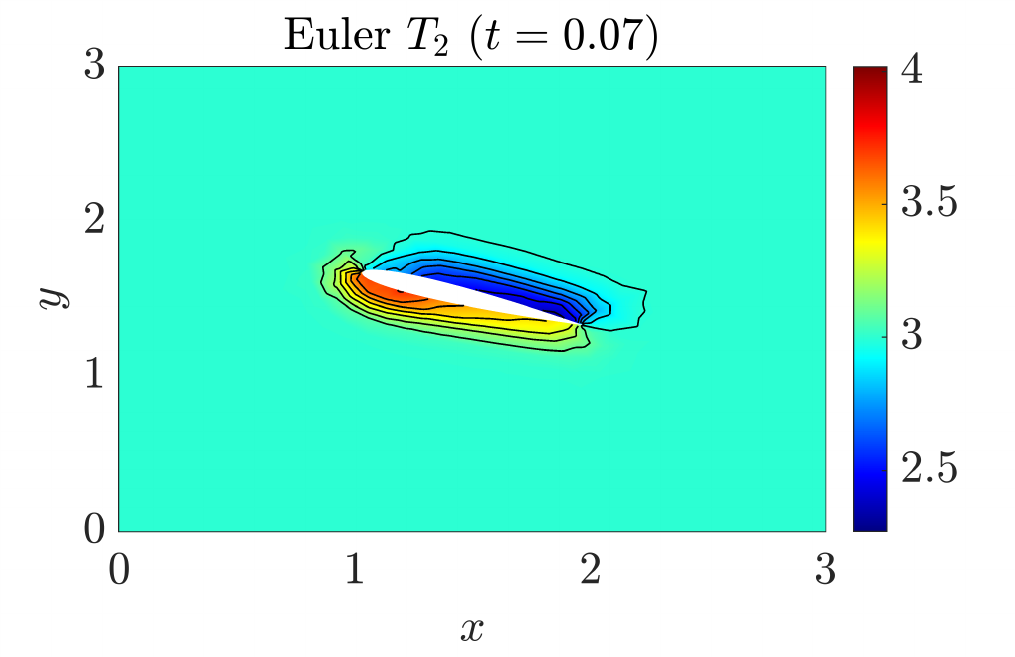}{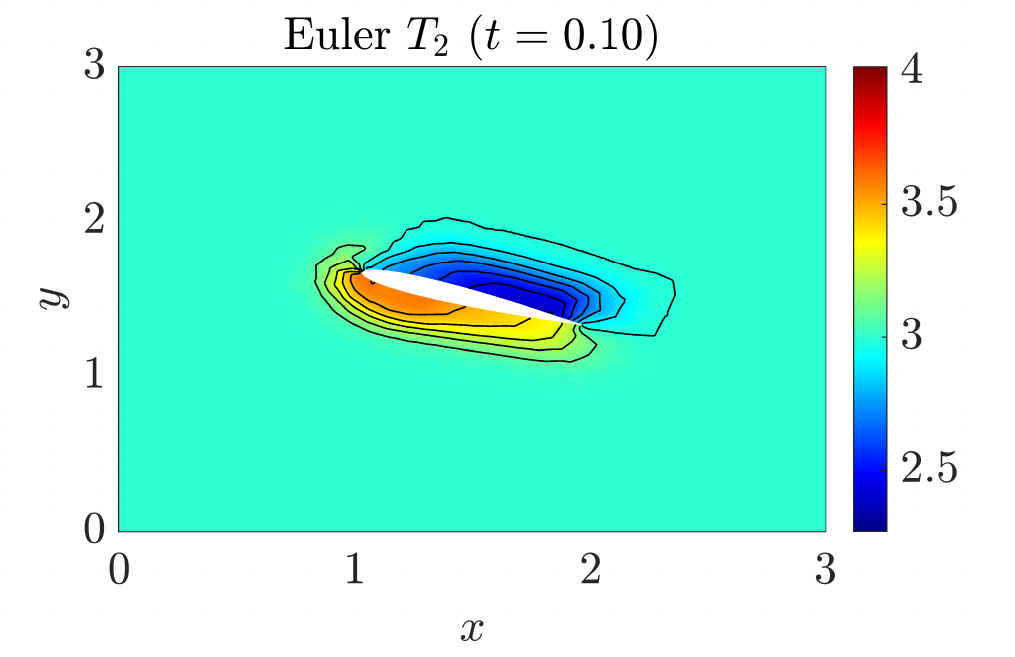}
\caption{\hyperref[Test3]{Test~3}, species 2. First row: partition of the domain. Second to fourth row: temperature computed with the kinetic, hybrid and Euler solvers.}
\label{Fig7-Temperature-Species2}
\end{figure}

The indicator gives a different partition for each of the two species. Figures~\ref{Fig7-Temperature-Species1} and~\ref{Fig7-Temperature-Species2} show the partition and the temperature of species 1 and 2 for the hybrid, kinetic and Euler solvers; the same colour scale is used in all panels of both figures, while the contour lines are equally spaced between the minimum and the maximum temperature of each species, so that the structure of the solution is visible for both species. The hybrid solution is very close to the kinetic one, and the kinetic region stays around the airfoil. The Euler solution is clearly different.

\begin{figure}[htbp]
\centering
\figrow{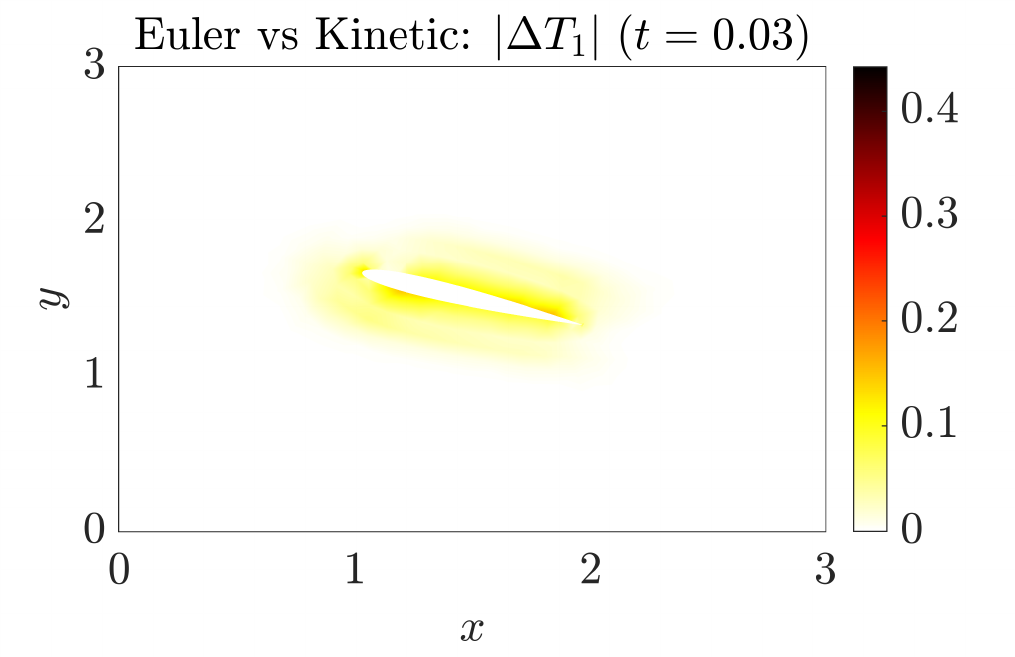}{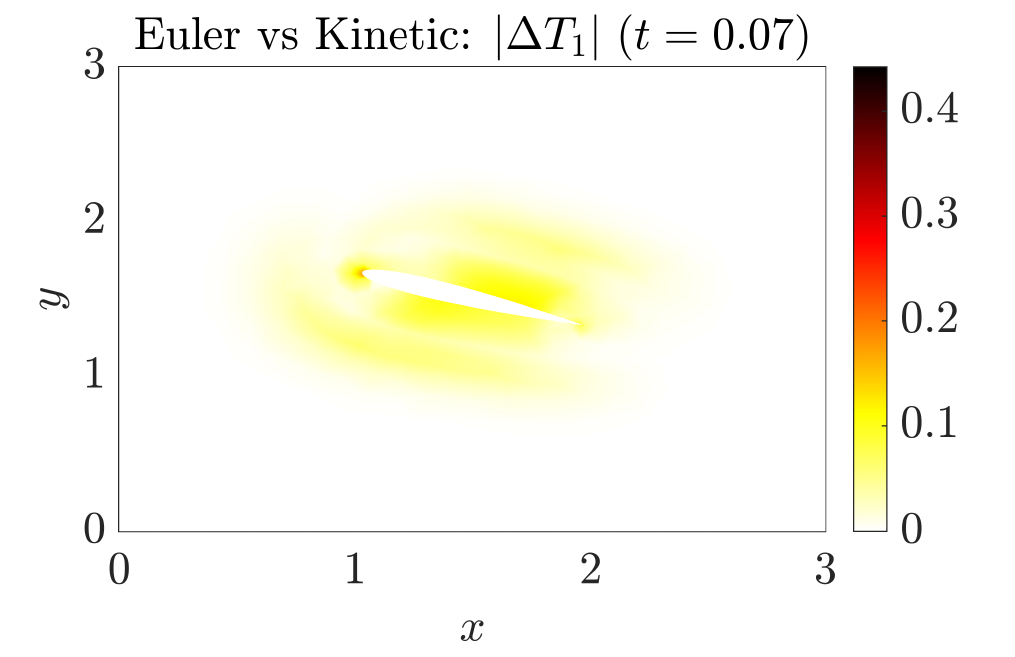}{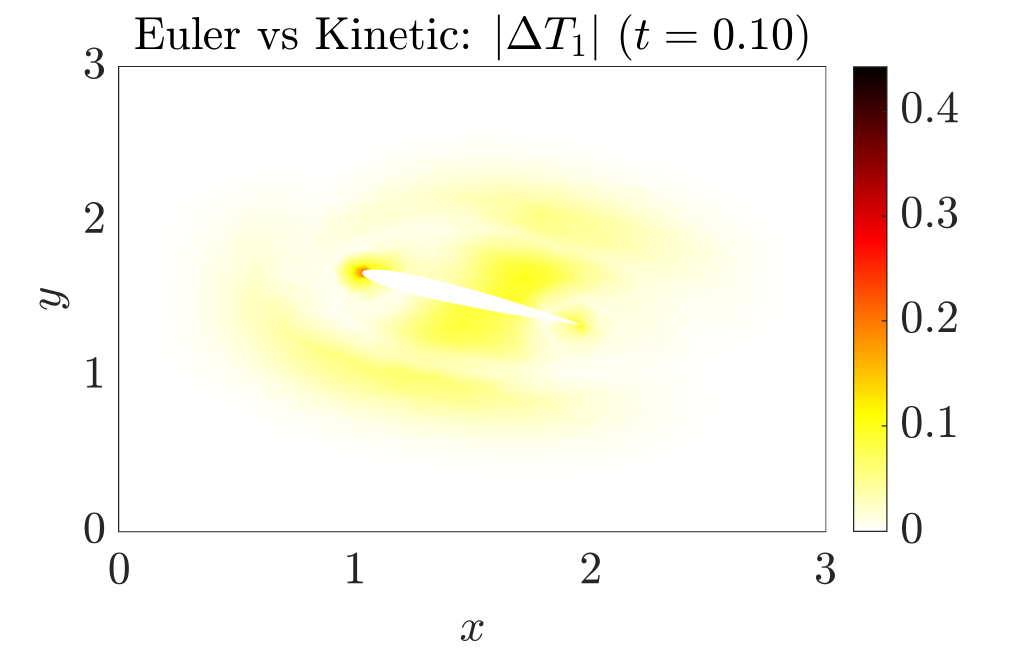}
\figrow{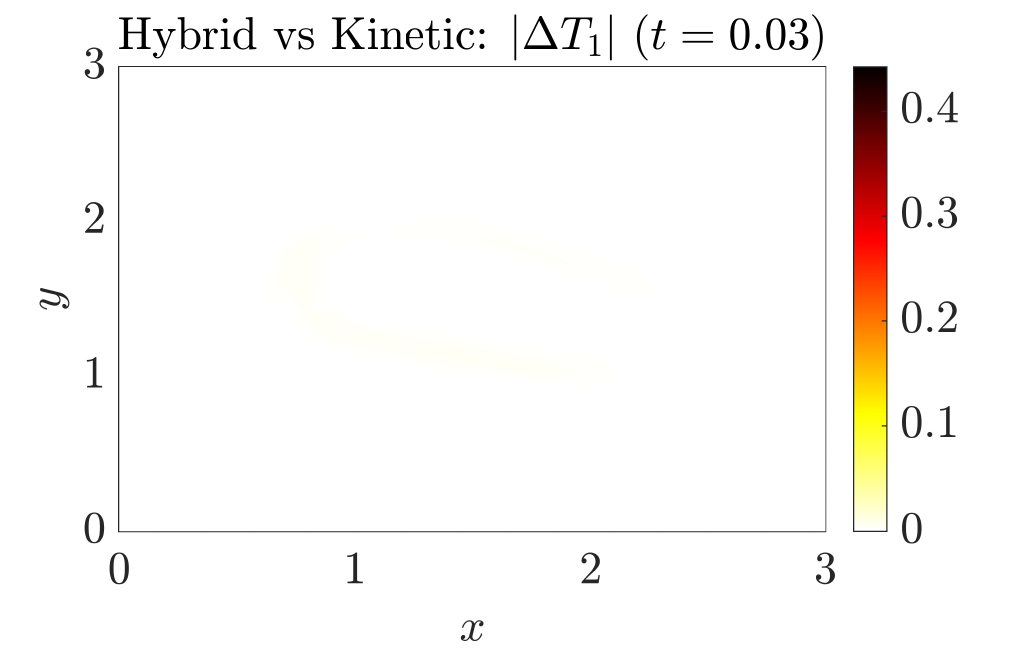}{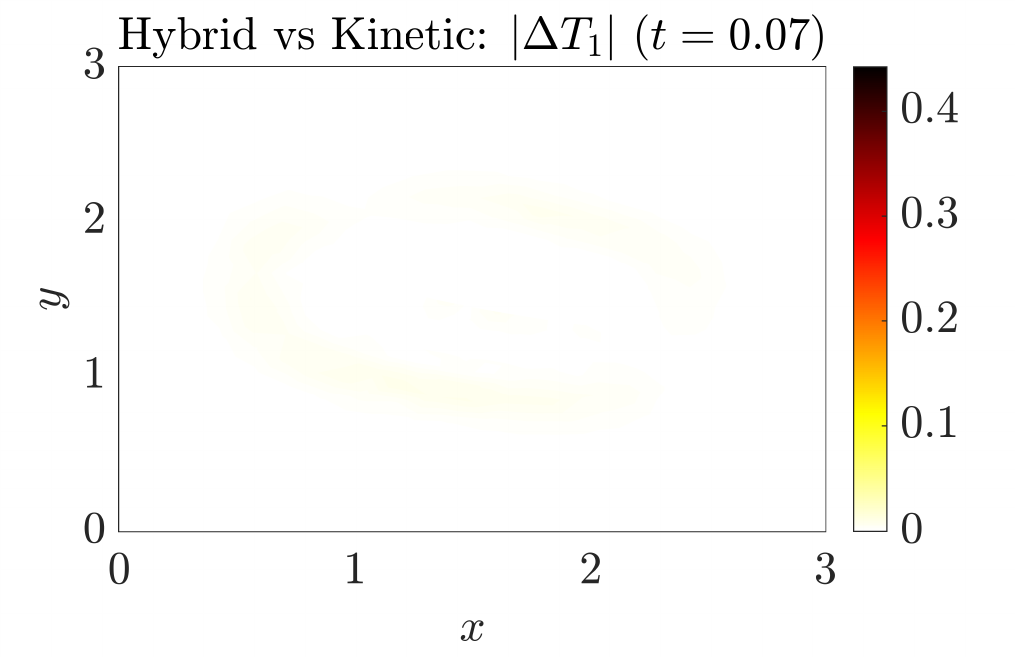}{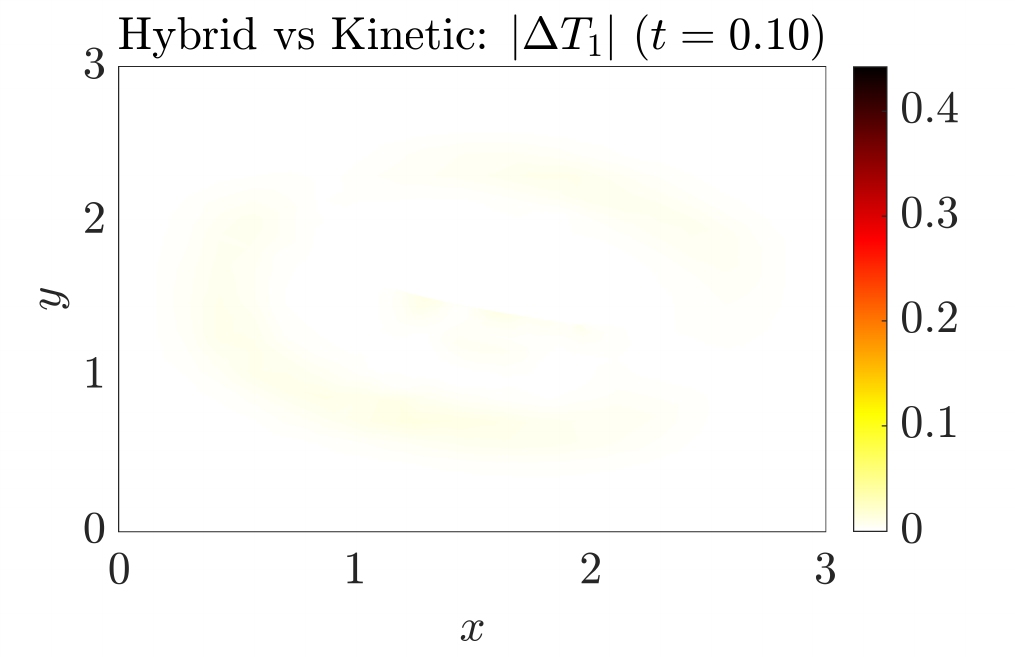}
\figrow{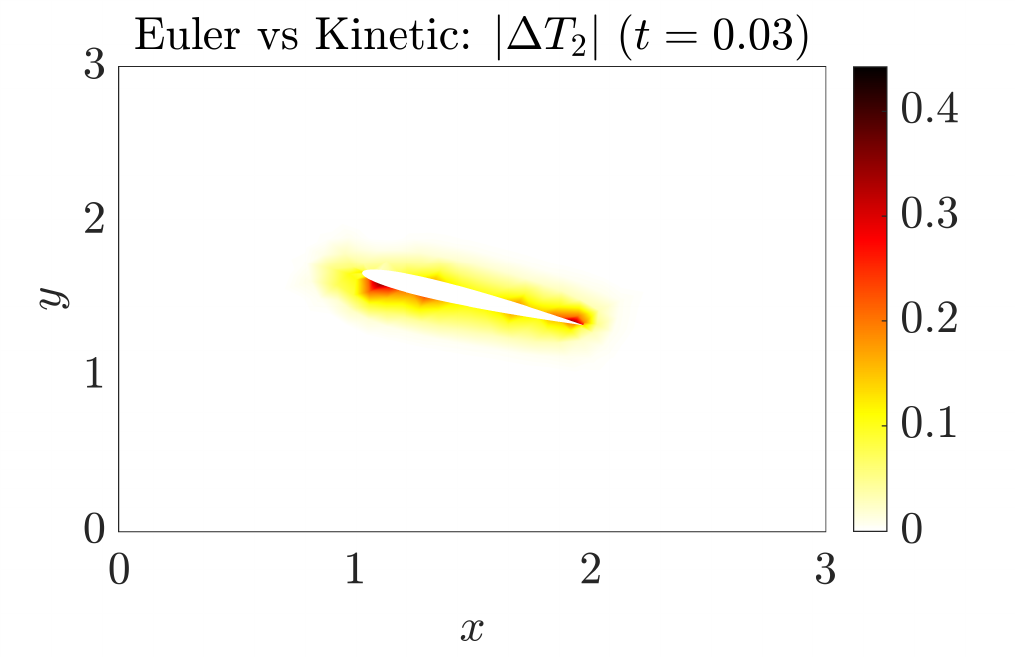}{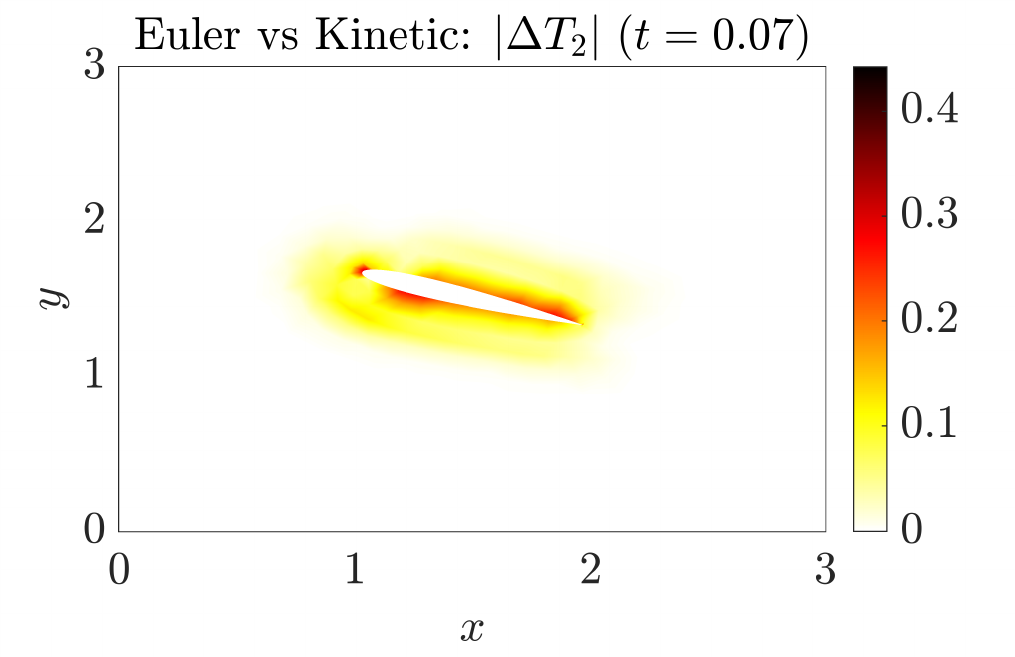}{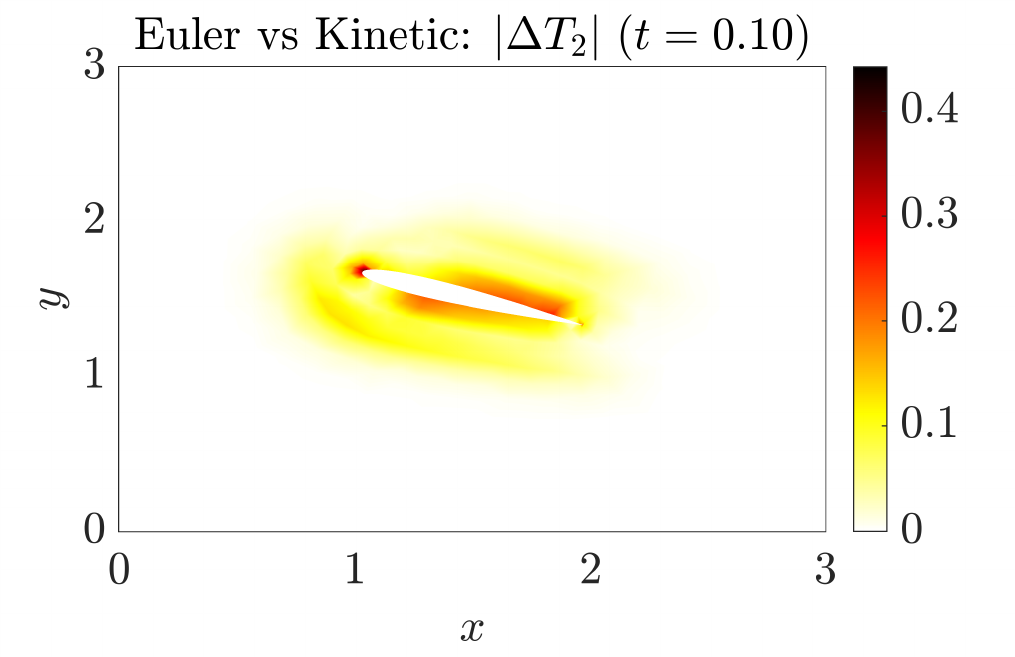}
\figrow{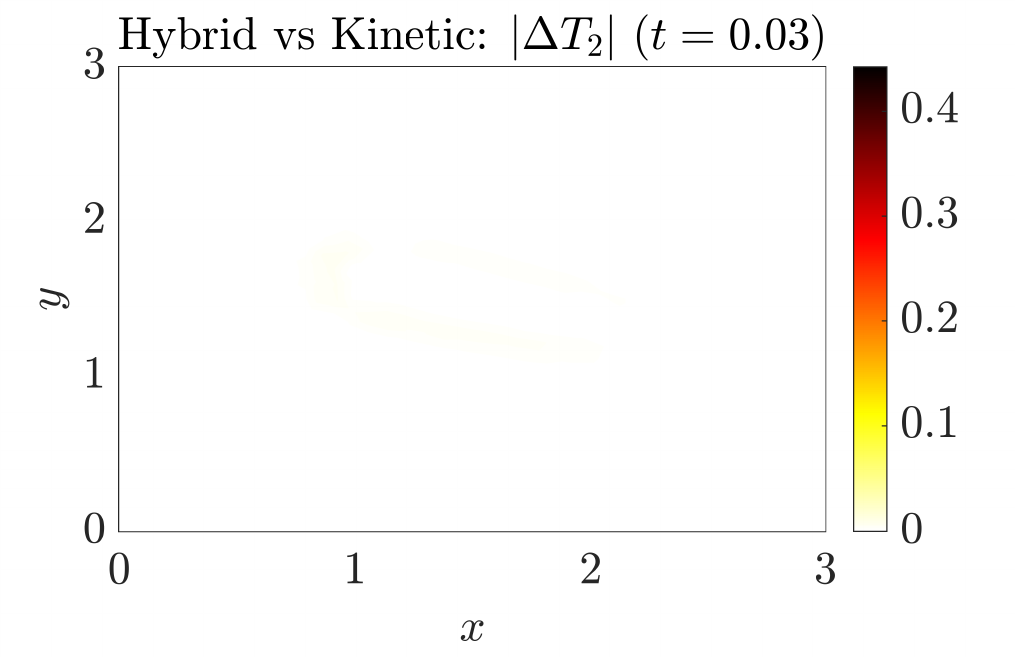}{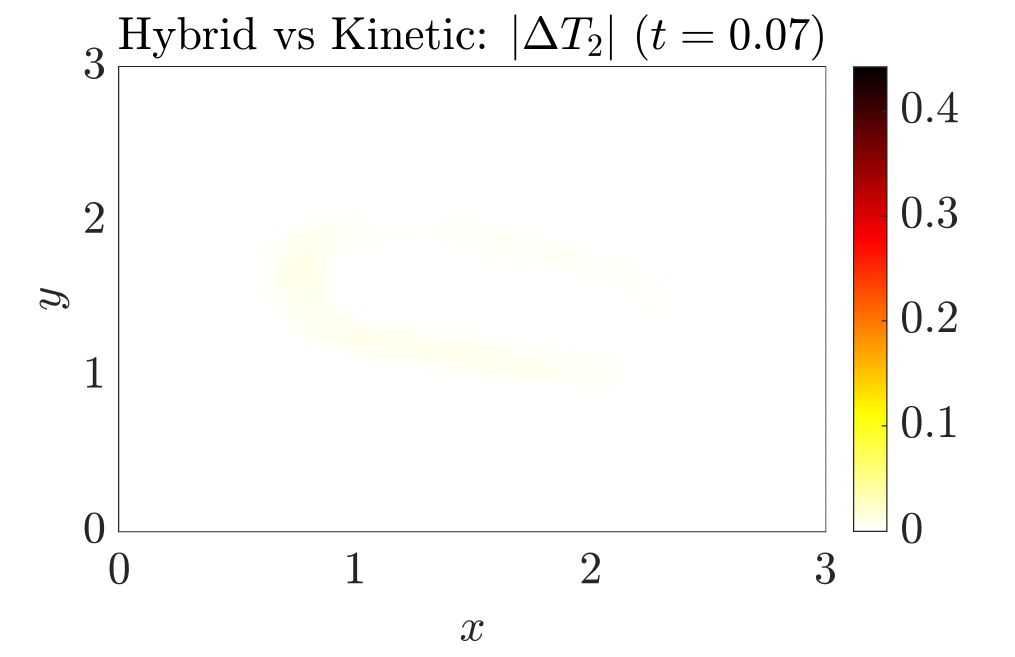}{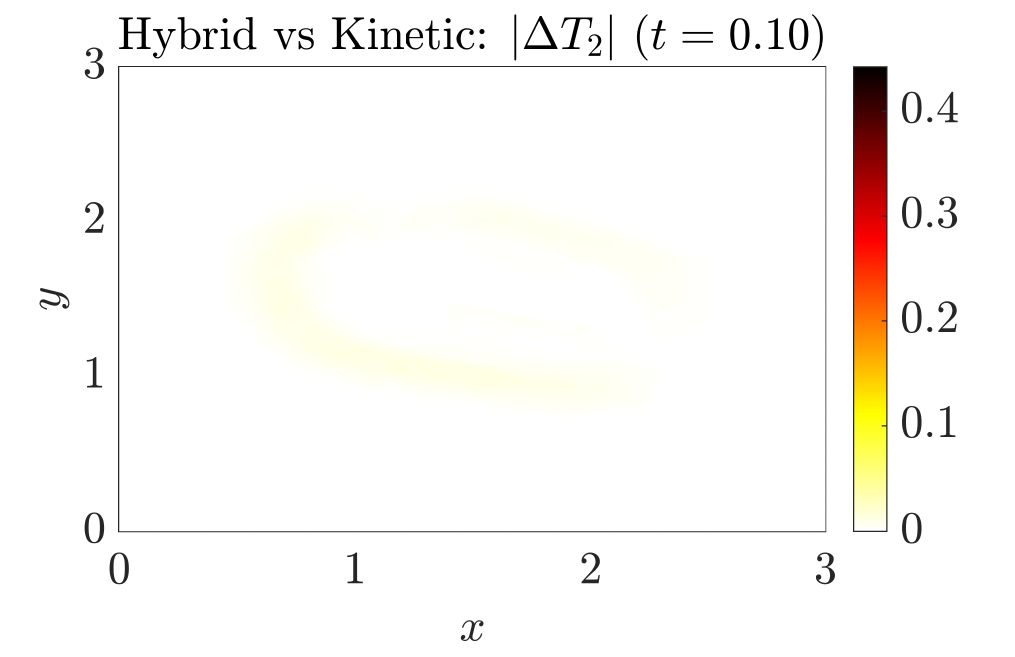}
\caption{\hyperref[Test3]{Test~3}. Absolute difference between the temperature computed with the Euler or hybrid solver and the kinetic one. First and second row: species 1, Euler and hybrid solver. Third and fourth row: species 2, Euler and hybrid solver. The same linear colour scale is used in all panels.}
\label{Fig7-Temperature-Difference}
\end{figure}

Figure~\ref{Fig7-Temperature-Difference} shows the pointwise difference $|T_i-T_i^{\mathrm{kin}}|$, $i=1,2$, taking the kinetic solution as reference. The difference is clearly visible for the Euler solver, while for the hybrid solver it remains small everywhere.

The computational times are 4475~s for the kinetic solver, 981.7~s for the hybrid one (speedup about 4.6) and 136.9~s for the Euler solver. As a measure of the distance from the kinetic solution we use the $L^1$ error on the temperature of each species, averaged in space and time,
\begin{equation}\label{eq:errT}
E(T_i)=\frac{1}{t_f\,|\Omega|}\sum_{n=1}^{N_t}\Delta t^n\sum_{c=1}^{N_P}|P_c|\,\big|T_{i,c}^n-T_{i,c}^{n,\mathrm{kin}}\big|,
\end{equation}
where $t_f$ is the final time, $N_t$ the number of time steps, $\Delta t^n=t^n-t^{n-1}$ the $n$-th time step, $P_c$, $c=1,\dots,N_P$, the triangles of the mesh, with area $|P_c|$, and $|\Omega|=\sum_c|P_c|$ the area of the computational domain; $T_{i,c}^n$ and $T_{i,c}^{n,\mathrm{kin}}$ are the temperatures of species $i$ in the cell $P_c$ at time $t^n$, computed with the solver under consideration (Euler or hybrid) and with the kinetic solver, respectively.
For the Euler solver we obtain $E(T_1)=5.3\times10^{-3}$ and $E(T_2)=5.2\times10^{-3}$, and for the hybrid solver $E(T_1)=7.5\times10^{-4}$ and $E(T_2)=4.4\times10^{-4}$, that is, about seven and twelve times smaller, respectively.

\subsection{Test 4: spatially varying Knudsen number}\label{Test4}
Following \cite{cap-par-rey-2025, fil-rey-2014}, we consider a Knudsen number that varies in space. The domain is $[-L_x,L_x]\times[0,L_y]$, with $L_x=1$ and $L_y=0.112$, discretized with 154 triangles. The mass ratio is $m_2/m_1=2$. The initial data are not at equilibrium:
\begin{equation}
f_i(0,x,v)=\frac{1}{2}\Big(\mathcal M_i\big[\rho_i(x),\,\bar u_i,\,T_i(x)\big]+\mathcal M_i\big[\rho_i(x),\,-\bar u_i,\,T_i(x)\big]\Big),\qquad i=1,2,
\end{equation}
with $\bar u_1=(1,0)$, $\bar u_2=(0.5,0)$ and
\begin{equation}
\begin{aligned}
\rho_1(x)&=1+0.1\cos(\pi x), &\qquad \rho_2(x)&=1.5\,\rho_1(x),\\
T_1(x)&=2.1\big(1+0.1\cos(\pi x)\big), &\qquad T_2(x)&=1.9\big(1+0.1\cos(\pi x)\big),
\end{aligned}
\end{equation} 
where $\mathcal M_i[\rho_i,u_i,T_i]=n_i\big(\frac{m_i}{2\pi T_i}\big)^{d_v/2}\exp\big(-\frac{m_i}{2T_i}|v-u_i|^2\big)$, $n_i=\rho_i/m_i$, is the Maxwellian in velocity dimension $d_v=2$. Here $T_i(x)$ denotes the temperature of each Maxwellian component, whereas the temperature moment of the bimodal initial distribution is $T_i^{\mathrm{init}}(x)=T_i(x)+m_i|\bar u_i|^2/d_v$, with zero mean velocity. The Knudsen number is
\begin{equation}
\varepsilon(x)=10^{-5}+0.5\big(\arctan(1+200x)+\arctan(1-200x)\big),
\end{equation}
so that the gas is rarefied only in the central part of the domain. The final time is $t_f=0.05$, the velocity domain $[-16,16]^2$ is discretized with $200$ points per direction, and the thresholds are $\eta_1=5\times10^{-5}$, $\eta_2=10^{-5}$. Zero-gradient conditions are imposed on the boundary. 
Since the initial state is not at equilibrium, the first time step is computed with the kinetic model everywhere; afterwards the partition is updated automatically.

\begin{figure}[htbp]
\centering
\figrow{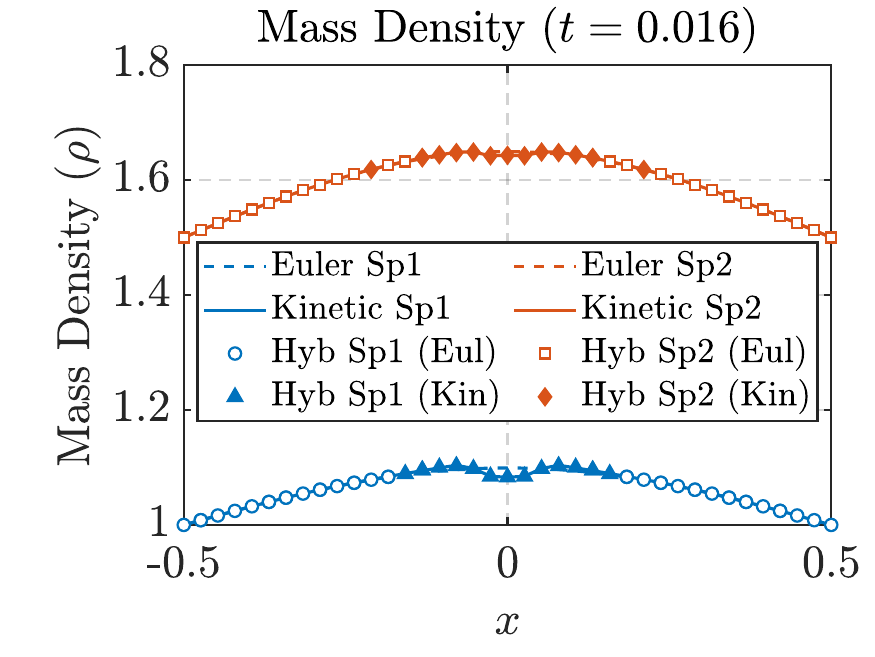}{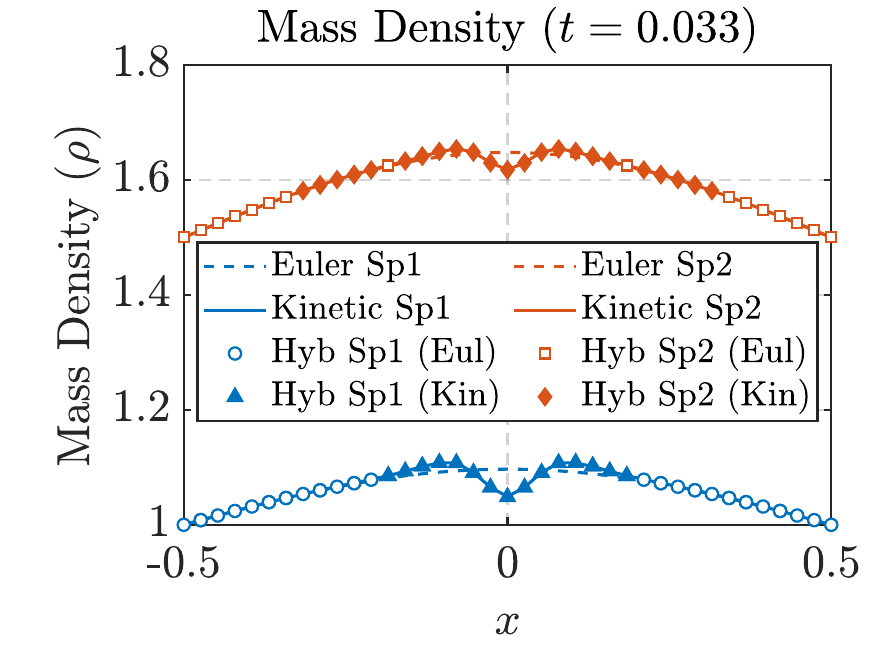}{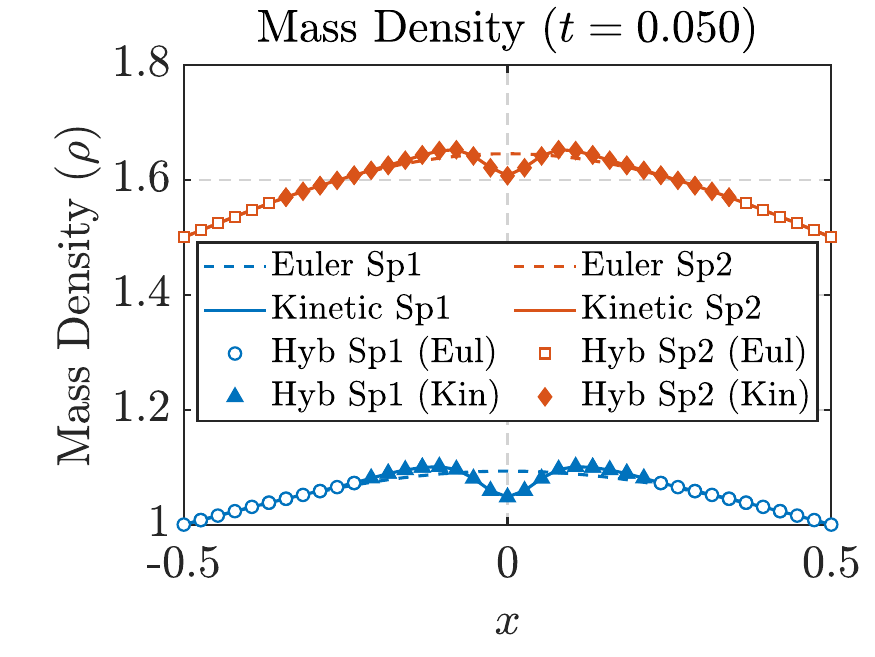}
\figrow{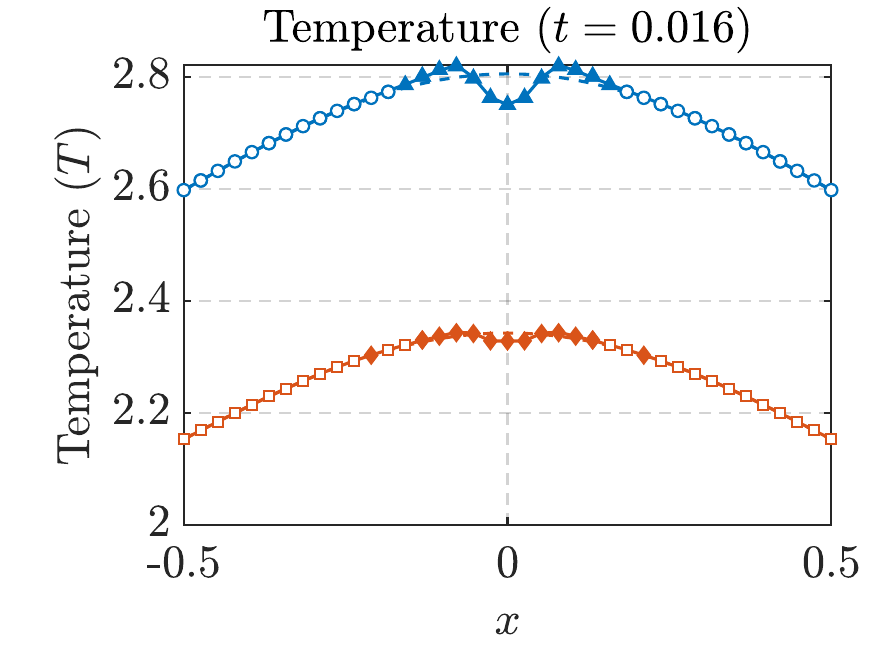}{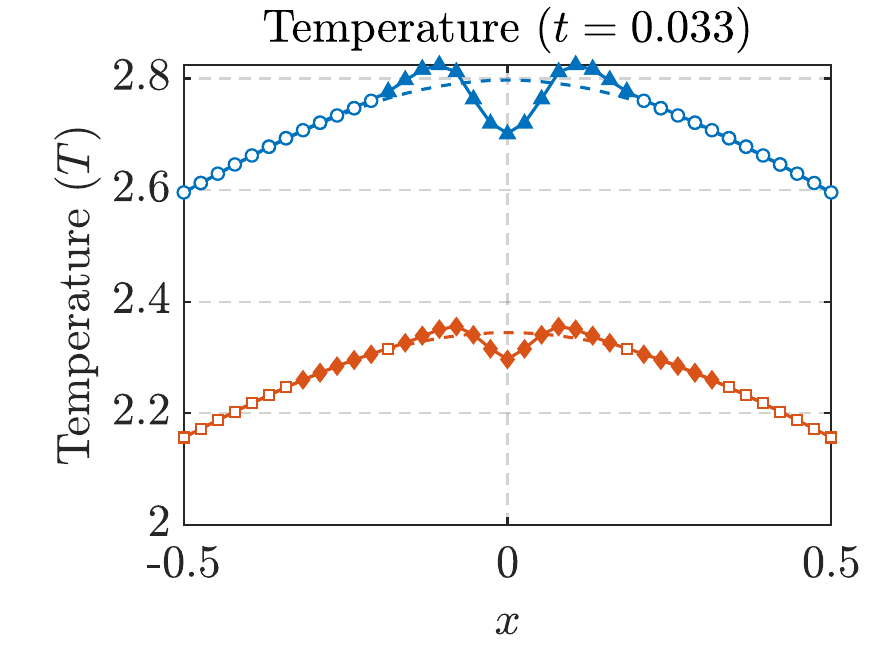}{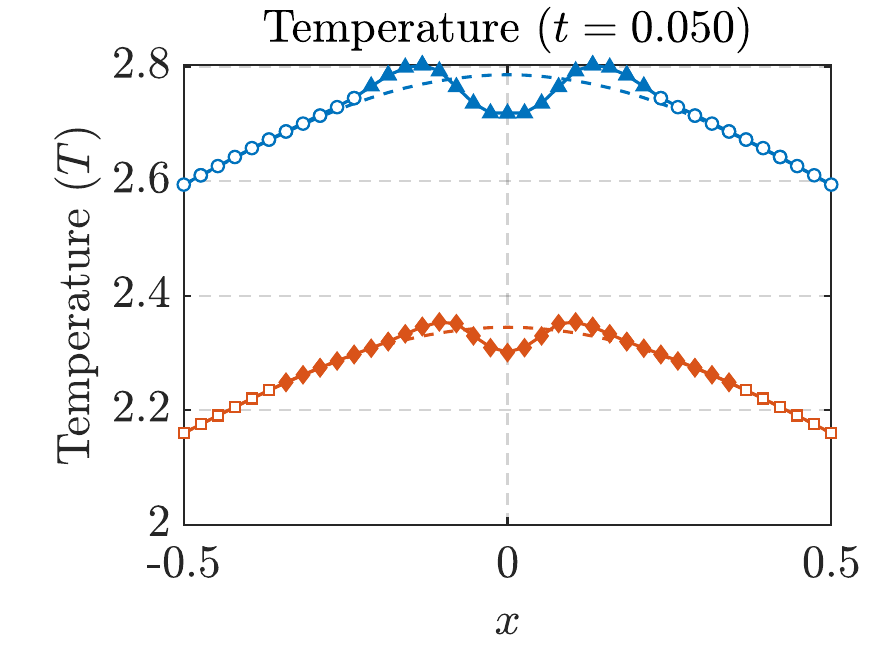}
\figrow{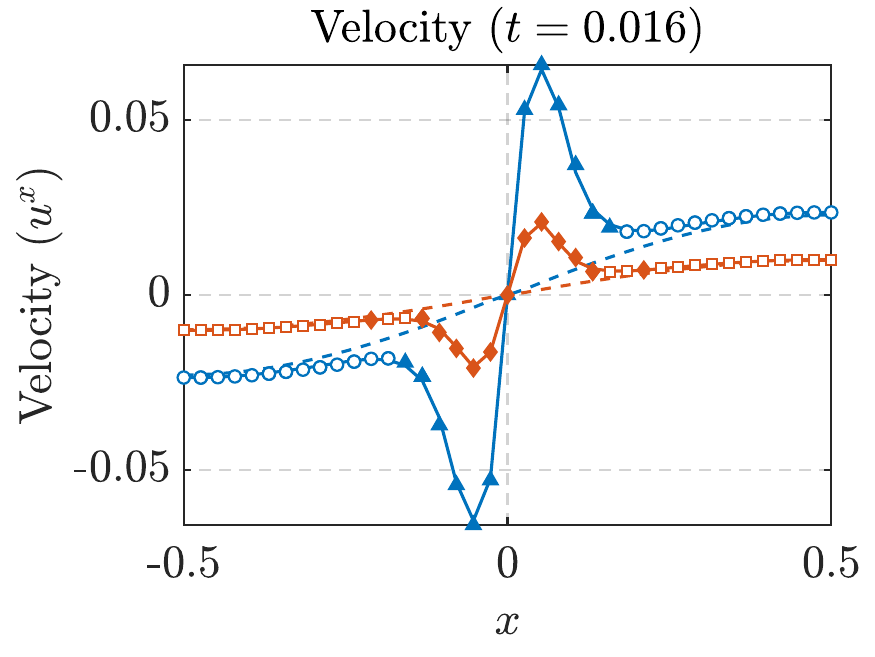}{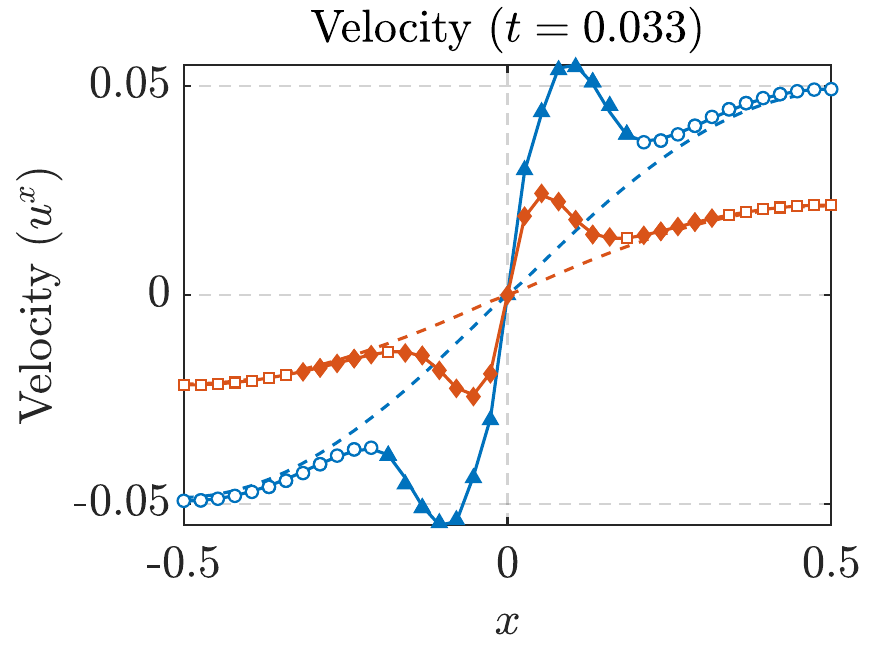}{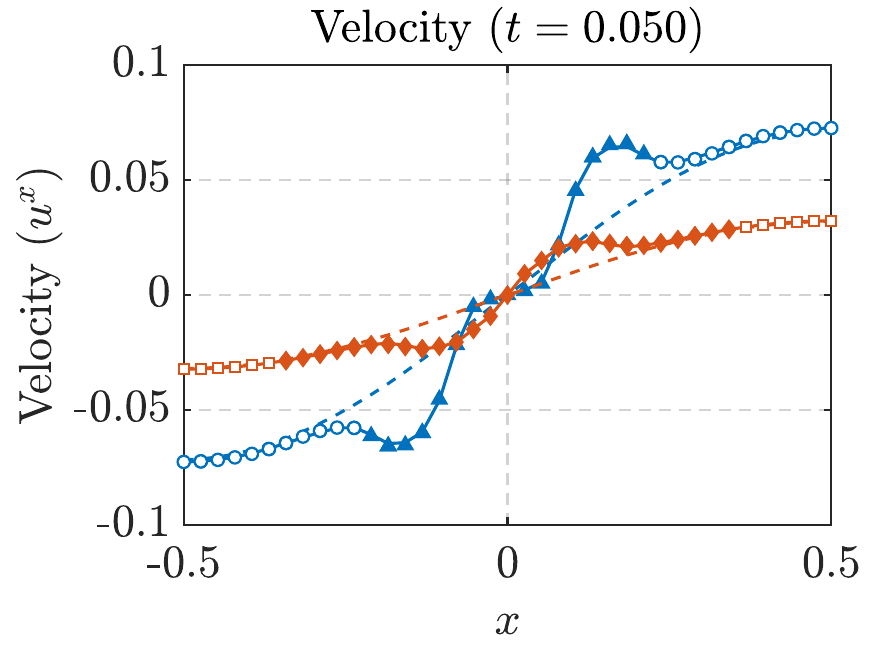}
\caption{\hyperref[Test4]{Test~4}. Mass density, temperature and velocity of both species along $x$ (projection along $y$) for the hybrid, Euler and kinetic solvers. For readability, the hybrid solution is shown only in a subset of cells, roughly equally spaced in $x$; open markers denote cells in the Euler regime and filled markers cells in the kinetic regime.}
\label{Fig15}
\end{figure}

Figure~\ref{Fig15} shows density, temperature and velocity of both species along $x$. The plots are restricted to $x\in[-0.5,0.5]$: outside this interval the gas is close to equilibrium and the three solutions cannot be distinguished. The kinetic solver is activated only in the central region. There the Euler solution differs clearly from the kinetic one, while the hybrid solution follows it closely.

The computational times are 0.9811~s (Euler), 111.8~s (kinetic) and 26.20~s (hybrid), a speedup of about 4.3. As in \hyperref[Test3]{Test~3}, we measure the distance from the kinetic solution with the space-time averaged $L^1$ error on the temperature of each species $i=1,2$,
\[
E(T_i)=\frac{1}{t_f\,|\Omega|}\sum_{n=1}^{N_t}\Delta t^n\sum_{c=1}^{N_P}|P_c|\,\big|T_{i,c}^n-T_{i,c}^{n,\mathrm{kin}}\big|,
\]
with the notation of \eqref{eq:errT}. For the Euler solver we obtain $E(T_1)=5.6\times10^{-3}$ and $E(T_2)=2.0\times10^{-3}$, and for the hybrid solver $E(T_1)=2.8\times10^{-4}$ and $E(T_2)=1.2\times10^{-4}$, about twenty and seventeen times smaller, respectively.

The switching thresholds are specified separately for each benchmark. The reported errors and speedups therefore are based on these choices.

\section{Computational cost}\label{secNumEff}
Table~\ref{table1} collects the computational times, the speedups and the errors $E(T_1)$ and $E(T_2)$ defined in \eqref{eq:errT} for Tests \hyperref[Test2]{2}, \hyperref[Test3]{3} and \hyperref[Test4]{4}. Errors are computed with respect to the kinetic solution. They are not reported for the shock tube: at $\varepsilon=10^{-5}$ all solvers converge to the same hydrodynamic solution, and the kinetic solution is not a better reference than the others.

In the rarefied tests the hybrid scheme reduces the temperature error of the Euler solver by about one order of magnitude, by a factor between 7 and 20, since the kinetic solver is used where non-equilibrium effects are significant. Kinetic computations for mixtures with different masses are expensive, because the velocity grid must resolve all species. In our tests the hybrid solver is between 3.4 and 4.6 times faster than the kinetic one. The gain depends on the number of active kinetic cell--species pairs, their reconstruction stencils and the overhead of the indicator and coupling.

\begin{table}[htbp]
\centering
\begin{tabular}{|cccc|cc|}
\toprule
\textbf{Test} & \textbf{Solver} & $E(T_1)$ & $E(T_2)$ & \textbf{Time} & \textbf{Time / hybrid time}\\
\midrule
\multirow{3}{*}{Shock tube} & Kinetic & -- & -- & 1403 s & 3.4\\
& Hybrid & -- & -- & 410.5 s & 1.0 \\
& Euler & -- & -- & 15.63 s & 0.038\\
\midrule
\multirow{3}{*}{NACA 0012} & Kinetic & -- & -- & 4475 s & 4.6 \\
& Hybrid & $7.5\times10^{-4}$ & $4.4\times10^{-4}$ & 981.7 s & 1.0\\
& Euler & $5.3\times10^{-3}$ & $5.2\times10^{-3}$ & 136.9 s & 0.14\\
\midrule
\multirow{3}{*}{Varying Kn} & Kinetic & -- & -- & 111.8 s & 4.3\\
& Hybrid & $2.8\times10^{-4}$ & $1.2\times10^{-4}$ & 26.20 s & 1.0\\
& Euler & $5.6\times10^{-3}$ & $2.0\times10^{-3}$ & 0.9811 s & 0.037\\
\bottomrule
\end{tabular}
\caption{Computational times, ratio between the time of each solver and that of the hybrid solver, and space-time averaged $L^1$ errors $E(T_1)$ and $E(T_2)$, defined in \eqref{eq:errT}, with respect to the kinetic solution, for Tests \hyperref[Test2]{2}, \hyperref[Test3]{3} and \hyperref[Test4]{4}. Errors are not reported for the shock tube, where all models converge to the same hydrodynamic solution.}
\label{table1}
\end{table}

\section{Conclusions}\label{secCon}
We have proposed a dynamic domain decomposition method for gas mixtures in which collisions between particles of the same species are dominant, so that the hydrodynamic limit has a velocity and a temperature for each species. The switching criterion extends the moment realizability approach of \cite{lev-mor-nad-1998}, used for a single gas in \cite{fil-rey-2014, fil-xio-2018, tiw-2000, li-son-wan-2021, xio-qiu-2017, cap-par-rey-2025, cap-ten-2025, cap-par-rey-2026}. For each species we derived a realizability matrix from the Chapman--Enskog expansion of the multi-species BGK model and used the comparison between its zeroth- and first-order approximations to decide locally between the multi-velocity, multi-temperature Euler equations and the kinetic equation. The indicator uses only macroscopic quantities.

Since the fluid model keeps a separate velocity and temperature for each species, the partition of the domain is computed separately for each species, and in the same cell some species may be described by the Euler equations and others by the kinetic model. This is useful when the species behave differently. The coupling between the two descriptions is done through moments and Maxwellians, and exchange terms between species in different regimes are computed from the moments.

The kinetic equation is solved with a second-order AP IMEX scheme on unstructured triangular meshes. We proved that, under suitable assumptions on the collision frequencies, the implicit relaxation step is uniquely solvable with positive temperatures for every time step, provided its input densities and temperatures are positive, and that in the limit $\varepsilon\to0$ the scheme reduces to a consistent discretization of the multi-velocity, multi-temperature Euler system. In the numerical tests the hybrid solution is close to the kinetic one, while the computational time is reduced by a factor between 3.4 and 4.6, on a standard laptop.

Several extensions are possible. In the companion paper \cite{PaperB} we consider the regime in which all collisions occur on the same time scale and the hydrodynamic limit has a single velocity and temperature; there the indicator has to be defined for the whole mixture, and the method is extended to third order. Further directions are the full Boltzmann operator for mixtures, polyatomic gases and chemical reactions, and parallel implementations for three-dimensional problems.

\section*{Acknowledgments}
The authors received funding from the European Union's Horizon Europe research and innovation program under the Marie Sk{\l}odowska-Curie Doctoral Network DataHyking (Grant No. 101072546).
DC and TR are also supported by the French government, through the UniCA$_{JEDI}$ Investments in the Future project managed by the National Research Agency (ANR) with the reference number ANR-15-IDEX-01.
The research of LP has been supported by the Royal Society under the Wolfson Fellowship ``Uncertainty quantification, data-driven simulations and learning of multiscale complex systems governed by PDEs''.  LP also acknowledges the  support by the Italian Ministry of University and Research (MUR) through the Fondo Italiano per la Scienza 2022–2023 (Bando FIS 2), under Decree No. 1236 of 01/08/2023, within the FIS2023-01334 Advanced Grant “Tackling Complexity: Advanced Numerical Approaches for Multiscale Systems with Uncertainties” (ADAMUS). This work has been written within the activities of GNCS group of INdAM (Italian National Institute of High Mathematics).

\section*{Declaration of generative AI and AI-assisted technologies}
During the preparation of this work the authors used a generative AI tool to improve the language and readability of the manuscript. After using this tool, the authors reviewed and edited the content as needed and take full responsibility for the content of the publication.

\bibliographystyle{acm}
\bibliography{biblio}
\end{document}